\documentclass{article}
\pdfoutput=1
\usepackage[utf8]{inputenc}
\usepackage{fancyhdr, hyperref, amsmath, amssymb, graphicx, wrapfig, amsthm, xcolor, bbold, enumitem, tcolorbox, wasysym, mathtools, tikz, scalerel, paracol, shortcuts}
\usepackage[margin=1in]{geometry}
\usepackage{fourier}
\usepackage{amsmath}
\usepackage{float}
\usepackage{anyfontsize}
\usepackage{mathscinet}
\usepackage{relsize}
\usepackage{subcaption}

\usepackage{titlesec}

\titleformat{\section}
  {\normalfont\Large\bfseries}
  {\thesection.}{0.2em}{\centering}

\title{A Well-Defined Jellyfish Algorithm for the Affine $E_7$ Subfactor Planar Algebra}
\author{Melody Molander}
\date{}
\begin{document}

\maketitle
\begin{abstract}
    In this paper, we contribute to the Kuperberg program by giving a diagrammatic presentation of generators and relations for the affine $E_7$ unshaded subfactor planar algebra. Using this presentation, we prove that its jellyfish algorithm is a well-defined surjection onto $\mbb{C}$. In particular, this shows that the jellyfish algorithm is an invariant on closed diagrams for this planar algebra.
\end{abstract}

\section{Introduction}

Subfactor planar algebras were introduced by Vaughan Jones \cite{Jon99} as a diagrammatic version of the standard invariant for subfactors. Subfactors are unital inclusions of factors, i.e., von Neumann algebras with trivial center. We will assume all subfactors are irreducible and all factors of type $II_1$. 

There are three main invariants of subfactors: a real number called the index, a principal graph, and the standard invariant (a unitary 2-category). In a landmark paper, Jones \cite{Jon83} proved that the set of indices for subfactors equals $\{4\cos^2(\pi/n):n\geq 3\}\cup [4,\infty]$. Subfactor classification has been the work of many authors and is completed up to index $5.25$. See \cite{JMS14, AMP23} for the history of this classification program. The standard invariant can be axiomatizated in many ways other than as a 2-category or a subfactor planar algebra. This includes Popa's $\lambda$-lattices \cite{Pop95} and Ocneanu's paragroups for finite depth subfactors \cite{Ocn88, EK98}. 

A planar algebra is a collection of vector spaces along with an action of the planar operad. When the planar algebra is evaluable, spherical, and positive definite, it is said to be ``subfactor". Given a subfactor planar algebra, Popa \cite{Pop95} proved, using the language of $\lambda$-lattices, that there is subfactor whose standard invariant is the planar algebra. On the other hand, given a subfactor with finite index, Jones proved that its standard invariant forms a subfactor planar algebra. The planar algebra language was found to be incredibly useful and pushed the subfactor classification program forward \cite{BPMS12, Han10, Pet10}. In particular, \cite{BPMS12} proved that the extended Haagerup graph is realizable as a principal graph of a subfactor.

Constructions of subfactor planar algebras are typically made using generators and relations presentations \cite{Big10, Mol24, MPS10, Pet10}. When a planar algebra is a standard invariant of a subfactor, its diagrams have checkerboard shading. Dropping the references to shading in the definition of planar algebra yields an \textbf{unshaded subfactor planar algebra}, which is still a mathematically rich object. From presentations of unshaded subfactor planar algebras, one can obtain diagrammatics for fusion categories \cite{Mol24} and construct link invariants \cite{MPS11}. Morrison, Peters, and Snyder \cite{MPS10} posed the following program to construct and explore diagrammatic presentations:
\\~

\noindent \tbf{The Kuperberg Program:} Provide a generators-and-relations presentation for every interesting planar algebra. Using this presentation, understand other properties of the associated categories, the planar algebra, and, in the subfactor case, the standard invariant. 
\\~

Much progress has been made towards this program. The Kuperberg program for index less than 4 has been completed \cite{MPS10, Big10}. For index greater than 4, this has been ongoing research \cite{Pet10, BPMS12, Han10, MP15}. In this paper, we focus on giving a presentation by generators and relations of the unshaded index 4 subfactor planar algebra with principal graph affine $E_7$. 

Index 4 subfactor planar algebras were classified by Popa \cite{Pop94}. Their principal graphs are necessarily the bipartite affine $ADE$ Dynkin diagrams or the $A_\infty$ Dynkin diagram. In the unshaded subfactor planar algebra case, while these do not give standard invariants of subfactors, there are further index 4 principal graphs beyond those found by Popa. For example, this author \cite{Mol24} proved that unshaded subfactor planar algebras can arise from all affine $A$ Dynkin diagrams (not just the bipartite ones). The Kuperberg program at index 4 for the principal graph affine $A$ was completed \cite{Mol24} and her thesis \cite{Mol25} includes the affine $D$ case. The $A_\infty$ subfactor planar algebras are the ubiquitous Temperley-Lieb planar algebras. 

When describing a planar algebra using diagrammatic generators and relations, the $k$th vector space can be described as linear combinations of planar diagrams using the generators and strands with $k$ strands on both the top and bottom boundary, modulo the skein relations. In order to show a planar algebra is ``subfactor", often the toughest criteria to check is that its $0$th vector space of closed diagrams is one-dimensional. Bigelow, Morrison, Peters, and Snyder's paper \cite{BPMS12} introduced a novel evaluation algorithm called \textbf{the jellyfish algorithm} to show that this vector space was at most one-dimensional for the extended Haagerup subfactor planar algebra. Their argument also proved that this algorithm works for other planar algebras with sufficiently nice relations. Bigelow and Penneys \cite{BP14} then gave criteria for which subfactor planar algebras of index greater than 4 can be constructed using the jellyfish algorithm. 

The jellyfish algorithm was an extremely clever evaluation algorithm. In quantum topology, typically evaluation algorithms reduce the ``common sense" complexity in each step. Examples include the Kauffman bracket which reduces the number of crossings of a link in every step, and the HOMFLYPT polynomial which reduces the number of crossings or the unknotting number in each step. However, the jellyfish algorithm may initially \textit{increase} complexity! Specifically, the first steps may drastically increase the number of generators in a diagram. Despite being counterintuitive, this initial increase in complexity puts the elements of the vector space in a nice form that can then be simplified through a graph theory argument.

In this paper, we show that the jellyfish algorithm can be used to bound this vector space dimension from below as well for the affine $E_7$ unshaded subfactor planar algebra. Thus, the jellyfish algorithm is an invariant on the set of closed diagrams for this planar algebra. In particular, we show that the algorithm gives a well-defined surjection onto $\mbb{C}$, despite the numerous choices made throughout the algorithm's process. Therefore, the vector space is at least one-dimensional. 

In \cite{BPMS12}, the authors show the vector space of closed diagrams is at least one-dimensional by locating it inside the graph planar algebra of the extended Haagerup principal graph. They first use numerical approximation methods to find the generator of the planar algebra inside of a 148375-dimensional vector space. Then to show that this generator satisfied the required relations, they use computer-assisted proofs using Mathematica. 

Instead of embedding the planar algebra into the graph planar algebra, one can prove this vector space is at least one-dimensional by showing that the jellyfish algorithm gives a well-defined surjection. The relations of the affine $E_7$ subfactor planar algebra allow generators to pass through strands as in a braided tensor category, which is not a relation afforded to the extended Haagerup planar algebra. Consequently, the complexity of the generators does not initially increase in the algorithm, but the number of crossings does. Using the braided skein theory, we can avoid any computer-assisted proofs that would have been used by embedding into the graph planar algebra. We hope these techniques can be generalized to get invariants  of closed diagrams for other subfactor planar algebras or tensor categories with sufficiently nice braided structure.

In section \ref{Background}, we first review the necessary background material, starting with planar algebras, then the principal graph, and concluding with an overview of the jellyfish algorithm. Then we prove results that hold for any index 4 unshaded subfactor planar algebra in section \ref{generalpafacts} that will be heavily used in the subsequent sections. Section \ref{sufficientE7} first defines a planar algebra using generators and relations then shows that these relations are sufficient to define the affine $E_7$ unshaded subfactor planar algebra (Theorem \ref{thm: E7sufficientmaintheorem}). To do so, we define the jellyfish algorithm in subsection \ref{e7jellyfish} and show it is a well-defined surjective function onto $\mbb{C}$. The proof of this is substantial and is carried out across many subsections. The section then finishes with describing the minimal projections in the planar algebra. The final section \ref{necessaryE7} shows that if one starts with the unshaded subfactor planar algebra of affine $E_7$, it has a presentation that is necessarily the same as the one found in the previous section.

\subsection{Acknowledgments}
The author thanks Stephen Bigelow, David Penneys, and Gregory Faurot for their helpful conversations and insight. The author would also like to thank the careful referee who found an error in the preprint as well as gave many useful suggestions, including guidance on the proof of Lemma \ref{lem: there-is-an-S2}.

\section{Background}\label{Background}

\subsection{Planar Algebras}
\begin{figure}[h!]
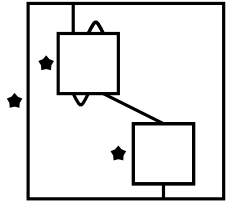

    \centering
   $\unshadedtangleexample$
    \caption{An example of a planar tangle.}
    \label{fig:unshadedtangleexample}
\end{figure}

We first give an overview of planar algebras. The interested reader can find more background \cite{Mol24, Mol25, MPS10, Pet10}. Let $X$ denote a black, unoriented strand. A \tbf{planar tangle} is a diagram of a marked square with a (possibly empty) finite collection of marked squares interior to it connected by nonintersecting strands, where the number of strands on the top and bottom of each square must be equal (see figure \ref{fig:unshadedtangleexample}). Label the interior squares $D_k$ and the exterior square $\mc{D}$. Let $\partial{D_k}$ (respectively $\partial \mc{D}$) equal the number of strands ending on the top boundary of $D_k$ (respectively $\mc{D}$). 

\begin{figure}[h!]
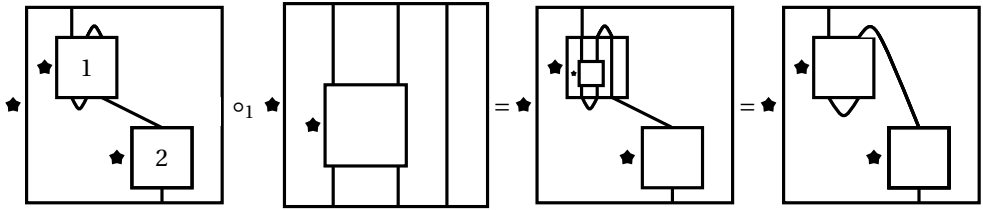

    \centering
       \scalebox{1}{
$\adjustbox{scale=0.94}{\unshadedtangleexamplelabelled} \text{ }\circ_1 \text{ }\adjustbox{scale=0.91}{\unshadedtanglecompositiona}=\adjustbox{scale=0.91}{\unshadedtangleexamplewithcompositiona}=\adjustbox{scale=0.91}{\unshadedtanglecompositionb}$
}
    \caption{Composition of planar tangles.}
    \label{fig:compositionofplanartangles}
\end{figure}

Let $T_1$ and $T_2$ be two planar tangles. Suppose for some interior square, $D_k^1$ of $T_1$,  $\partial D_k^1$ equals $\partial \mc{D}^2$, the boundary of the exterior square of $T_2$. Then we can compose $T_1$ and $T_2$ by inserting $T_2$ into the square $D_k^1$. See figure \ref{fig:compositionofplanartangles} for an example. 

Let $\{\mc{P}_k\}_{k\in \mathbb{Z}_{\geq 0}}$ be a collection of vector spaces over $\mathbb{C}$. For a square in a planar tangle with $k$ boundary points on the top and bottom, associate it to the vector space $\mc{P}_k$. In figure \ref{fig:unshadedtangleexample}, the bottom interior square and the outside square are associated to $\mc{P}_1$ and the top interior square is associated to $\mc{P}_3$. To every planar tangle with interior squares $\{D_k\}$ and exterior square $\mc{D}$, we associate a multilinear map $Z_T:\otimes_k \mc{P}_{\partial D_k}\to \mc{P}_{\partial \mc{D}}$.

\begin{figure}[h!]
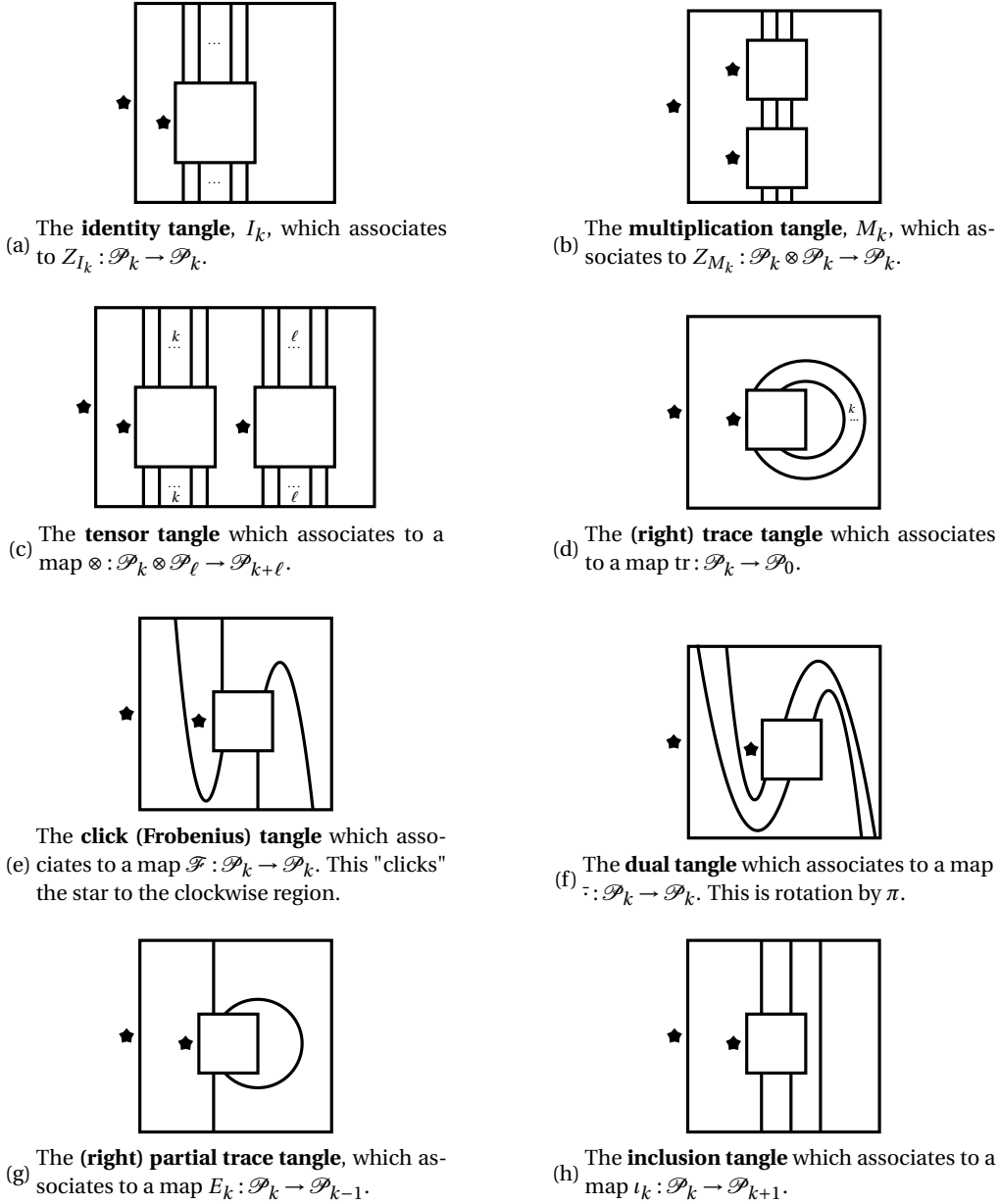

    \centering
    \begin{subfigure}{.45\textwidth}
    \centering
    \adjustbox{scale=0.97}{\unshadeidentitytangle}
    \caption{\parbox{0.75\textwidth}{The \textbf{identity tangle}, $I_k$, which associates to $Z_{I_k}:\mc{P}_k\to\mc{P}_k$.}}
    \label{fig:identitytangle}
\end{subfigure}%
\begin{subfigure}{.45\textwidth}
  \centering
 \adjustbox{scale=0.97}{\unshademultiplicationtangle}
  \caption{\parbox{0.75\textwidth}{The \tbf{multiplication tangle}, $M_k$, which associates to $Z_{M_k}:\mc{P}_k\otimes \mc{P}_k\to \mc{P}_k$.}}
  \label{fig: multiplication-tangle}
\end{subfigure}%
\\~\\

\begin{subfigure}{.45\textwidth}
\centering
   \adjustbox{scale=0.97}{ \unshadetensortangle}
  \caption{\parbox{0.74\textwidth}{The \tbf{tensor tangle} which associates to a map $\otimes: \mc{P}_k\otimes \mc{P}_\ell\to \mc{P}_{k+\ell}$.}}
  \label{fig: tensor-tangle}
\end{subfigure}% 
\vspace{0.1cm}\begin{subfigure}{.45\textwidth}
  \centering
   \adjustbox{scale=0.97}{ \unshadetracetangle}
  \caption{\parbox{0.75\textwidth}{The \tbf{(right) trace tangle} which associates to a map $\text{tr}:\mc{P}_k\to \mc{P}_{0}$.}}
  \label{fig: trace-tangle}
\end{subfigure}%
\\~\\

\begin{subfigure}{.45\textwidth}
\centering
   \adjustbox{scale=0.97}{ \unshadeclicktangle}
  \caption{\parbox{0.75\textwidth}{The \tbf{click (Frobenius) tangle} which associates to a map $\mc{F}:\mc{P}_k\to \mc{P}_k$. This "clicks" the star to the clockwise region.}}
  \label{fig: click-tangle}
\end{subfigure}% 
\begin{subfigure}{.45\textwidth}
  \centering
    \adjustbox{scale=0.97}{\unshadedualtangle}
  \caption{\parbox{0.75\textwidth}{The \tbf{dual tangle} which associates to a map $\overline{\cdot}:\mc{P}_k\to \mc{P}_k$. This is rotation by $\pi$.}}
  \label{fig: dual-tangle}
\end{subfigure}%
\\~\\

\begin{subfigure}{.45\textwidth}
\centering
  \adjustbox{scale=0.97}{  \unshadepartialtracetangle}
  \caption{\parbox{0.75\textwidth}{The \tbf{(right) partial trace tangle}, which associates to a map $E_k:\mc{P}_k\to \mc{P}_{k-1}$.}}
  \label{fig: click-tangle}
\end{subfigure}% 
\begin{subfigure}{.45\textwidth}
  \centering
   \adjustbox{scale=0.97}{ \unshadeinclusiontangle}
  \caption{\parbox{0.75\textwidth}{The \tbf{inclusion tangle} which associates to a map $\iota_k: \mc{P}_k\to \mc{P}_{k+1}$.}}
  \label{fig: dual-tangle}
\end{subfigure}%
\caption{A collection of common planar tangles.}
\label{fig: planar-tangles-examples}
\end{figure}

 We will use several common tangles in this paper, which are given in figure \ref{fig: planar-tangles-examples}. Many of these maps have subscripts dropped. In those cases, it should be clear from context the vector spaces in the domain and codomain. For both the right trace tangle and right partial trace tangle, there is an analogous left version, denoted $\te{tr}^\ell$ and $E_k^\ell$, respectively.

Planar algebras were first defined by Jones in \cite{Jon99}. A \tbf{planar algebra} is a collection of vector spaces over $\mbb{C}$, $\{\mc{P}_0, \mc{P}_1, \mc{P}_2,...\}$, along with for every planar tangle $T$,  an associated a multilinear map $Z_T: \otimes_k \mc{P}_{\partial D_k}\to \mc{P}_{\partial \mc{D}}$ such that
\begin{enumerate}[label=(\roman*)]
    \item planar isotopic tangles give the same map;
    \item composition of tangles corresponds to the composition of their maps;
    \item for every $k\geq 0$, the identity tangle is the identity on $\mc{P}_k$.
\end{enumerate}

We frequently call the $k$th vector space $\mc{P}_k$ of a planar algebra $\mc{P}$ its \tbf{$k$th box space}.

\begin{ex}
    As a first example, we define the \tbf{Temperley-Lieb planar algebra}, $\mc{TL}$. Choose the vector spaces to be the Temperley-Lieb algebras $\{\mc{TL}_k\}_{k \in \mbb{Z}_{\geq 0}}$ over $\mbb{C}$, with a chosen value $\delta>0$. The Temperley-Lieb algebras \cite{TL71} were defined diagrammatically by Kauffman \cite{Kau87}, which we describe here. Elements of $\mc{TL}_k$ are linear combinations of diagrams with $k$ points on the top and bottom with nonintersecting strands connecting those $k$ points. See figure \ref{fig: elements-of-TL3} for examples of elements in $\mc{TL}_3$. Multiplication, $fg$, of two elements $f$ and $g$ in $\mc{TL}_k$ is given by stacking $f$ on top of $g$. When a ``bubble" is formed, it can be removed for a factor of $\delta$ (see figure \ref{fig: multiplication-in-TL}). As an algebra, $\mc{TL}_k$ is generated by elements $X^{\otimes k}$ and $\{e_i\}_{1\leq i \leq k-1}$, where $e_i$ is the diagram with straight strands except a cup and cap in the $(i, i+1)$ position. 
    
    Now the associated maps to the tangles can be defined through insertion. For example, the map $Z_{M_k}:\mc{TL}_k\to \mc{TL}_k$ associated to the multiplication tangle is the map given by inserting $f\in \mc{TL}_k$ in to the top interior square and $g \in \mc{TL}_k$ into the bottom interior square. The resulting diagram is $fg$. See figure \ref{fig: multiplication-tangle-TL}. 
\end{ex}

\begin{figure}[h!]
    \centering
    \begin{overpic}[unit=1mm, scale=0.4]{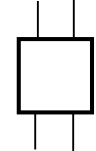}
        \put(4, 6.2){$R$}
    \end{overpic} 
    \caption{An example of a potential generator of a planar algebra. This generator would be in $\mc{P}_2$.}
    \label{fig: generator-of-arbPA}
\end{figure}

\begin{figure}[h!]
    \centering
        \begin{overpic}[unit=1mm, scale=0.4]{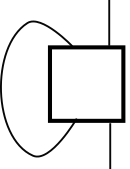}
        \put(8, 7.5){$R$}
    \end{overpic}, 
 \begin{overpic}[unit=1mm, scale=0.4]{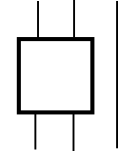}
        \put(4.3, 6.2){$R$}
    \end{overpic}
    \caption{If a planar algebra $\mc{P}$, is generated by $R$ in figure \ref{fig: generator-of-arbPA}, then these would be elements in $\mc{P}_1$ and $\mc{P}_3$, respectively.}
    \label{fig: elements-in-arbPA}
\end{figure}

The preceding example shows that to best understand the associated multilinear maps to planar tangles, the vector spaces should be described using diagrams. Then the map is just the diagram obtained by inserting diagrams into the interior squares. Therefore, it is best to describe each vector space using diagrammatic generators and skein relations. The \tbf{generators} of a planar algebra $\mc{P}=\{\mc{P}_k\}_{k\geq 0}$ are a collection of diagrams that can be combined, along with strands, in any arbitrary planar fashion, to obtain any element of its vector spaces. Then any element of $\mc{P}_k$ is a linear combination of diagrams that are made up of generators and nonintersecting strands with $k$ boundary points on the top and bottom of the diagrams. As an example, if a planar algebra is generated by some diagram $R$ in $\mc{P}_2$ subject to some skein relations (see figure \ref{fig: generator-of-arbPA}), then figure \ref{fig: elements-in-arbPA} gives an element of $\mc{P}_1$ and $\mc{P}_3$, respectively.

 \begin{figure}[h]
        \centering
        $f=$  \adjustbox{scale=0.91}{\unshadefmorph} \quad $g=$  \adjustbox{scale=0.91}{\unshadegmorph}
        \caption{Two elements of $\mc{TL}_3$.}
        \label{fig: elements-of-TL3}
    \end{figure}

    \begin{figure}[h]
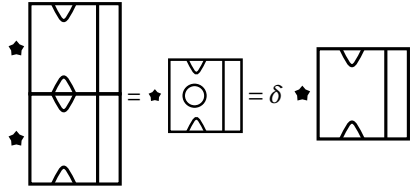

        \centering
         \adjustbox{scale=0.91}{\unshadegsquared}=\hspace{0.1cm} \adjustbox{scale=0.91}{\unshadegbubble} $= \delta$ \te{ } \adjustbox{scale=0.91}{\unshadegmorph}
        \caption{Multiplication of $g^2$ (from figure \ref{fig: elements-of-TL3}) in $\mc{TL}_3$.}
        \label{fig: multiplication-in-TL}
    \end{figure}

    \begin{figure}[h!]
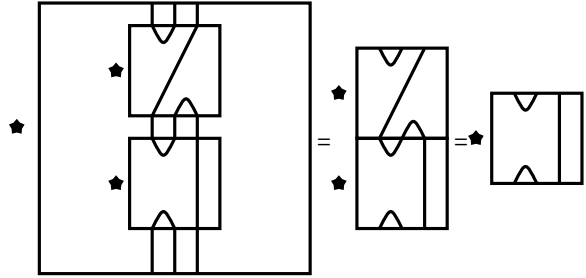

        \centering
        \adjustbox{scale=0.91}{\unshadefgmorphtangle}=\adjustbox{scale=0.91}{\unshadefgmorph}=\adjustbox{scale=0.91}{\unshadegmorph}
        \caption{The multiplication tangle map in the Temperley-Lieb planar algebra with $f$ and $g$ of figure \ref{fig: elements-of-TL3} plugged in.}
        \label{fig: multiplication-tangle-TL}
    \end{figure}

\begin{rem*}
    The Temperley-Lieb planar algebra is said to have no generators, since every element of every vector space consists of linear combinations of diagrams with only nonintersecting strands.
\end{rem*} 

We define a diagram $D_k\in \mc{P}_k$ to be \tbf{uncuppable} (respectively, \tbf{uncappable}) if for all $1\leq j \leq k-1$, $e_j \in \mc{P}_k$, we have that $D_ke_j=0$ (respectively, $e_j D_k=0$). We define a diagram $D_k\in \mc{P}_k$ to be \tbf{unsidecappable} if the left and right partial trace of $D_k$ are both zero. If $D_k\in \mc{P}_k$ is uncuppable, uncappable, and unsidecappable, we say $D_k$ is \tbf{U.C.C.S.}

For this paper, our planar algebras will all be unshaded. However, in order to be standard invariants of subfactors, they would also need a checkerboard shading, see \cite{BPMS12, Jon99, Mol24, Mol25, Pet10} for examples of planar algebras with checkerboard shading and their connections to subfactors. Beyond a checkerboard shading, subfactor planar algebras $\mc{P}=\{\mc{P}_k\}_{k\geq 0}$ need to satisfy the below properties.
\begin{enumerate}
    \item $\mc{P}$ is \textbf{evaluable}: $\mc{P}_0$ is one-dimensional and for all $k> 0$, $\mc{P}_k$ is finite-dimensional. 
    \item $\mc{P}$ is \tbf{spherical}: the left and right trace map on $\mc{P}_1$ are equal. 
    \item $\mc{P}$ is \tbf{positive definite}: there exists an antilinear adjoint map on each vector space, denoted $*:\mc{P}_k \to \mc{P}_k$, that is compatible with reflection over a horizontal line on planar tangles. Further, the sesquilinear form on each $\mc{P}_k$ defined by $\langle f,g \rangle \coloneq \te{tr}(fg^*)$ is positive definite. 
\end{enumerate}

An (unshaded) planar algebra that is evaluable, spherical, and positive definite is said to be an (unshaded) \tbf{subfactor planar algebra.} We will commonly drop ``unshaded" and call both shaded and unshaded cases ``subfactor planar algebras".  

\begin{ex}
    The Temperley-Lieb planar algebra with chosen $\delta>0$ with $\delta^2\geq 4$ is an unshaded subfactor planar algebra.
\end{ex}

Jones \cite{JonPA} proved that any subfactor of index 4 or greater satisfies that its subfactor planar algebra contains the Temperley-Lieb planar algebra where $\delta^2$ is the index of the subfactor. We then refer to $\delta^2$ as the \tbf{index} of the planar algebra.

\subsection{The Principal Graph}

We assume some familiarity with tensor categories. See \cite{EGNO15}, Chapter 2 and 4 for more details on this topic. From a subfactor planar algebra, we can obtain a tensor category. Let $\mc{P}=\{\mc{P}_k\}_{k\geq 0}$ be a subfactor planar algebra. A \tbf{projection} in $\mc{P}$ is an element $P$ in the planar algebra such that $P^2=P=P^*$. A \tbf{morphism} $f\in \te{Hom}(P,Q)$ between two planar algebra elements $P$ and $Q$ is a diagram whose bottom boundary is $Q$ and top boundary is $P$. We say two projections $P$ and $Q$ are \tbf{isomorphic} if and only if there exists a morphism $f\in \te{Hom}(P,Q)$ such that $ff^*=P$ and $f^*f=Q$. 

\begin{ex}
    The strand $X$ is a projection in $\mc{TL}$. Further, $\mc{TL}_k$ is a collection of morphisms in $\te{End}(X^{\otimes k})$. 
\end{ex}

\begin{rem*}
    For any subfactor planar algebra, $\mc{P}$, $\mc{P}_k=\te{End}(X^{\otimes k})$. 
\end{rem*}

First create a rigid tensor category from $\mc{P}$ whose objects are projections and morphisms are morphisms in $\mc{P}$. In this category, take $\otimes$ to be horizontal concatenation, with unit object $\emptyset$, and dual (which we denote by $\overline{\cdot}$) to be rotation by $-\pi$. Then take the abelian envelope of that category and get the category $\mc{C}_\mc{P}$. That is, objects are direct sums of projections and 
\begin{align*}
    \text{Hom}(P_1 \oplus ... \oplus P_m, Q_1 \oplus ... \oplus Q_n)=\oplus_{i,j} \text{Hom}(P_i, Q_j).
\end{align*}

A projection $P$ in $\mc{P}$ is said to be \tbf{minimal} if $\te{End}(P)$ is one-dimensional. In the Temperley-Lieb planar algebra, the minimal projections are called the \textbf{Jones-Wenzl projections}. For this planar algebra, there is exactly one minimal projection in each box space. We denote $f^{(k)}$ to be the minimal projection in $\mc{TL}_k$. Jones-Wenzl projections are uncuppable, uncappable, and satisfy \tbf{Wenzl's relation} \cite{Wen87}: $f^{(0)}=\emptyset$, $f^{(1)}=X$, and
\begin{align*}
    f^{(k)}\otimes X \cong f^{(k-1)}\otimes f^{(k+1)}.
\end{align*}

Since any subfactor planar algebra of index 4 or greater contains Temperley-Lieb, the planar algebra of this paper (which is of index 4) will contain these projections (but they may no longer be minimal). The trace of Jones-Wenzl projections is well-known. First, we define the $q$-deformation of integers. The \tbf{$n$th quantum number} is:
    \begin{align*}
        [n]_q = \frac{q^n-a^{-n}}{q-q^{-1}}=q^{n-1}+q^{n-3}+...+q^{-(n-3)}+q^{-(n-1)}.
    \end{align*}
We can always define $q$ such that $\delta=[2]_q=q+q^{-1}$. When clear, we drop the subscript $q$. The Jones-Wenzl projections satisfy that $\te{tr}(f^{(k)})=[k+1]_q$. In particular, $\te{tr}(f^{(0)})=\te{tr}(\emptyset)=1$, $\te{tr}(f^{(1)})=\te{tr}(X)=\delta$. 

For any subfactor planar algebra $\mc{P}$, the category $\mc{C}_{\mc{P}}$ is semisimple, i.e., the tensor product of any two objects can be decomposed as a direct sum of minimal projections. How objects decompose can be reflected in a graph called the \tbf{principal graph}. We describe a construction of the principal graph. This is not the original construction, but Bisch \cite{Bis97} showed this yields the same graph. Vertices of a principal graph are isomorphism classes of minimal projections. The number of edges between two vertices $P$ and $Q$ is the number of times $Q$ appears in the direct sum decomposition: 
\begin{equation}\label{eq: principal-graph}
    P\otimes X \cong \oplus_{Q\in N(P)} Q,
\end{equation}
where $N(P)$ is the set of neighbors of $P$ in the principal graph. 

\begin{ex}
    From Wenzl's relation, we can see that the principal graph of the Temperley-Lieb planar algebra is the $A_\infty$ Dynkin diagram.
\end{ex}

\begin{rem*}
    In the language of category theory, minimal projections are the simple objects of $\mc{C}_\mc{P}$ and the principal graph is its fusion graph with respect to $X$.
\end{rem*}

By taking the trace of both sides of (\ref{eq: principal-graph}), one finds that the principal graph of a subfactor planar algebra encodes its index. That is, 
    \begin{align*}
        \delta \te{tr}(P)=\sum_{Q\in N(P)} \te{tr}(Q),
    \end{align*}
which is called \tbf{the trace formula}. Additionally, we have the following useful fact.

\begin{prop} \label{prop: JonesdirectsumoftensorsofX} (Jones \cite{JonPA}, Theorem 3.1.3) Let $\mc{P}$ be a subfactor planar algebra. Then $X^{\otimes k}\cong \oplus_{Q\in T_k(\emptyset)}Q$, where $T_k(\emptyset)$ is the multiset of endpoints length $k$ from $\emptyset$, including multiples of the same element. 
\end{prop}

This paper will solely be focused on an index 4 unshaded subfactor planar algebra. Index 4 (shaded) subfactor planar algebras have been classified by Popa \cite{Pop94}. Popa found that principal graphs of index 4 (shaded) subfactor planar algebras are the $A_\infty$  Dynkin diagram (when the planar algebra is Temperley-Lieb) and the bipartite affine $ADE$ Dynkin diagrams. He also shows that when the principal graph is affine $E_7$ (which we sometimes will denote $\til{E}_7$), there is exactly 1 (shaded) subfactor planar algebra, up to planar algebra isomorphism.

\subsection{The Jellyfish Algorithm}
When defining a planar algebra by generators and relations, showing that $\mc{P}_0$ is one-dimensional is commonly the toughest criteria to check. For the planar algebra in this paper, we will show it is at most one-dimensional by showing there is an evaluation algorithm. In particular, we will describe a \tbf{jellyfish algorithm}, as originally defined in \cite{BPMS12}.

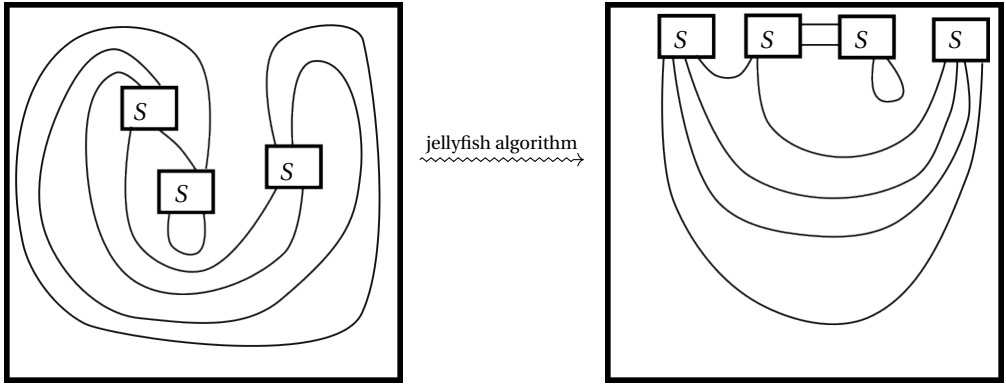
\begin{figure}[h!]
    \centering
         \begin{eqnarray*}
        \adjustbox{scale=1.1}{\begin{overpic}[unit=1mm, scale=0.5]{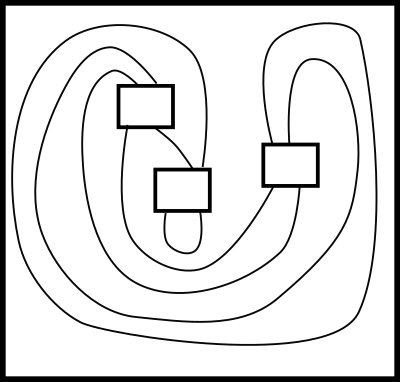}
     \put(22.5,24){$S$}
     \put(18,35){$S$}
     \put(37,27.5){$S$}
     \end{overpic}}
        \quad
    {\makebox[0pt][l]{\raisebox{2.8cm}{$\xrsquigarrow{\text{jellyfish algorithm}}$}}}\quad \quad\quad \quad\quad \quad\quad \quad
     \adjustbox{scale=1.1}{\begin{overpic}[unit=1mm, scale=0.5]{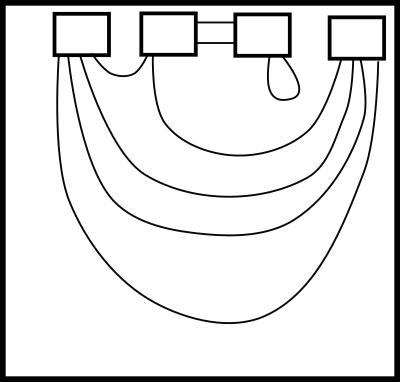}
     \put(21,44.5){$S$}
     \put(9.5,44.5){$S$}
     \put(33.5,44.5){$S$}
     \put(46,44.5){$S$}
     \end{overpic}}
     \end{eqnarray*}
    \caption{A schematic illustration of the first steps of the jellyfish algorithm.}
    \label{fig: general-jfa}
\end{figure}

Previously, in \cite{BPMS12}, the jellyfish algorithm was applied to a planar algebra with a single generator, say $S\in \mc{P}_k$. (Note that this paper was handling checkerboard shaded subfactor planar algebras, so had $S$ in both shadings, $S_{\pm}\in \mc{P}_{k,\pm}$). We outline the process in this case. We commonly refer to a generator such as in figure \ref{fig: generator-of-arbPA} as a ``box". Let $D\in \mc{P}_0$. The idea of the jellyfish algorithm is to recognize $D$ as being a diagram in a square ``tank". When you can pass a box past a strand as well as click strands around boxes, you can float the boxes to the top of the tank with their strands hanging down (resembling a jellyfish). When clicking the strands around a box, you may need to multiply by a scalar. When moving a box past a strand, you may get a sum of diagrams with more boxes that are moved past the strand and below the boxes is some picture in Temperley-Lieb (we are soon to define ``trains"  and then can restate this to mean any box with a strand over it can be rewritten as a linear combination of trains). We allow drawing strands all on the bottom of the box or if the box is attached to its neighbor, on the right or left sides. It is implied that these strands rainbow over the box on the right. We call a diagram at this point in \tbf{jellyfish form}. Below, in figure \ref{fig: general-jfa}, is an example where the diagram on the left is floated to the top with the right-hand side picture being one of its resulting summands in jellyfish form. 

 Now, since strands cannot cross, a graph theoretical argument given in \cite{BPMS12} (Theorem 3.8) gives that either an $S$ is cupped or there exists at least one pair of adjacent boxes that are connected by at least $k$ strands. That is, the two $S$-boxes form an $S^2$. If there is no cup, and it is true that $S^2$ is a linear combination of diagrams with at least one less $S$, then eventually we can reduce the diagram to a linear combination of diagrams with zero or one $S$-box. In the summands where there are no $S$-boxes, the diagram is closed and in Temperley-Lieb, so can be evaluated to a number by popping bubbles. In the summands where there is exactly one $S$-box, $S$ must be cupped, capped, or side-capped somewhere since the diagram is closed. If $S$ is known to be U.C.C.S, then these summands are zero. This would then give that $D$ can be evaluated. 

Below, in Proposition \ref{prop: BMPS}, we will list the relations that would be sufficient to have a jellyfish algorithm on a single generator. The first two relations allow the boxes to be floated to the top of the tank, and the latter two relations allow the diagram in jellyfish form to be evaluated. Before stating this propostion, however, we give a few definitions.

Recall that in a rigid monoidal category, every object $a$, has a dual $\overline{a}$, with \textbf{evaluation} and \textbf{coevaluation} maps $\te{ev}_a:\overline{a}\otimes a \to \emptyset$, $\te{coev}_a:\emptyset \to a\otimes \overline{a}$ satisfying "zig-zag" identities \cite{EGNO15} (Definition 2.10.1). Diagrammatically, $\te{ev}_X$ is a cap and $\te{coev}_X$ is a cup. 

\begin{figure}[h!]
\centering
    \begin{subfigure}{.33\textwidth}
    \centering
    \begin{align*}
       \adjustbox{scale=0.91}{\begin{overpic}[unit=1mm, scale=0.5]{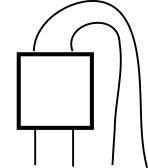}
     \put(5,9){$D_k$}
     \put(5,3){$...$}
     \put(15.4,3){$...$}
     \end{overpic}}
\end{align*}
    \caption{The ``rainbowed" diagram $\mc{R}(D_k)$.}
    \label{fig: rainbowed-box}
\end{subfigure}%
\begin{subfigure}{.33\textwidth}
  \centering
    \begin{align*}
    \adjustbox{scale=0.91}{     \begin{overpic}[unit=1mm, scale=0.5]{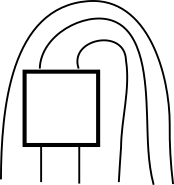}
     \put(5.5,9){$D_k$}
     \put(6,3){$...$}
     \put(16,3){$...$}
     \end{overpic} }
    \end{align*}
     \caption{\parbox{0.75\textwidth}{The ``over the rainbow", diagram $\mc{O}(\mc{R}(D_k))$.}}
    \label{fig: over-the-rainbow}
\end{subfigure}
\begin{subfigure}{.33\textwidth}
  \centering
    \begin{align*}
     \adjustbox{scale=0.91}{\begin{overpic}[unit=1mm, scale=0.5]{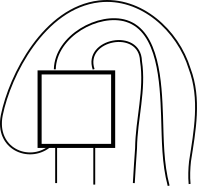}
     \put(7.5,9){$D_k$}
     \put(8,3){$...$}
     \put(18,3){$...$}
     \end{overpic}}
\end{align*}
     \caption{The diagram $\mc{R}(\mc{F}(D_k))$.}
    \label{fig: R-F-of-diagram}
\end{subfigure}
\end{figure}

In this paper, we need some common diagrams. Notice that for $D_k\in \mc{P}_k$, we can write the following diagram in figure \ref{fig: rainbowed-box} as $\te{ev}_X(X\otimes \te{ev}_X\otimes X)(X \otimes X \otimes \te{ev}_X \otimes X \otimes X)...(D_k\otimes X^{\otimes k})$. We will define $\mc{R}_k(D_k)\coloneq\te{ev}_X(X\otimes \te{ev}_X\otimes X)(X \otimes X \otimes \te{ev}_X \otimes X \otimes X)...(D_k\otimes X^{\otimes k})$. We can think of $\mc{R}$ as ``rainbowing" the strands on top of $D_k$ to the right. We will drop the subscript when the box space is clear. Another common diagram is a big cap over the rainbowed-box. That is, if $D_k \in \mc{P}_k$, $e_1(X\otimes \mc{R}(D_k)\otimes X)$. This looks like figure \ref{fig: over-the-rainbow} and is denoted as $\mc{O}(\mc{R}(D_k))$. Notice that for $D_k\in \mc{D}_k$ that $\mc{R}(\mc{F}(D_k))$ diagrammatically looks like figure \ref{fig: R-F-of-diagram}.

\begin{figure}[h!]
     \centering
     \begin{align*}  
    \adjustbox{scale=0.91}{\begin{overpic}[unit=1mm, scale=0.5]{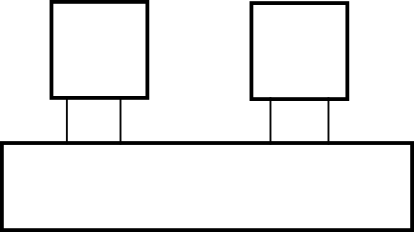}
     \put(12,23){$S$}
     \put(10.5, 15){$...$}
     \put(25, 15){$...$}
     \put(25,5){$D$}
     \put(38,23){$S$}
     \put(37, 15){$...$}
     \end{overpic}}
\end{align*}
     \caption{A train, where  $D\in \mc{TL}$ (and implicitly the strands of $S$ are rainbowed over, but drawn for convenience below the box). }
     \label{fig: train}
 \end{figure}

\begin{figure}[h]
    \centering
    \begin{align*}
     \adjustbox{scale=0.91}{\begin{overpic}[unit=1mm, scale=0.5]{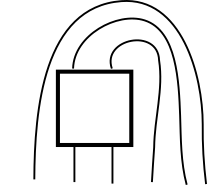}
     \put(11,9){$S$}
     \put(10.5,2){$...$}
     \put(20.5,2){$...$}
     \end{overpic}}
     \quad 
     {\makebox[0pt][l]{\raisebox{0.8cm}{$=\mathlarger{\mathlarger{\mathlarger{\sum_i}}}$}}}
     \quad \quad \quad \quad 
    \adjustbox{scale=0.91}{\begin{overpic}[unit=1mm, scale=0.5]{train.png}
     \put(12,23){$S$}
     \put(10.5, 15){$...$}
     \put(25, 15){$...$}
     \put(25,5){$D_i$}
     \put(38,23){$S$}
     \put(37, 15){$...$}
     \end{overpic}}
\end{align*}
    \caption{The relation showing a cap over the rainbowed $S$-box is a linear combination of trains, where each $D_i \in \mc{TL}$.}
    \label{fig:overstrand is a linear combinator of strands}
\end{figure}

 Let $\mc{P}$ be a planar algebra and $S\in \mc{P}_k$ be a box. A \tbf{train} (in $S$) \cite{BP14} is a diagram where the top boundary is a finite collection of rainbowed $S$-boxes (in particular, there are no strands capped over the $S$-boxes) on top of a diagram in Temperley-Lieb. That is, diagrammatically a train in $S$ looks like figure \ref{fig: train}. We say $S$ satisfies a \tbf{jellyfish relation} if $\mc{O}(\mc{R}(S))$ can be written as a linear combination of trains. Diagrammatically the jellyfish relation is shown in figure \ref{fig:overstrand is a linear combinator of strands}. While the definition of train first appeared in \cite{BP14}, it's equivalent notion in terms of $\lambda$-lattices first appeared in \cite{Pop95}. If a $k$-box $S$ in a planar algebra $\mc{P}$ has a jellyfish relation, we say $S$ is a \tbf{jellyfish generator (at depth $k$)}.

 {\begin{prop}\label{prop: BMPS} (Bigelow--Morrison--Peters--Snyder \cite{BPMS12}, Theorem 3.8) A planar algebra generated by a single box $S\in \mc{P}_k$ satisfying the below criteria 
    \begin{enumerate}[label=(\roman*)]
    \item There exists an $a \in \mbb{C}$ such that $\mc{F}(S)=aS$,
    \item $S$ has a jellyfish relation,
    \item there exists $b,c\in \mbb{C}$ such that $S^2=bS+cf^{(k)}$, and
    \item $S$ is U.C.C.S.,
\end{enumerate}
can be evaluated using the jellyfish algorithm.
\end{prop}}

\section{Results on General Index 4 Planar Algebras}\label{generalpafacts}
We first state and prove some results true for any index 4 planar algebras. 

\begin{rem*}\label{rem: Wenzlrelationassum}
    Wenzl's relation can be written as 
    \begin{align*}
        f^{(k+1)}=f^{(k)}\otimes X-\f{[k]_q}{[k+1]_q}(f^{(k)}\otimes X)e_k(f^{(k)}\otimes X).
    \end{align*}
\end{rem*}

\begin{lem}\label{Joneswenzlformulaupto3}
    For any index 4 planar algebra, $[k]_q=k$, $f^{(2)}=X\otimes X - \f{1}{2}e_1$, $f^{(3)}=X^{\otimes 3}-\f{2}{3}e_1-\f{2}{3}e_2+\f{1}{3}e_2e_1+\f{1}{3}e_1e_2$, and $\te{tr}(f^{(k)})=k+1$ for all $k\geq 0$.
\end{lem}

\begin{proof}
    Recall that $\te{tr}(f^{(k)})=[k+1]_q$ and $[k]_q=q^{k-1}+q^{k-3}+...+q^{-(k-3)}+q^{-(k-1)}$. We first prove by induction on $k\geq 0$ that $\te{tr}(f^{(k)})(=[k+1]_q)=k+1$.  We know that $\te{tr}(f^{(0)})=[1]_q=1$ and $\delta=[2]_q=q+q^{-1}=2$, which gives that $q=1$. Suppose the claim is true for $k\geq 1$. Then $[k+1]_q=q^{(k+1)-1}+[k-1]_q+q^{-((k+1)-1)}=1+(k-1)+1=k+1$. 

  By Wenzl's relation, $f^{(2)}=f^{(1)}\otimes X\cong \f{1}{2}(f^{(1)}\otimes X)e_1(f^{(1)}\otimes X)$. Since $f^{(1)}=X$, we get $f^{(2)}=X\otimes X - \f{1}{2}e_1$. 

     We employ the same method for $f^{(3)}$. Wenzl's relation gives that $f^{(3)}=f^{(2)}\otimes X-\f{2}{3}(f^{(2)}\otimes X)e_2(f^{(2)}\otimes X)$. Expanding out $f^{(2)}$ gives 
     \begin{eqnarray*}
         f^{(3)}&=(X\otimes X - \f{1}{2}e_1)\otimes X-\f{2}{3}((X\otimes X - \f{1}{2}e_1)\otimes X)e_2((X\otimes X - \f{1}{2}e_1)\otimes X)\\
         &=X^{\otimes 3}-\f{2}{3}e_1-\f{2}{3}e_2+\f{1}{3}e_2e_1+\f{1}{3}e_1e_2
     \end{eqnarray*}
    as we wished.
\end{proof}

We define a \tbf{leaf} of a graph to be a vertex with only 1 neighbor. When the context is clear, we may also refer to a \textbf{leaf} as the minimal projection associated to a leaf vertex. We call the following equality the \tbf{leaf property}. 

\begin{lem}\label{lem: leafproperty}
    Let $L\in \mc{P}_k$ be a leaf of a principal graph for any planar algebra of index 4. Then $L\otimes X=2 (L\otimes X)e_{k}(L\otimes X)$.
\end{lem}

\begin{proof}
    Let $P$ be the only neighbor of $L$ on the principal graph. Then $L\otimes X \cong P$. This gives $\te{Hom}(L\otimes X, L\otimes X)\cong \te{Hom}(P,P)$ is 1-dimensional since $P$ is a minimal projection. Since $(L\otimes X)e_{k}(L\otimes X)$ is in $\te{Hom}(L\otimes X, L\otimes X)$, this gives there exists an $a \in \C$ such that $(L\otimes X)e_{k}(L\otimes X)=a(L\otimes X)^2=a(L\otimes X)$. Taking the trace of both sides gives $a=\f{1}{2}$.
\end{proof}

We will frequently use that for all $k\geq 1$, the coefficient of the $X^{\otimes k}$ diagram in the expansion of $f^{(k)}$ is always 1. The rest of the terms involve cups and caps. See \cite{Mor17}, Lemma 2.1 for a proof of this fact.  

\begin{lem}
    Let $\mc{P}$ be an arbitrary planar algebra. If $D_k \in \mc{P}_k$ is uncuppable and uncappable, then $D_k(f^{(j)}\otimes X^{\otimes (k-j)})=(f^{(j)}\otimes X^{\otimes (k-j)})D_k=D_k$ for all $1\leq j \leq k$. If further $D_k$ is U.C.C.S., then $\mc{R}(D_k)(f^{(j)}\otimes X^{\otimes (2k-j)})=(f^{(j)}\otimes X^{\otimes (2k-j)})\mc{R}(D_k)=\mc{R}(D_k)$ for all $1\leq j \leq 2k$. 
\end{lem}

\begin{proof}
    Let $\mc{P}$ be an arbitrary planar algebra. Let $D_k\in \mc{P}_k$ be uncuppable and uncappable. All terms except the identity of $f^{(j)}$ has a cup and cap. The identity has coefficient 1. Thus $D_k(f^{(j)}\otimes X^{\otimes (k-j)})=D_k(X^{\otimes j}\otimes X^{(k-j)})=D_k$. The same idea proves the rest of the equalities. 
\end{proof}

The above properties proved for uncuppable, uncappable, and U.C.C.S. diagrams are called \tbf{absorption rules}. We say $D_k$ (respectively, $\mc{R}(D_k)$) \tbf{absorbs} $f^{(j)}$ (respectively, $f^{(2k-j)}$). 

\begin{lem}\label{lem: generalrecursiveabsorptionrule}
    If $Q$ is a projection in the $k$th box space of an arbitrary planar algebra and $P$ is defined by $P\coloneq Q\otimes X - a (Q\otimes X)e_k(Q \otimes X)$ for some $a\in \mbb{C}$, then we have the following absorption rule: $P(Q\otimes X)=(Q\otimes X)P=P$. 
\end{lem}

\begin{proof}
    By the definition of $P$, $P(Q\otimes X)=(Q\otimes X)^2-a(Q \otimes X)e_k(Q\otimes X)^2=Q\otimes X -a(Q\otimes X)e_k(Q \otimes X)=P$. The case $(Q\otimes X)P$ is proved the same way. 
\end{proof}

We call the rule proved in Lemma \ref{lem: generalrecursiveabsorptionrule}, the \tbf{recursive absorption rule}.

Define $g_{n,i}$ to be the diagram in Temperley-Lieb with $n$ strands on the top and bottom, a cap in the $(n-1,n)$ position and a cup in the $(i,i+1)$ position. Define $g_{n,n}$ to be the identity on $n$ strands. Below are $g_{4,4}, g_{4,3}, g_{4,2},$ and $g_{4,1}$ respectively.
\begin{align*}
     \adjustbox{scale=1.2}{ \begin{overpic}[unit=1mm, scale=0.5]{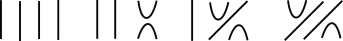}
     \put(8,0){,}
     \put(21,0){,}
     \put(33,0){,}
     \end{overpic}}
\end{align*}

We will occasionally employ the following theorem due to Frenkel and Khovanov \cite{FK97}, adapted for our index 4 case. Note that the statement below is somewhat different from that in their paper because of different conventions. (Also see Morrison \cite{Mor17}, Proposition 3.3, for another formulation of this statement.)

\begin{prop}\label{prop: MorrisonJWexpansion} (Frenkel--Khovanov \cite{FK97}, Theorem 3.5) Let $\mc{P}$ be an arbitrary planar algebra of index 4. For all $n\in \N$, the Jones-Wenzl projections satisfy
    \begin{align*}
        f^{(n+1)}=\f{n+1}{n+1}g_{n+1,n+1}(f^{(n)}\otimes X)-\f{n}{n+1}g_{n+1,n}(f^{(n)}\otimes X)+...\pm \f{1}{n+1}g_{n+1,1}(f^{(n)}\otimes X).
    \end{align*}
\end{prop}

\section{Sufficient Relations for Affine \texorpdfstring{$E_7$}{E7} Unshaded Subfactor Planar Algebras}\label{sufficientE7}

We begin by giving the sufficient relations for the affine $E_7$ (unshaded) subfactor planar algebra. In the next section, we show that these relations are necessary. Each box has a marked point on its left middle boundary. If drawn, it is denoted with a star.

\begin{defn}\label{defn: ourKauffmanbracket}
    Define the following:
    \begin{equation} \label{eq: firstcrossingrel}
       \adjustbox{scale=0.97}{\begin{overpic}[unit=1mm, scale=0.7]{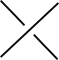}
    \end{overpic} } \quad {\makebox[0pt][l]{\raisebox{0.5cm}{{$\coloneq i$}}}}\quad \quad\adjustbox{scale=0.975}{\begin{overpic}[unit=1mm, scale=0.7]{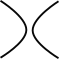}
    \end{overpic}}
    {\raisebox{0.5cm}{{$-i$}}}\hspace{1mm}\adjustbox{scale=0.975}{\begin{overpic}[unit=1mm, scale=0.7]{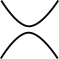}
    \end{overpic}}
\hspace{1mm}.
    \end{equation}
        As a consequence, this also defines,
        \begin{equation}\label{eq: secondcrossingrel}
       \begin{overpic}[unit=1mm, scale=0.7]{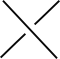}
    \end{overpic}  \quad {\makebox[0pt][l]{\raisebox{0.5cm}{{$\coloneq i$}}}}\quad \quad\begin{overpic}[unit=1mm, scale=0.7]{bsmoothing.png}
    \end{overpic}
    {\raisebox{0.5cm}{{$-i$}}}\hspace{1mm}\begin{overpic}[unit=1mm, scale=0.7]{asmoothing.png}
    \end{overpic}
\hspace{1mm}.
    \end{equation}

    We will denote the crossing in equation \ref{eq: firstcrossingrel} by $\mathfrak{X}$. The diagram in the expansion of $\mathfrak{X}$ with the coefficient $i$ is called an \tbf{A-smoothing}, denoted by $\mathfrak{A}$, and the diagram with the coefficient $-i$ is called a \tbf{B-smoothing}, denoted by $\mathfrak{B}$. The crossing in equation \ref{eq: secondcrossingrel} is denoted by $\mathfrak{Y}$. 
\end{defn}

We will focus this section on proving the below theorem.

\begin{thm}\label{thm: E7sufficientmaintheorem}
    Let $\mc{P}$ be the planar algebra generated by a self-adjoint box $S\in \mc{P}_4$: 
    \begin{align*}
        \adjustbox{scale=1}{ \begin{overpic}[unit=1mm, scale=0.4]{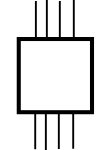}
        \put(4.5, 6.4){$S$}
    \end{overpic}},
    \end{align*}
    with the * extended anti-linearly and on diagrams, modulo the following relations:
    \begin{enumerate}[label=(\roman*)]
         \item (\textbf{the bubble relation})  $\adjustbox{scale=1}{\bubble}=2$,
          \item (\textbf{$S$ is U.C.C.S.})
      \quad {\makebox[0pt][l]{\raisebox{-0.6cm}{{$\begin{aligned}[t] \adjustbox{scale=1}{ \begin{overpic}[unit=1mm, scale=0.4]{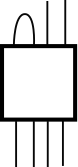}
        \put(2.7, 7.4){$S$}
    \end{overpic} }
    \quad {\makebox[0pt][l]{\raisebox{0.6cm}{{$=...=$}}}}\quad\quad\quad
    \adjustbox{scale=1}{\begin{overpic}[unit=1mm, scale=0.4]{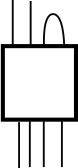}
        \put(2.9, 7.5){$S$}
    \end{overpic}}
     \quad {\makebox[0pt][l]{\raisebox{0.6cm}{{$=$}}}}\quad
      \adjustbox{scale=1}{ \begin{overpic}[unit=1mm, scale=0.4]{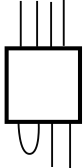}
        \put(3.4, 7.5){$S$}
    \end{overpic} }
    \quad {\makebox[0pt][l]{\raisebox{0.6cm}{{$=...=$}}}}\quad\quad\quad
     \adjustbox{scale=1}{\begin{overpic}[unit=1mm, scale=0.4]{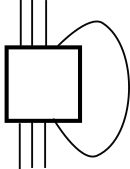}
        \put(3, 7.5){$S$}
    \end{overpic} }\quad {\makebox[0pt][l]{\raisebox{0.6cm}{{$=$}}}}
    \quad 
    \adjustbox{scale=1}{ \begin{overpic}[unit=1mm, scale=0.4]{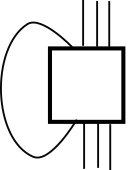}
        \put(8.4, 7.5){$S$}
    \end{overpic}} \quad {\makebox[0pt][l]{\raisebox{0.6cm}{{$=0$,}}}}
\end{aligned}$}}}}
    \item (\textbf{$S^2$ relation}) 
          {\makebox[0pt][l]{\raisebox{-0.6cm}{{$\begin{aligned}[t] \adjustbox{scale=1}{\begin{overpic}[unit=1mm, scale=0.4]{S-4strand.png}
        \put(4.5, 6.4){$S^2$}
    \end{overpic}}
     \quad {\raisebox{0.6cm}{{$=6$}}}
     \adjustbox{scale=1}{\begin{overpic}[unit=1mm, scale=0.4]{S-4strand.png}
        \put(3.9, 6.4){$f^{(4)}$}
    \end{overpic}}
    \quad  {\raisebox{0.6cm}{{$+$}}}
     \adjustbox{scale=1}{\begin{overpic}[unit=1mm, scale=0.4]{S-4strand.png}
        \put(4.5, 6.4){$S$}
    \end{overpic}}
    \end{aligned}$}}}}
    \item (\textbf{Jellyfish relation})
     {\makebox[0pt][l]{\raisebox{-0.6cm}{{$\begin{aligned}[t] \adjustbox{scale=1}{\begin{overpic}[unit=1mm, scale=0.4]{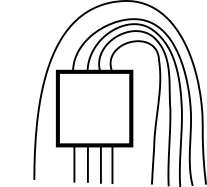}
       \put(8.5, 6.8){$S$}
    \end{overpic}}\quad {\raisebox{0.6cm}{{$=$}}}
    \adjustbox{scale=1}{ \begin{overpic}[unit=1mm, scale=0.4]{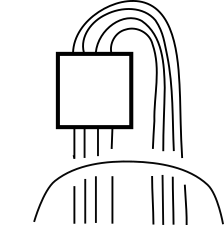}
        \put(8.4, 12.9){$S$}
    \end{overpic}}
    \end{aligned}$}}}}
    \end{enumerate}
    Then $\mc{P}$ is an unshaded subfactor planar algebra with the below principal graph:
    \begin{equation}\label{eq: generalE7graph}
    \arbEsevengraph
\end{equation}
    where $P_1 = f^{(1)}=X$, $P_2 =f^{(2)}$, $P_3 =f^{(3)}$, $P_4=\f{3}{5}f^{(4)}-\f{1}{5}S$, $P_5 = P_4 \otimes X - \f{4}{3}(P_4 \otimes X)e_4 (P_4\otimes X)$, $P_6= P_5 \otimes X - \f{3}{2}(P_5 \otimes X)e_5(P_5\otimes X)$, and $Q_4 =\f{2}{5}f^{(4)}+\f{1}{5}S$.
\end{thm}

Define the planar algebra $\mc{P}$ and elements $P_1, ..., P_6,$ and $Q_4$ as in the preceding theorem. We first show that these elements are projections.

\begin{lem}
    $Q_4$, $P_1$, $P_2$, $P_3$, and $P_4$ are projections. Further, they are all uncuppable and uncappable. 
\end{lem}

\begin{proof}
    $Q_4$ is a linear combination of two self-adjoint, uncuppable, and uncappable elements so is itself self-adjoint, uncuppable, and uncappable. Further, using that $S$ is U.C.C.S. and the $S^2$ relation,
    \begin{align*}
        Q_4^2&=\f{1}{25}(4f^{(4)}+4S+S^2)=\f{1}{25}(4f^{(4)}+4S+6f^{(4)}+S)
        =\f{1}{25}(10f^{(4)}+5S)=\f{2}{5}f^{(4)}+\f{1}{5}S=Q_4
    \end{align*}
    so $Q_4$ is a projection.

    Since $P_1$, $P_2$, and $P_3$ are Jones-Wenzl projections, all that is left to show is that $P_4$ is a uncuppable and uncappable projection. $P_4$ is a linear combination of two self-adjoint, uncuppable, and uncappable elements so is itself self-adjoint, uncuppable, and uncappable. Further, using that $S$ is U.C.C.S. and the $S^2$ relation, we obtain
    \begin{align*}
        P_4^2&=\f{1}{25}(9f^{(4)}-6S+S^2)=\f{1}{25}(9f^{(4)}-6S+6f^{(4)}+S)=\f{1}{25}(15f^{(4)}-5S)=\f{3}{5}f^{(4)}-\f{1}{5}S=P_4
    \end{align*}
    so $P_4$ is also a projection.
\end{proof}

\begin{lem}\label{lem: e7sufficientpartialtraces}
    We have the following partial traces:
    \begin{enumerate}[label=(\roman*)]
        \item $E_2(f^{(2)})=\f{3}{2}f^{(1)}=\f{3}{2}X$,
        \item $E_3(f^{(3)})=\f{4}{3}f^{(2)}$,
        \item $E_4(f^{(4)})=\f{5}{4}f^{(3)}$,
        \item $E_5(P_5)=\f{2}{3}P_4$,
        \item $E_4(Q_4)=\f{1}{2}f^{(3)}$, and
        \item $E_4(P_4)=\f{3}{4}f^{(3)}$.
    \end{enumerate}
\end{lem}

\begin{proof}
    By Wenzl's relation, 
    $E_2(f^{(2)})=2f^{(1)}-\f{1}{2}{f^{(1)}}^2=\f{3}{2}f^{(1)}$, $E_3(f^{(3)})=2f^{(2)}-\f{2}{3}{f^{(2)}}^2=\f{4}{3}f^{(2)}$, and $E_4(f^{(4)})=2f^{(3)}-\f{3}{4}{f^{(3)}}^2=\f{5}{4}f^{(3)}$. Using that $P_4$ is a projection, $E_5(P_5)=2P_4-\f{4}{3}P_4^2=2P_4-\f{4}{3}P_4=\f{2}{3}P_4$. Using that $S$ is U.C.C.S., $E_4(Q_4)=\f{2}{5}E_4(f^{(4)})+\f{1}{5}E_4(S)=\f{2}{5}\cdot \f{5}{4}f^{(3)}=\f{1}{2}f^{(3)}$. Similarly, $E_4(P_4)=\f{3}{5}E_4(f^{(4)})-\f{1}{5}E_4(S)=\f{3}{5}\cdot \f{5}{4}f^{(3)}=\f{3}{4}f^{(3)}$. 
\end{proof}

\begin{lem}
    $P_5$ is a projection. 
\end{lem}

\begin{proof}
    Since $P_4$ is self-adjoint, so is $P_5$. Further,
    \begin{align*}
        P_5^2&=(P_4\otimes X)^2-\f{4}{3}(P_4 \otimes X)^2e_4(P_4 \otimes X)-\f{4}{3}(P_4 \otimes X)e_4(P_4 \otimes X)^2+\f{16}{9}(P_4 \otimes X)(X^{\otimes 3}\otimes \te{coev}_X)E_4(P_4)(X^{\otimes 3}\otimes \te{ev}_X)(P_4 \otimes 4)\\
        &=(P_4\otimes X)-\f{4}{3}(P_4 \otimes X)e_4(P_4 \otimes X)-\f{4}{3}(P_4 \otimes X)e_4(P_4 \otimes X)
        +\f{16}{9}\cdot \f{3}{4}(P_4 \otimes X)(X^{\otimes 3}\otimes \te{coev}_X)f^{(3)}(X^{\otimes 3}\otimes \te{ev}_X)(P_4 \otimes 4)\\
        &=(P_4\otimes X)-\f{8}{3}(P_4 \otimes X)e_4(P_4 \otimes X)+\f{4}{3}(P_4 \otimes X)(X^{\otimes 3}\otimes \te{coev}_X)f^{(3)}(X^{\otimes 3}\otimes \te{ev}_X)(P_4 \otimes 4)
    \end{align*}
    $P_4$ is uncuppable and uncappable, so the last term in the above sum becomes $\f{4}{3}(P_4 \otimes X)e_4(P_4 \otimes X)$.
    This gives $P_5^2=P_4 \otimes X-\f{4}{3}(P_4 \otimes X)e_4(P_4 \otimes X)=P_5$, as we wished.  
\end{proof}

\begin{lem}\label{lem: e7sufficientpartialtraceofp6}
    The partial trace of $P_6$ is $\f{1}{2}P_5$. That is, $E_6(P_6)=\f{1}{2}P_5$. Also $P_6$ is a projection.
\end{lem}

\begin{proof}
    Since $P_5^2=P_5$, we have that $E_6(P_6)=2P_5-\f{3}{2}P_5^2=\f{1}{2}P_5$. We can calculate the square of $P_6$ using the partial trace of $P_5$ and that $P_5$ is a projection:
    \begin{align*}
        P_6^2&=(P_5 \otimes X)^2-\f{3}{2}(P_5\otimes X)^2e_5(P_5 \otimes X)-\f{3}{2}(P_5\otimes X)e_5(P_5 \otimes X)^2+\f{9}{4}(P_5 \otimes X)(X^{\otimes 4}\otimes \te{coev}_X)E_5(P_5)(X^{\otimes 4}\otimes \te{ev}_X)(P_5 \otimes X)\\
        &=(P_5 \otimes X)-\f{6}{2}(P_5\otimes X)e_5(P_5 \otimes X)+\f{9}{4}\cdot \f{2}{3}(P_5 \otimes X)(X^{\otimes 4}\otimes \te{coev}_X)P_4(X^{\otimes 4}\otimes \te{ev}_X)(P_5 \otimes X)
    \end{align*}
    \begingroup
    \allowdisplaybreaks
    The recursive absorption rule of Lemma \ref{lem: generalrecursiveabsorptionrule} gives that
    \begin{align*}
        P_6^2&=(P_5 \otimes X)-\f{6}{2}(P_5\otimes X)e_5(P_5 \otimes X)+\f{9}{4}\cdot \f{2}{3}(P_5 \otimes X)(X^{\otimes 4}\otimes \te{coev}_X)(X^{\otimes 4}\otimes \te{ev}_X)(P_5 \otimes X)\\
        &=(P_5 \otimes X)-\f{6}{2}(P_5 \otimes X)e_5(P_5 \otimes X)+\f{3}{2}(P_5 \otimes X)e_5(P_5 \otimes X)\\
        &=P_5 \otimes X -\f{3}{2}(P_5 \otimes X)e_5(P_5 \otimes X)=P_6.
    \end{align*}
    \endgroup
    So indeed $P_6$ is a projection.
\end{proof}

\begin{lem}\label{lem: E7sufficienttraces}
    The traces of the projections are as follows.
    \begin{enumerate}[label=(\roman*)]
        \item $\te{tr}(\emptyset)=1$,
        \item $\te{tr}(P_1)=2$,
        \item $\te{tr}(P_2)=3$,
        \item $\te{tr}(P_3)=4$,
        \item $\te{tr}(P_4)=3$,
        \item $\te{tr}(P_5)=2$,
        \item $\te{tr}(P_6)=1$, and
        \item $\te{tr}(Q_4)=2$.
    \end{enumerate}
\end{lem}

    \begin{proof}
        By convention $\te{tr}(\emptyset)=1$. By the bubble relation, $\te{tr}(P_1)=2$. We use the partial traces found in Lemma \ref{lem: e7sufficientpartialtraces} and \ref{lem: e7sufficientpartialtraceofp6} to compute  
        \begingroup
    \allowdisplaybreaks
    \begin{align*}
        \te{tr}(P_2)&=\te{tr}(E_2(P_2))=\f{3}{2}\te{tr}(P_1)=3,\\
        \te{tr}(P_3)&=\te{tr}(E_3(P_3))=\f{4}{3}\te{tr}(P_2)=4,\\
        \te{tr}(P_4)&=\te{tr}(E_4(P_4))=\f{3}{4}\te{tr}(P_3)=3,\\
        \te{tr}(P_5)&=\te{tr}(E_5(P_5))=\f{2}{3}\te{tr}(P_4)=2,\\
        \te{tr}(P_6)&=\te{tr}(E_6(P_6))=\f{1}{2}\te{tr}(P_5)=1, \te{ and }\\
        \te{tr}(Q_4)&=\te{tr}(E_4(Q_4))=\f{1}{2}\te{tr}(P_3)=2,
    \end{align*}
    \endgroup
    as we wished.
\end{proof}

\begin{lem}
   We have the following relationships between projections:
   \begin{enumerate}[label=(\roman*)]
       \item $P_5(P_4\otimes X)=(P_4\otimes X)P_5=P_5$,
       \item $P_6(P_5 \otimes X)=(P_5\otimes X)P_6=P_6$,
       \item $P_4(f^{(3)}\otimes X)=(f^{(3)}\otimes X)P_4=P_4$, and
       \item $P_6(P_4 \otimes X^{\otimes 2})=(P_4\otimes X^{\otimes 2})P_6=P_6$.
   \end{enumerate}
\end{lem}

\begin{proof}
    The first two parts follow from the recursive absorption rule. The third part follows from $P_4$ being uncuppable and uncappable. The fourth part is clear using the previous parts as
    \begin{align*}
        P_6(P_4 \otimes X^{\otimes 2})=P_6(P_5\otimes X)(P_4\otimes X^{\otimes 2})=P_6(P_5 \otimes X)=P_6
    \end{align*}
    A symmetric argument also shows $(P_4 \otimes X^{\otimes 2})P_6=P_6$. 
\end{proof}

\begin{lem}
    $P_5$ and $P_6$ are uncuppable and uncappable.
\end{lem}

\begin{proof}
    For $P_5,$ if $1\leq k\leq 3$, then $e_kP_5=P_5e_k=0$, since $P_4$ is uncuppable and uncappable. What's left to show is that $e_4P_5=P_5e_4=0$. By symmetry, we will only compute $e_4P_5$. Notice,
    \begin{align*}
        e_4P_5&=e_4(P_4 \otimes X)-\f{4}{3}(X^{\otimes 3}\otimes \te{coev}_X) E_4(P_4)(X^{\otimes 3}\otimes \te{ev}_X)(P_4\otimes X)=e_4(P_4 \otimes X)-\f{4}{3}\cdot \f{3}{4}e_4(P_4\otimes X)=0.
    \end{align*}

    For $P_6$, we use that $P_5$ is uncuppable and uncappable to obtain that $e_kP_6=P_6e_k=0$ for $1\leq k \leq 4$. By symmetry, all that is left to show is $e_5P_6=0$. Using the partial trace of $P_5:$
    \begin{align*}
        e_5P_6&=e_5(P_5 \otimes X)-\f{3}{2}(X^{\otimes 4}\otimes \te{coev}_X)E_5(P_5)(X^{\otimes 4}\otimes \te{ev}_X)(P_5 \otimes X)\\
        &=e_5(P_5 \otimes X)-\f{3}{2}\cdot \f{2}{3}(X^{\otimes 4}\otimes \te{coev}_X)P_4 (X^{\otimes 4} \otimes \te{ev}_X)(P_5 \otimes X)
    \end{align*}
    Using the absorption rule of Lemma \ref{lem: generalrecursiveabsorptionrule}, we obtain $e_5P_6=0$.
\end{proof}

\begin{lem}\label{lem: e7sufficienttensorproductsofP}
    In $\mc{P}$, we have the following tensor products:
    \begin{enumerate}[label=(\roman*)]
        \item $\emptyset \otimes X \cong P_1$,
        \item $P_1 \otimes X \cong \emptyset \oplus P_2$,
        \item $P_2 \otimes X \cong P_1 \oplus P_3$,
        \item $P_3 \otimes X\cong P_2 \oplus P_4 \oplus Q_4$,
        \item $P_4 \otimes X \cong P_3 \oplus P_5$, and 
        \item $P_5 \otimes X \cong P_4 \oplus P_6$.
    \end{enumerate}
\end{lem}

\begin{proof}
    Part (i), (ii), and (iii) follow from Wenzl's relation. For part (iv), we can notice that $P_4+Q_4=f^{(4)}$, so this also follows from Wenzl's relation.

    For part (v), we claim the explicit isomorphisms are:
    \begin{align*}
        &A_1\coloneq \begin{bmatrix}
            (P_4 \otimes X)(X^{\otimes 3}\otimes \te{coev}_X)P_3 & (P_4 \otimes X)P_5
        \end{bmatrix}: P_3 \oplus P_5 \to P_4 \otimes X, \te{ and }\\
        &B_1\coloneq \begin{bmatrix}
        \f{4}{3}P_3(X^{\otimes 3}\otimes \te{ev}_X)(P_4 \otimes X) \\ P_5(P_4 \otimes X)
        \end{bmatrix}: P_4 \otimes X \to P_3 \oplus P_5
    \end{align*}
    Notice, we can more conveniently recognize
    \begin{align*}
        A_1=\begin{bmatrix}
            (P_4\otimes X)(X^{\otimes 3}\otimes  \te{coev}_X) & P_5 
        \end{bmatrix}
        \te{ and }
        B_1 =\begin{bmatrix}
            \f{4}{3} (X^{\otimes 3}\otimes \te{ev}_X)(P_4 \otimes X) \\ P_5
        \end{bmatrix}
    \end{align*}
    Then, we can see that $A_1B_1=\f{4}{3}(P_4 \otimes X)e_4(P_4\otimes X)+P_5=P_4\otimes X$. On the other hand, 
    \begin{align*}
        B_1A_1=\begin{bmatrix}
            \f{4}{3}E_4(P_4) & \f{4}{3}(X^{\otimes 4}\otimes \te{ev}_X)(P_4 \otimes X)P_5 \\ 
            P_5(P_4 \otimes X)(X^{\otimes 3}\otimes \te{coev}_X) & P_5
        \end{bmatrix}
    \end{align*}
    By the recursive absorption rule, the off-diagonal entries are 0. Since $E_4(P_4)=\f{3}{4}f^{(3)}$, we obtain $B_1A_1=P_3 \oplus P_5$. 

    For part (vi) we claim the explicit isomorphisms are:
    \begin{align*}
        &A_2\coloneq \begin{bmatrix}
            (P_5\otimes X)(X^{\otimes 4}\otimes \te{coev}_X)P_4 & (P_5 \otimes X)P_6
        \end{bmatrix}=\begin{bmatrix}
            (P_5 \otimes X)(X^{\otimes 4}\otimes \te{coev}_X) & P_6
        \end{bmatrix}\\
        &B_2 \coloneq \begin{bmatrix}
            \f{3}{2} P_4(X^{\otimes 4}\otimes \te{ev}_X)(P_5 \otimes X)\\
            P_6(P_5 \otimes X)
        \end{bmatrix} = \begin{bmatrix}
            \f{3}{2}(X^{\otimes 4}\otimes \te{ev}_X)(P_5 \otimes X) \\ P_6
        \end{bmatrix}.
    \end{align*}
    So $A_2\in \te{Hom}(P_4 \oplus P_6, P_5 \otimes X)$ and $B_2 \in \te{Hom}(P_5 \otimes X, P_4 \oplus P_6)$. We also have that 
    \begin{align*}
        A_2B_2=\f{3}{2}(P_5 \otimes X)e_5(P_5 \otimes X)+P_6 = P_5 \otimes X
    \end{align*}
    and 
    \begin{align*}
        B_2A_2 =\begin{bmatrix}
            \f{3}{2}E_5(P_5) & \f{3}{2} (X^{\otimes 4}\otimes \te{ev}_X)(P_5 \otimes X)P_6 \\
            P_6(P_5 \otimes X)(X^{\otimes 4}\otimes \te{coev}_X) & P_6
        \end{bmatrix}
    \end{align*}
    By the recursive absorption rule, the off-diagonals are zero. Since $E_5(P_5)=\f{2}{3}P_4$, we obtain $B_2A_2=P_4 \oplus P_6$. Thus, $P_4\oplus P_6 \cong P_5 \otimes X$. 
\end{proof}

\subsection{Results Involving Crossings in an Index 4 Planar Algebra}
Next, we give a collection of proofs for general planar algebras with a defined crossing. Recall a classical theorem of Reidemeister \cite{Rei27} that says two knots (or links) are equivalent if any only if they differ by a finite sequence of the below \textbf{Reidemeister moves} and planar isotopies:
\begin{align*}
{\makebox[0pt][l]{\raisebox{-0.4cm}{{$\begin{aligned}[t]
    &{\raisebox{0.4cm}{{\tbf{Reidemeister I:}}}} \hspace{3mm} 
    \adjustbox{scale=1.1}{\begin{overpic}[unit=1mm, scale=0.7]{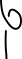}
    \end{overpic}}
    \quad {\raisebox{0.4cm}{{$\leftrightarrow$}}}\quad\adjustbox{scale=1.1}{\begin{overpic}[unit=1mm, scale=0.7]{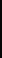}
    \end{overpic}}\quad {\raisebox{0.4cm}{{$\leftrightarrow$}}}\quad
    \adjustbox{scale=1.1}{\begin{overpic}[unit=1mm, scale=0.7]{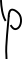}
    \end{overpic}}
    \\
    &{\raisebox{0.2cm}{{\tbf{Reidemeister II:}}}} \hspace{3mm} \adjustbox{scale=1.1}{\begin{overpic}[unit=1mm, scale=0.7]{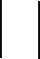}
    \end{overpic}}\quad {\raisebox{0.4cm}{{$\leftrightarrow$}}}\quad
   \adjustbox{scale=1.1}{ \begin{overpic}[unit=1mm, scale=0.7]{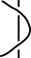}
    \end{overpic}}
    \end{aligned}$}}}}\\
    & {\raisebox{0.4cm}{{\tbf{Reidemeister III:}}}} \hspace{3mm}\adjustbox{scale=1.1}{ \begin{overpic}[unit=1mm, scale=0.7]{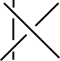}
    \end{overpic}}\quad {\raisebox{0.4cm}{{$\leftrightarrow$}}}\quad
    \adjustbox{scale=1.1}{\begin{overpic}[unit=1mm, scale=0.7]{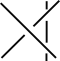}
    \end{overpic}} \quad
    {\raisebox{0.4cm}{\text{and}}}\quad 
   \adjustbox{scale=1.1}{ \begin{overpic}[unit=1mm, scale=0.7]{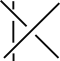}
    \end{overpic}}\quad {\raisebox{0.4cm}{{$\leftrightarrow$}}}\quad
    \adjustbox{scale=1.1}{\begin{overpic}[unit=1mm, scale=0.7]{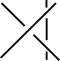}
    \end{overpic}}
\end{align*}

\begin{lem}
    Definition \ref{defn: ourKauffmanbracket} agrees with the bubble relation and Reidemeister II and III.  
\end{lem}

\begin{proof}
    Since we have for $q=1$, $\te{tr}(X)=q+q^{-1}$, and $\mathfrak{X}=\sqrt{-q}\mathfrak{A}+\sqrt{-q}^{-1}\mathfrak{B}$, we exactly have the Kauffman bracket, first defined by Kauffman in \cite{Kau87}. This is well known in knot theory to satisfy Reidemeister II and III. 
\end{proof}

We do not have Reidemeister I. Instead, we have the following. 

\begin{lem}\label{lem: cuponcrossingresolved}
    If an arbitrary planar algebra has crossings defined as in Definition \ref{defn: ourKauffmanbracket}, then we have the following rules for the crossings: 
    \begin{align*}
        \adjustbox{scale=1.1}{\begin{overpic}[unit=1mm, scale=0.7]{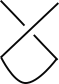}
    \end{overpic}}\quad {\raisebox{0.8cm}{{$=-i$}}}\hspace{2mm}
    {\raisebox{0.8cm}{\adjustbox{scale=1.1}{\begin{overpic}[unit=1mm, scale=0.7]{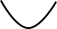}
    \end{overpic}}}}
    \quad {\raisebox{0.8cm}{{\text{ and }}}}\hspace{2mm}
    \adjustbox{scale=1.1}{\begin{overpic}[unit=1mm, scale=0.7]{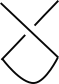}
    \end{overpic}}\quad {\raisebox{0.8cm}{{$=i$}}}\hspace{2mm}
    {\raisebox{0.8cm}{\adjustbox{scale=1.1}{\begin{overpic}[unit=1mm, scale=0.7]{cup.png}
    \end{overpic}}}}.
    \end{align*}
\end{lem}

\begin{proof}
    A direct application of Definition \ref{defn: ourKauffmanbracket}.
\end{proof}

\begin{lem}\label{lem: boxwithovercrossingonbottom}
    Let $R$ be an uncuppable box in the $n$th box space of an arbitrary planar algebra with crossings defined as in Definition \ref{defn: ourKauffmanbracket}. Then
    \begin{align*}
      \adjustbox{scale=1.1}{  \begin{overpic}[unit=1mm, scale=0.5]{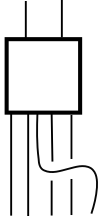}
        \put(7.2,2){\tiny{$...$}}
        \put(6.8,-2){\tiny{$(m)$}}
        \put(1.8,2){\tiny{$...$}}
        \put(4,17){$R$}
    \end{overpic}}\quad {\raisebox{0.8cm}{{$=(-i)^mR$}}}
    \quad {\raisebox{0.8cm}{{\text{ and }}}}\hspace{2mm}
    \adjustbox{scale=1.1}{ \begin{overpic}[unit=1mm, scale=0.5]{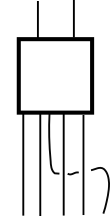}
        \put(8.9,2){\tiny{$...$}}
        \put(7.2,-2){\tiny{$(m)$}}
        \put(3.4,2){\tiny{$...$}}
        \put(5.3,17){$R$}
    \end{overpic}}\quad {\raisebox{0.8cm}{{$=i^mR$}}}
    \end{align*}
    where $0\leq m\leq n-1$ indicates the number of strands crossed over or under.
\end{lem}

\begin{proof}
    We will prove both cases by induction on $m$. The case when $m=0$ is clear. Assume the lemma is true for $0\leq m <n-1$. Then for $m+1$, resolve the left most crossing. When strand is overcrossing $m+1$ strands, the crossing is of the form $\mathfrak{Y}$. Recall $\mathfrak{Y}=i\mathfrak{B}-i\mathfrak{A}$. The $B$-smoothing will cup $R$, and thus that term will be zero. The $A$-smoothing will result in the strand overcrossing $m$ strands, which by induction is $(-i)^mR$. Multiplying by the $-i$ coefficient will then result in $(-i)^{m+1}R$. The undercrossing case is nearly the same. In this case, we have a crossing of the form $\mathfrak{X}$. The $A$-smoothing has a coefficient of $i$ and the $B$-smoothing with cap $R$.
\end{proof}

The same idea is used for the next proof, so will be omitted. 

\begin{lem}\label{lem: boxwithcrossingofstrandtotheleft}
    Let $R$ be an uncuppable box in the $n$th box space of an arbitrary planar algebra with crossing defined as in Definition \ref{defn: ourKauffmanbracket}. Then
     \begin{align*}
        \adjustbox{scale=1.1}{ \begin{overpic}[unit=1mm, scale=0.5]{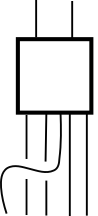}
        \put(9.5,2){\tiny{$...$}}
        \put(3.2,-2){\tiny{$(m)$}}
        \put(4,2){\tiny{$...$}}
        \put(5.5,17){$R$}
    \end{overpic}}\quad {\raisebox{0.8cm}{{$=i^mR$}}}
    \quad {\raisebox{0.8cm}{{\text{ and }}}}\hspace{2mm}
    \adjustbox{scale=1.1}{ \begin{overpic}[unit=1mm, scale=0.5]{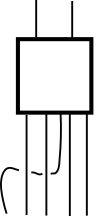}
        \put(9.6,2){\tiny{$...$}}
        \put(3,-2){\tiny{$(m)$}}
        \put(4,2){\tiny{$...$}}
        \put(5.5,17){$R$}
    \end{overpic}}
    \quad {\raisebox{0.8cm}{{$=(-i)^m R$}}}
    \end{align*}
    where $0\leq m\leq n$ indicates the number of strands crossed over or under.
\end{lem}

\begin{lem}\label{lem: overcrossingresolve}
    When resolving all the crossings of the below diagram,
    \begin{equation}\label{eqn: overcrossing2strands}
        \adjustbox{scale=0.9}{\begin{overpic}[unit=1mm, scale=0.7]{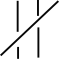}
        \put(4,2){\tiny{$...$}}
         \put(4,-2){\tiny{$(c)$}}
        \end{overpic}}
    \end{equation}
    where $c\geq 1$ is the number of strands crossed over, the coefficient of the term
     \begin{equation}\label{eqn: captwodiagcup}
        \adjustbox{scale=0.9}{\begin{overpic}[unit=1mm, scale=0.7]{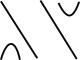}
        \put(8,2){\tiny{$...$}}
         \put(8,-2){\tiny{$(c)$}}
        \end{overpic}}
    \end{equation}
    is $(-i)^c$. Similarly, when resolving all the crossings of 
    \begin{equation}\label{eqn: overcrossing2strands2}
        \adjustbox{scale=0.9}{\begin{overpic}[unit=1mm, scale=0.7]{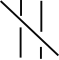}
        \put(4.5,2){\tiny{$...$}}
         \put(4,-2){\tiny{$(c)$}}
        \end{overpic}}
    \end{equation}
    where $c\geq 1$ is the number of strands crossed over, the coefficient of the term
     \begin{equation}\label{eqn: captwodiagcup2}
        \adjustbox{scale=0.9}{\begin{overpic}[unit=1mm, scale=0.7]{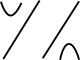}
        \put(4,2){\tiny{$...$}}
         \put(2.5,-2){\tiny{$(c)$}}
        \end{overpic}}
    \end{equation}
    is $i^c$. 
    
\end{lem}

\begin{proof}
    We prove by induction on $c$. When $c=1$, there is a single crossing whose resolution is given by Definition \ref{defn: ourKauffmanbracket}. Suppose the result is true for $c\geq 1$. By resolving the left-most crossing we obtain
    \begin{equation}\label{eqn: overcrossing2strandsexpansion}
         \adjustbox{scale=0.9}{\begin{overpic}[unit=1mm, scale=0.7]{overcrossing2strands.png}
        \put(4,2){\tiny{$...$}}
         \put(2.5,-2){\tiny{$(c+1)$}}
        \end{overpic}}
        \quad {\raisebox{0.4cm}{{$=i$}}}\hspace{2mm}
    \adjustbox{scale=0.9}{\begin{overpic}[unit=1mm, scale=0.7]{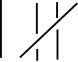}
    \put(7.6,2){\tiny{$...$}}
         \put(7.5,-2){\tiny{$(c)$}}
    \end{overpic}}
    \quad {\raisebox{0.4cm}{{$-i$}}}\hspace{2mm}
     \adjustbox{scale=0.9}{\begin{overpic}[unit=1mm, scale=0.7]{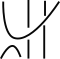}
     \put(5.5,2){\tiny{$...$}}
         \put(5.4,-2){\tiny{$(c)$}}
    \end{overpic}}
    \end{equation}
    The first term on the right-hand side will not have a term of the form (\ref{eqn: captwodiagcup}) since the left-most strand is straight with no crossings. Thus, we focus our attention to the second term on the right-hand side of the above equality. Notice that the term that will result in the form (\ref{eqn: captwodiagcup}) for $c+1$ will be the term where diagram in (\ref{eqn: overcrossing2strands}) resolves to the diagram in (\ref{eqn: captwodiagcup}) for $c$. By induction, the coefficient of this term is $(-i)^c$. Thus, due to the $-i$ in front of the term in equation \ref{eqn: overcrossing2strandsexpansion}, the coefficient will be $(-i)^{c+1}$ for the desired term. The other case is nearly identical in proof.
\end{proof}

\begin{lem}\label{lem: switchevennumberofcrossings}
    Switching an even number of crossings in any diagram of the planar algebra $\mc{P}$ results in the same diagram.
\end{lem}

\begin{proof}
    Suppose there are two local neighborhoods, $N_1$ and $N_2$, that just consist of a crossing. Denote the crossings or smoothings in these neighborhoods by multiplication. That is, $\mathfrak{X}\mathfrak{Y}$ indicates that in $N_1$ is an $\mathfrak{X}$ crossing and in $N_2$ is a $\mathfrak{Y}$ crossing. It is an easy computation to see that $\mathfrak{X}\mathfrak{X}=\mathfrak{Y}\mathfrak{Y}=- \mathfrak{B}\mathfrak{B}+\mathfrak{B}\mathfrak{A}+\mathfrak{A}\mathfrak{B}-\mathfrak{A}\mathfrak{A}=-\mathfrak{X}\mathfrak{Y}=-\mathfrak{Y}\mathfrak{X}$.
\end{proof}

\subsection{Equivalences of Diagrams in $\mathcal{P}$.}

In this subsection, we first introduce notation which will be used for the rest of the paper. Then we compute equivalences of a variety of diagrams. None of these diagrams are closed. We reserve subsection \ref{sec: closed-diagrams} for the computations of closed diagrams.

\begin{defn}
    We introduce the following notation for convenience. Let $R$ be a box in the $n$th box space of a planar algebra. 
    \begin{enumerate}
        \item We will denote the following for $0\leq m\leq n$:
        \begin{align*}
       \adjustbox{scale=1.1}{  \begin{overpic}[unit=1mm, scale=0.4]{box.png}
        \put(4, 6.5){$R$}
        \put(4.3, 14){$...$}
        \put(4.1, 2){$...$}
        \put(1.9,-2){\tiny{$(2n-m)$}}
    \end{overpic}}\quad {\raisebox{0.4cm}{{$\coloneq$}}}\hspace{2mm}
         \adjustbox{scale=1.1}{ \begin{overpic}[unit=1mm, scale=0.4]{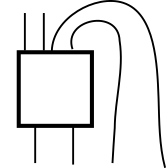}
        \put(4.5,2){\tiny{$...$}}
         \put(11.5,-2){\tiny{$(n-m)$}}
         \put(13.5, 2){\tiny{$...$}}
         \put(4, 7){$R$}
        \end{overpic}}
        \quad {\raisebox{0.4cm}{{\text{and}}}}\hspace{2mm}
       \adjustbox{scale=1.1}{  \begin{overpic}[unit=1mm, scale=0.4]{box.png}
        \put(4, 6.5){$R$}
        \put(4.3, 14){$...$}
        \put(4.1, 2){$...$}
        \put(2,-2){\tiny{$(n-m)$}}
    \end{overpic}}\quad {\raisebox{0.4cm}{{$\coloneq$}}}\hspace{2mm}
   \adjustbox{scale=1.1}{  \begin{overpic}[unit=1mm, scale=0.4]{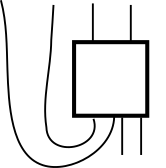}
        \put(2,14){\tiny{$...$}}
         \put(0.8,19){\tiny{$(m)$}}
         \put(10.5, 15){\tiny{$...$}}
         \put(10.5, 8){$R$}
        \end{overpic}}.
    \end{align*} 
    \item We will sometimes denote $n$ black, unoriented parallel strands by $\strand n$.
    \item We will denote the following:
     \begin{align*}
       \adjustbox{scale=1.1}{  \begin{overpic}[unit=1mm, scale=0.4]{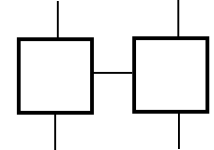}
        \put(4, 6.2){$R$}
        \put(16.5, 6.2){$R$}
        \put(11, 6.2){\tiny{$k$}}
        \put(0,-2){\tiny{$2n-k-\ell$}}
        \put(14, -2){\tiny{$2n-k-m$}}
    \end{overpic}}\quad {\raisebox{0.4cm}{{$\coloneq$}}}\hspace{2mm}
        \adjustbox{scale=1.1}{  \begin{overpic}[unit=1mm, scale=0.4]{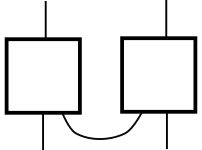}
         \put(13,-2){\tiny{$2n-k-m$}}
         \put(3, 13){\tiny{$\ell$}}
         \put(10, 2){\tiny{$k$}}
         \put(3.2, 6.2){$R$}
         \put(14, 13){\tiny{$m$}}
         \put(15.5, 6.2){$R$}
         \put(-2,-2){\tiny{$2n-k-\ell$}}
        \end{overpic}}.
    \end{align*} 
    \end{enumerate}
\end{defn}

Now we return to proofs regarding the planar algebra $\mc{P}$ defined in Theorem \ref{thm: E7sufficientmaintheorem}.

\begin{lem}\label{lem: e7sufficientclickrelation}
    In the planar algebra $\mc{P}$, $\mc{F}(S)=S$. 
\end{lem}

\begin{proof}
    This is equvalent to showing $\mc{R}(\mc{F}(S))=\mc{R}(S)$. By using the jellyfish relation, $\mc{R}(\mc{F}(S))$ is equivalent to the diagram
    \begin{align*}
        \adjustbox{scale=0.95}{\begin{overpic}[unit=1mm, scale=0.5]{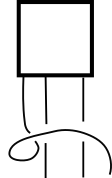}
        \put(7,-2){\tiny{$(7)$}}
        \put(7,2){\tiny{$...$}}
        \put(6,17){$S$}
    \end{overpic}}\hspace{0.1cm}.
    \end{align*}
    By Lemma \ref{lem: cuponcrossingresolved}, the leftmost crossing can be resolved to become
    \begin{align*}
        {\raisebox{1cm}{$-i$}}
        \adjustbox{scale=0.95}{\begin{overpic}[unit=1mm, scale=0.5]{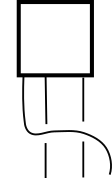}
        \put(7,-2){\tiny{$(7)$}}
        \put(7,2){\tiny{$...$}}
        \put(6,17){$S$}
    \end{overpic}}
    \end{align*}
    Since $S$ is U.C.C.S., Lemma \ref{lem: boxwithovercrossingonbottom} gives that the above is equal to
    \begin{align*}
       \adjustbox{scale=0.95}{ \begin{overpic}[unit=1mm, scale=0.5]{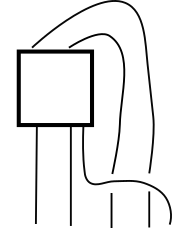}
        \put(15.5,-2){\tiny{$(4)$}}
        \put(6,2){\tiny{$...$}}
        \put(15.5,2){\tiny{$...$}}
        \put(6,17){$S$}
    \end{overpic}}\hspace{0.1cm}.
    \end{align*}
    Resolving the leftmost crossing is a linear combination of two terms, one of which is zero since it is sidecapping the $S$. The other term is  \begin{align*}
        \raisebox{1cm}{$-i$}
        \adjustbox{scale=0.95}{\begin{overpic}[unit=1mm, scale=0.5]{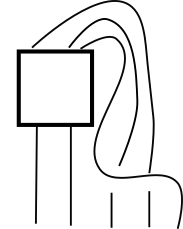}
        \put(16,-2){\tiny{$(3)$}}
        \put(6,2){\tiny{$...$}}
        \put(16,2){\tiny{$...$}}
        \put(6,17){$S$}
    \end{overpic}}.
    \end{align*}
    Using Lemma \ref{lem: boxwithovercrossingonbottom} again gives $(-i)^4\mc{R}(S)=\mc{R}(S)$. 
\end{proof}

\begin{lem}\label{lem: canuntwisttwoSboxes}
    In the planar algebra $\mc{P}$, the following equality holds for all $0 \leq c \leq 8$.
     \begin{equation}\label{eqn: twoboxesconnectedoneboxbelowonleft}\adjustbox{scale=0.97}{
       \adjustbox{scale=0.97}{{\begin{overpic}[unit=1mm, scale=0.6]{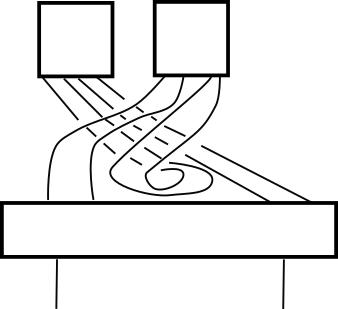}
        \put(10.5, 41.5){$S$}
        \put(29.5, 41.5){$S$}
        \put(10, 20){\tiny{$...$}}
        \put(16,11){$f^{(16-2c)}$}
        \put(20, 20){\tiny{$...$}}
        \put(19, 21){\tiny{$(c)$}}
        \put(39,20){\tiny{$...$}}
        \put(22, 2){$...$}
    \end{overpic}}}}\quad {\raisebox{0.9cm}{{$=(-1)^c$}}}\hspace{2mm}
    \adjustbox{scale=0.97}{
         \begin{overpic}[unit=1mm, scale=0.6]{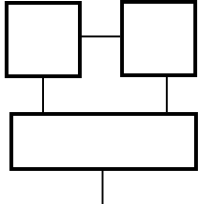}
         \put(28,17){\tiny{$8-c$}}
         \put(-2, 17){\tiny{$8-c$}}
         \put(5, 25){$S$}
         \put(15, 24){\tiny{$c$}}
         \put(24, 25){$S$}
         \put(11,8){$f^{(16-2c)}$}
        \end{overpic}}
    \end{equation}
\end{lem}

\begin{proof}
    Fix $0 \leq c \leq 8$. First use Lemma \ref{lem: cuponcrossingresolved} to resolve the innermost loop. Then, resolve the all of the crossing from the rightmost strand of the right $S$-box. It is clear that there is only one resolution of these crossings that does not cup an $S$-box or cap the $f^{(16-2n)}$. This resolution is where each crossing becomes a $B$-smoothing. Each $B$-smoothing results in a factor of $-i$. Along with the $-i$ factor from the twist, we obtain $(-i)^8=1$ times the below diagram: 
    \begin{align*}
       \adjustbox{scale=0.97}{{ \begin{overpic}[unit=1mm, scale=0.6]{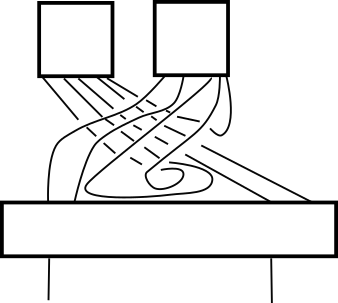}
        \put(10, 41){$S$}
        \put(29, 41){$S$}
        \put(9, 20){\tiny{$...$}}
        \put(16,11){$f^{(16-2c)}$}
        \put(18, 21){\tiny{$...$}}
        \put(15.5, 18.5){\tiny{$(c-1)$}}
        \put(37,21){\tiny{$...$}}
        \put(22, 2){$...$}
    \end{overpic}}}
    \end{align*}
    
      Continuing in this fashion, we obtain
      \begin{align*}
       \adjustbox{scale=0.97}{{\begin{overpic}[unit=1mm, scale=0.6]{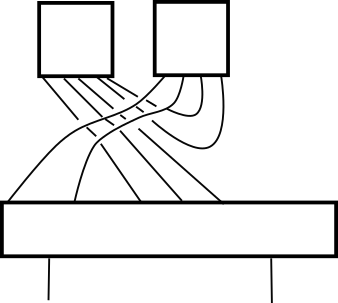}
        \put(10, 41){$S$}
        \put(29, 41){$S$}
        \put(7, 20){\tiny{$...$}}
        \put(16,11){$f^{(16-2c)}$}
        \put(20, 21){\tiny{$...$}}
        \put(30, 27){\tiny{$(c)$}}
        \put(31,26){\tiny{$...$}}
        \put(22, 2){$...$}
        \put(4, 18){\tiny{$(8-c)$}}
    \end{overpic}}}
    \end{align*}

    Now we look at the rightmost strand connecting the right $S$-box to the Jones-Wenzl projection, $f^{(16-2c)}$. Again, there is only one resolution of the crossings that will not cup an $S$-box or cap $f^{(16-2c)}$. This resolution is given by an $A$-smoothing on the lower $8-c$ crossings and a $B$-smoothing on the $c$ upper crossings. This results in the scalar $(-1)^c$ multiplied by the below diagram. 
    
    \begin{equation*}
        \adjustbox{scale=0.97}{{\begin{overpic}[unit=1mm, scale=0.6]{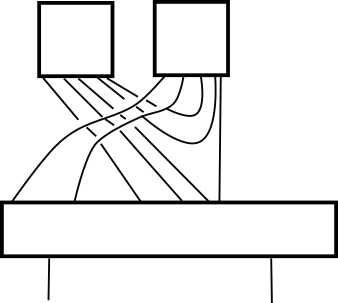}
        \put(10, 41){$S$}
        \put(29, 41){$S$}
        \put(8, 20){\tiny{$...$}}
        \put(16,11){$f^{(16-2c)}$}
        \put(20, 21){\tiny{$...$}}
        \put(28.5, 27.5){\tiny{$(c)$}}
        \put(29,26){\tiny{$...$}}
        \put(22, 2){$...$}
        \put(4.5, 18){\tiny{$(7-c)$}}
    \end{overpic}}}
    \end{equation*}
    Continuing this process, we obtain $((-1)^c)^{8-c}=(-1)^c$ times the below diagram 

    \begin{equation*}
           \adjustbox{scale=0.97}{\begin{overpic}[unit=1mm, scale=0.6]{twoboxesconnectedoneboxbelow.png}
         \put(28,17){\tiny{$8-c$}}
         \put(-2, 17){\tiny{$8-c$}}
         \put(5, 25){$S$}
         \put(15, 24){\tiny{$c$}}
         \put(24, 25){$S$}
         \put(11,8){$f^{(16-2c)}$}
        \end{overpic}},
    \end{equation*}      
        as we wished.
\end{proof}

Using Lemma \ref{lem: switchevennumberofcrossings}, the following is now immediate.
\begin{lem}
    In the planar algebra $\mc{P}$, the following equality holds for all $0\leq c \leq 8$
    \begin{align*}
        \adjustbox{scale=0.97}{{\begin{overpic}[unit=1mm, scale=0.6]{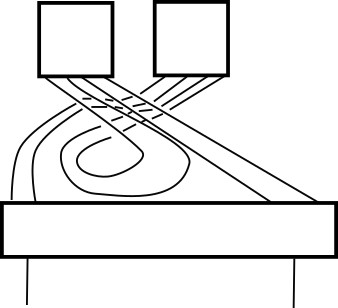}
        \put(10, 41){$S$}
        \put(29, 41){$S$}
        \put(2.5, 20){\tiny{$...$}}
        \put(16,11){$f^{(16-2c)}$}
        \put(24, 20){\tiny{$...$}}
        \put(24, 22){\tiny{$(c)$}}
        \put(40,20){\tiny{$...$}}
        \put(22, 2){$...$}
    \end{overpic}}}\quad {\raisebox{0.4cm}{{$=(-1)^c$}}}\hspace{2mm}
         \adjustbox{scale=0.97}{\begin{overpic}[unit=1mm, scale=0.6]{twoboxesconnectedoneboxbelow.png}
         \put(28,17){\tiny{$8-c$}}
         \put(-2, 17){\tiny{$8-c$}}
         \put(5, 25){$S$}
         \put(15, 24){\tiny{$c$}}
         \put(23, 25){$S$}
         \put(11,8){$f^{(16-2c)}$}
        \end{overpic}}
    \end{align*}
\end{lem}

\begin{lem}\label{lem: SSf10iszero}
    In $\mc{P}$, the below diagram is zero for $c=1$ and $3$.
     \begin{equation}\label{eqn: SSf10}
        \adjustbox{scale=0.97}{ \begin{overpic}[unit=1mm, scale=0.6]{twoboxesconnectedoneboxbelow.png}
         \put(28,17){\tiny{$8-c$}}
         \put(-2, 17){\tiny{$8-c$}}
         \put(5, 25){$S$}
         \put(15, 24){\tiny{$c$}}
         \put(23, 25){$S$}
         \put(11,8){$f^{(16-2c)}$}
        \end{overpic}}
    \end{equation}
\end{lem}

\begin{proof}
    The proof is nearly identical for $c=1$ and $3$. We just show the $c=3$ case. By Lemma \ref{lem: canuntwisttwoSboxes}, we have
     \begin{align*}
       {\adjustbox{scale=0.97}{ {%
        \begin{overpic}[unit=1mm, scale=0.6]{twist1.png}
        \put(10, 41){$S$}
        \put(29, 41){$S$}
        \put(10, 20){\tiny{$...$}}
        \put(16,11){$f^{(10)}$}
        \put(20, 20){\tiny{$...$}}
        \put(19, 21){\tiny{$(3)$}}
        \put(39,20){\tiny{$...$}}
        \put(22, 2){$...$}
    \end{overpic}
    }}}
    \quad {\raisebox{0.6cm}{{$=-1$}}}\hspace{2mm}
         {\adjustbox{scale=0.97}{ {%
          \begin{overpic}[unit=1mm, scale=0.6]{twoboxesconnectedoneboxbelow.png}
         \put(28,17){\tiny{$5$}}
         \put(0, 17){\tiny{$5$}}
         \put(5, 25){$S$}
         \put(15, 24){\tiny{$3$}}
         \put(23, 25){$S$}
         \put(11,8){$f^{(10)}$}
        \end{overpic}
        }}}
    \end{align*}
    On the other hand, using the jellyfish relation, the right-hand $S$-box can be braided under the strands of the left-hand $S$-box to obtain
    \begin{align*}
       {\adjustbox{scale=0.97}{ {%
        \begin{overpic}[unit=1mm, scale=0.6]{twist1.png}
        \put(10, 41){$S$}
        \put(29, 41){$S$}
        \put(10, 20){\tiny{$...$}}
        \put(16,11){$f^{(10)}$}
        \put(20, 20){\tiny{$...$}}
        \put(19, 21){\tiny{$(3)$}}
        \put(39,20){\tiny{$...$}}
        \put(22, 2){$...$}
    \end{overpic}
    }}}
    \quad {\raisebox{0.4cm}{{$=$}}}\hspace{2mm}
         {\adjustbox{scale=0.97}{ {%
          \begin{overpic}[unit=1mm, scale=0.6]{twoboxesconnectedoneboxbelow.png}
         \put(28,17){\tiny{$5$}}
         \put(0, 17){\tiny{$5$}}
         \put(5, 25){$S$}
         \put(15, 24){\tiny{$3$}}
         \put(23, 25){$S$}
         \put(11,8){$f^{(10)}$}
        \end{overpic}
        }}}
    \end{align*}
    which gives the result.
\end{proof}

\begin{lem}\label{lem: SSSf14one}
    In the planar algebra $\mc{P}$, the below diagram is zero.
     \begin{equation}\label{eq: 3boxesand1belowstrandsonbottom}
       \adjustbox{scale=0.97}{  \begin{overpic}[unit=1mm, scale=0.6]{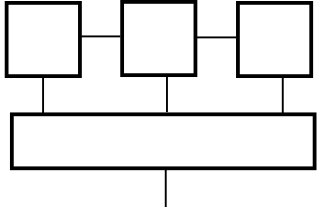}
         \put(24,26){$S$}
         \put(5, 26){$S$}
         \put(41, 26){$S$}
         \put(21, 18){\tiny{$3$}}
         \put(41, 18){\tiny{$6$}}
         \put(23, 9){$f^{(14)}$}
          \put(3, 18){\tiny{$5$}}
            \put(34, 28){\tiny{$2$}}
            \put(15, 28){\tiny{$3$}}
        \end{overpic}}
    \end{equation}  
\end{lem}

\begin{proof}
    Consider the diagram from Lemma \ref{lem: SSf10iszero} for $c=3$. On the one hand, the diagram is zero. On the other hand, we can expand the $f^{(10)}$ to get a linear combination of diagrams giving 
    \begin{equation}\label{eq: S3S}
      \raisebox{1cm}{$0=$} \raisebox{0.5cm}{ \adjustbox{scale=0.97}{  \begin{overpic}[unit=1mm, scale=0.6]{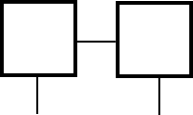}
         \put(23, 10){$S$}
         \put(15, 8){\tiny{$3$}}
         \put(5, 10){$S$}
        \end{overpic} }}\hspace{2mm}
    \raisebox{1cm}{$+ \mathlarger{\mathlarger{\sum_i}} a_i$}\hspace{2mm}
    \adjustbox{scale=0.97}{ \begin{overpic}[unit=1mm, scale=0.6]{twoboxesconnectedoneboxbelow.png}
         \put(28,17){\tiny{$5$}}
         \put(4, 17){\tiny{$5$}}
         \put(5, 25){$S$}
         \put(15, 24){\tiny{$3$}}
         \put(23, 25){$S$}
         \put(13,8){$T_i$}
        \end{overpic}}
    \end{equation}
    where $T_i$ are in Temperley-Lieb with at least one cup, and all $a_i \in \mbb{C}$. Each term in the sum that cups either $S$ becomes zero. The nonzero terms will be the ones where $T_i$ has only cups in the middle, forming an $S^2$. Say $T_i$ has $d_i$ cups. Therefore, 
    \begin{equation}\label{eqn: writeS3SaslinearcombooflessS}
\raisebox{0.5cm}{ \adjustbox{scale=0.97}{  \begin{overpic}[unit=1mm, scale=0.6]{twoboxesconnecteddown.png}
         \put(23, 10){$S$}
         \put(15, 8){\tiny{$3$}}
         \put(5, 10){$S$}
        \end{overpic} }}\hspace{2mm}
    \raisebox{1cm}{$=- \mathlarger{\mathlarger{\sum_i}} a_i$}\hspace{2mm}
     \adjustbox{scale=0.97}{  \begin{overpic}[unit=1mm, scale=0.6]{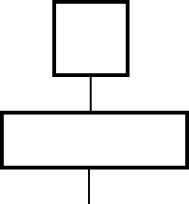}
          \put(13, 25){$S$}
         \put(17, 17){\tiny{$8$}}
         \put(13, 9){$R_i$}
        \end{overpic}}
         \hspace{2mm}\raisebox{1cm}{$-6 \mathlarger{\mathlarger{\sum_i}} a_i$}\hspace{2mm}
 \adjustbox{scale=0.97}{ \begin{overpic}[unit=1mm, scale=0.6]{twoboxesstackedwithbottomstrand.png}
          \put(12, 25){$f^{(4)}$}
         \put(17, 17){\tiny{$8$}}
         \put(13, 9){$R_i$}
        \end{overpic}}
    \end{equation}
    where $R_i$ are the same as $T_i$ except with $d_i-1$ cups in the middle. Note that $R_i$ has $8$ strands on top and $10$ strands on bottom.

        Use this formula to replace the two leftmost $S$-boxes in the diagram (\ref{eq: 3boxesand1belowstrandsonbottom}). This gives that the diagram (\ref{eq: 3boxesand1belowstrandsonbottom}) equals
        
        \begin{equation}\label{eq: 211expansion}
        \scalebox{0.85}{\raisebox{2.5cm}{$-\mathlarger{\mathlarger{\sum_i}} a_i$}\hspace{2mm}
        \adjustbox{scale=0.97}{  \begin{overpic}[unit=1mm, scale=0.6]{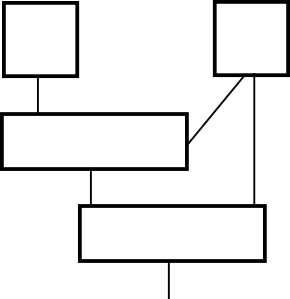}
         \put(31,32){\tiny{$2$}}
         \put(1, 32){\tiny{$8$}}
         \put(5, 40){$S$}
         \put(38, 24){\tiny{$6$}}
         \put(39, 40){$S$}
         \put(13,23){$R_i$}
         \put(16, 17){\tiny{$8$}}
         \put(22, 9){$f^{(14)}$}
        \end{overpic}}
        \hspace{2mm}\raisebox{2.5cm}{$-6 \mathlarger{\mathlarger{\sum_i}} a_i$}\hspace{2mm}
        \adjustbox{scale=0.97}{  \begin{overpic}[unit=1mm, scale=0.6]{211.png}
         \put(31,32){\tiny{$2$}}
         \put(1, 32){\tiny{$8$}}
         \put(3, 40){$f^{(4)}$}
         \put(38, 24){\tiny{$6$}}
         \put(39, 40){$S$}
         \put(13,23){$R_i$}
         \put(16, 17){\tiny{$8$}}
         \put(22, 9){$f^{(14)}$}
        \end{overpic}}}
        \end{equation}

        Recall that we assumed $T_i$ was not the identity, so $T_i$ and thus $R_i$, has some cap. Since $f^{(14)}$ and $S$ are uncappable, there is a cap between strands 8 and 9. In this case, there are two other options of $R_i$, the rest of the strands are straight, or there is a cap between strands 7 and 10. 

        Consider the terms in the left sum of equation \ref{eq: 211expansion}. If $R_i$ has a cap between strands 8 and 9 and otherwise has no caps then we get the term 
        \begin{align*}
          \adjustbox{scale=0.97}{   \begin{overpic}[unit=1mm, scale=0.6]{twoboxesconnectedoneboxbelow.png}
         \put(28,17){\tiny{$7$}}
         \put(0, 17){\tiny{$7$}}
         \put(5.5, 25){$S$}
         \put(15, 24){\tiny{$1$}}
         \put(23.5, 25){$S$}
         \put(11.5,9){$f^{(14)}$}
        \end{overpic}}\hspace{1mm},
        \end{align*}
        which by Lemma \ref{lem: SSf10iszero} is zero. If $R_i$ has another cap between strands 7 and 10 then there must be an additional cup on top, which cups $S$ and is zero. The first summand is thus zero.

        Consider the terms in the second sum of equation \ref{eq: 211expansion}. If $R_i$ has a cap between strands 8 and 9 and otherwise no caps, then we get the term
         \begin{align*}
           \adjustbox{scale=0.97}{  \begin{overpic}[unit=1mm, scale=0.6]{twoboxesconnectedoneboxbelow.png}
         \put(28,17){\tiny{$7$}}
         \put(0, 17){\tiny{$7$}}
         \put(4, 25){$f^{(4)}$}
         \put(15, 24){\tiny{$1$}}
         \put(23.5, 25){$S$}
         \put(11.5,9){$f^{(14)}$}
        \end{overpic}}\hspace{1mm}.
        \end{align*}
        Expanding $f^{(4)}$, we can see that all the terms will cap $f^{(14)}$, so this diagram is zero. If $R$ also has a cap between strands 7 and 10, then it must have another cup. The only place a nonzero cup could happen would be between strands 4 and 5, which would give a partial trace of $f^{(4)}$. Recall that $E_4(f^{(4)})=\f{5}{4}f^{(3)}$. The $f^{(3)}$ will still cap $f^{(14)}$, so this term is also zero.
\end{proof}

A symmetric argument gives the next lemma. 

\begin{lem}\label{lem: SSSf14two}
     In the planar algebra $\mc{P}$, the below diagram is zero.
      \begin{align*}
       \adjustbox{scale=0.97}{ \begin{overpic}[unit=1mm, scale=0.6]{3boxesand1belowstrandsonbottom.png}
         \put(24,26){$S$}
         \put(5, 26){$S$}
         \put(41, 26){$S$}
         \put(21, 18){\tiny{$3$}}
         \put(41, 18){\tiny{$5$}}
         \put(23, 9){$f^{(14)}$}
          \put(3, 18){\tiny{$6$}}
            \put(34, 28){\tiny{$3$}}
            \put(15, 28){\tiny{$2$}}
        \end{overpic}}
    \end{align*}  
\end{lem}

\begin{lem}\label{lem: SSSSf14}
    In the planar algebra, the below diagram is zero.  
    \begin{equation}\label{eqn: SSSSf14}
       \adjustbox{scale=0.97}{  \begin{overpic}[unit=1mm, scale=0.6]{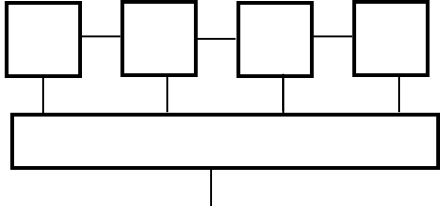}
         \put(5, 25){$S$}
         \put(23,25){$S$}
         \put(42, 25){$S$}
         \put(60,25){$S$}
            \put(34, 9){$f^{(14)}$}
           \put(3, 17){\tiny{$5$}}
           \put(15, 28){\tiny{$3$}}
         \put(22, 17){\tiny{$2$}}
         \put(52, 28){\tiny{$3$}}
         \put(42, 17){\tiny{$2$}}  
         \put(60,18){\tiny{$5$}} 
         \put(32.5, 28){\tiny{$3$}}
        \end{overpic}}
    \end{equation}  
\end{lem}

\begin{proof}
    Using equation \ref{eqn: writeS3SaslinearcombooflessS} on the two left $S$-boxes gives that diagram in (\ref{eqn: SSSSf14}) equals
    \begin{align*}
    \scalebox{0.85}{\raisebox{1cm}{$\mathlarger{\mathlarger{\sum_i}}a_i$}\hspace{2mm}
        \adjustbox{scale=0.97}{\begin{overpic}[unit=1mm, scale=0.6]{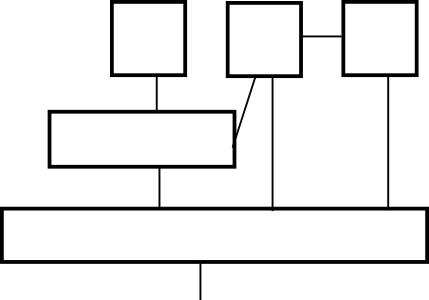}
         \put(20, 25){$R_i$}
         \put(23,40){$S$}
         \put(41, 40){$S$}
         \put(60,40){$S$}
            \put(32, 9){$f^{(14)}$}
           \put(37, 32){\tiny{$3$}}
         \put(22, 32){\tiny{$8$}}
         \put(51, 43){\tiny{$3$}}
         \put(41, 30){\tiny{$2$}}  
        \put(59,30){\tiny{$5$}} 
         \put(22, 17){\tiny{$7$}}
        \end{overpic}}
        \hspace{2mm}\raisebox{1cm}{$-6\mathlarger{\mathlarger{\sum_i}}a_i$}\hspace{2mm}
        \adjustbox{scale=0.97}{ \begin{overpic}[unit=1mm, scale=0.6]{311.png}
         \put(20, 25){$R_i$}
         \put(21,40){$f^{(4)}$}
         \put(41, 40){$S$}
         \put(60,40){$S$}
            \put(32, 9){$f^{(14)}$}
           \put(37, 32){\tiny{$3$}}
         \put(22, 32){\tiny{$8$}}
         \put(51, 43){\tiny{$3$}}
         \put(41, 30){\tiny{$2$}}  
         \put(59,30){\tiny{$5$}} 
         \put(22, 17){\tiny{$7$}}
        \end{overpic}}
        }
    \end{align*}  
    where $R_i$ has 8 strands on top and 10 strands on bottom and $a_i\in \mathbb{C}$. $R_i$ must have some cap.
    
    Consider a diagram in the left sum. Since $S$ is U.C.C.S., there must be a cap between strands 7 and 8. If the diagram is nonzero, there can be no other caps as then there will be another cup, which will cup $S$ resulting in zero. Then the rest of the strands are straight and the resulting diagram is 
    \begin{align*}
        \adjustbox{scale=0.97}{\begin{overpic}[unit=1mm, scale=0.6]{3boxesand1belowstrandsonbottom.png}
         \put(24,26){$S$}
         \put(5, 26){$S$}
         \put(43, 26){$S$}
         \put(21, 18){\tiny{$3$}}
         \put(41, 18){\tiny{$5$}}
         \put(23, 9){$f^{(14)}$}
          \put(3, 18){\tiny{$6$}}
            \put(34, 28){\tiny{$3$}}
            \put(15, 28){\tiny{$2$}}
        \end{overpic}}
    \end{align*}  
    which we already know is zero. 

    Consider a diagram in the right sum. Again, there must be a cap between strands 7 and 8. If the rest of the strands are straight then we obtain the following diagram: 

     \begin{align*}
        \adjustbox{scale=0.97}{\begin{overpic}[unit=1mm, scale=0.6]{3boxesand1belowstrandsonbottom.png}
         \put(24,26){$S$}
         \put(5, 26){$f^{(4)}$}
         \put(42, 26){$S$}
         \put(21, 18){\tiny{$3$}}
         \put(41, 18){\tiny{$5$}}
         \put(23, 9){$f^{(14)}$}
          \put(3, 18){\tiny{$6$}}
            \put(34, 28){\tiny{$3$}}
            \put(15, 28){\tiny{$2$}}
        \end{overpic}}.
    \end{align*}  

    Expanding $f^{(4)}$ will cap $f^{(14)}$, so this term is zero. If $R_i$ had another cap, which, in order to be nonzero, would need to be between strands 6 and 7, then there must be a cup between strands 4 and 5. This gives a partial trace of $f^{(4)}$, which is a scalar multiple of $f^{(3)}$. Expanding $f^{(3)}$ will cap $f^{(14)}$, giving zero.
\end{proof}

\begin{lem}\label{lem: twoboxesconnectedoneboxbelow}
    For $4<c<8$, the below diagram is zero.
    \begin{equation}
         \adjustbox{scale=0.97}{\begin{overpic}[unit=1mm, scale=0.6]{twoboxesconnectedoneboxbelow.png}
         \put(28,17){\tiny{$8-c$}}
         \put(-2, 17){\tiny{$8-c$}}
         \put(5, 25){$S$}
         \put(15, 24){\tiny{$c$}}
         \put(24, 25){$S$}
         \put(11,8){$f^{(16-2c)}$}
        \end{overpic}}
    \end{equation}
\end{lem}

\begin{proof}
    Fix $4<c<8$. Since $c-4>0$, the two $S$ boxes form an $S^2$-box with some cup, cap, or side cap. We can then expand this $S^2$-box by using that $S^2=6f^{(4)}+S$. Since $S$ is U.C.C.S., the second term is zero. On the other hand, the first term is 6 times the following diagram (where the $f^{(4)}$ is rainbowed with all endpoints on the bottom):
    \begin{align*}
        \adjustbox{scale=0.97}{ \begin{overpic}[unit=1mm, scale=0.4]{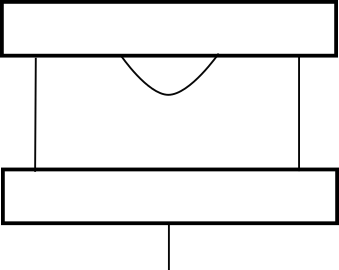}
         \put(33,17){\tiny{$8-c$}}
         \put(15, 16){\tiny{$c-4$}}
         \put(15, 24.5){$f^{(4)}$}
         \put(-3, 17){\tiny{$8-c$}}
         \put(11,6.5){$f^{(16-2c)}$}
        \end{overpic}}
    \end{align*}
    When expanding $f^{(4)}$ , every term will cap $f^{(16-2c)}$ and thus will also be zero.
\end{proof}

\subsection{The Leaf Relations}
In this section, we prove the necessary leaf relations for affine $E_7$ hold in $\mathcal{P}$.

\begin{lem}
    $\mc{R}(Q_4 \otimes X)f^{(10)}-2\mc{R}((Q_4 \otimes X)e_4(Q_4 \otimes X))f^{(10)}=0$.
\end{lem}

\begin{proof}
    Since $f^{(10)}$ is uncappable, $\mc{R}(Q_4 \otimes X)f^{(10)}=0$. Recall that $Q_4=\f{2}{5}f^{(4)}+\f{1}{5}S$, so expanding the second term gives $ \mc{R}((Q_4 \otimes X)e_4(Q_4 \otimes X))f^{(10)}$ equals
    \begin{align*}
       &\f{4}{25}\mc{R}((f^{(4)}\otimes X)e_4(f^{(4)}\otimes X))f^{(10)}+\f{2}{25}\mc{R}((f^{(4)}\otimes X)e_4(S\otimes X))f^{(10)}\\
       &+\f{2}{25}\mc{R}((S\otimes X)e_4(f^{(4)}\otimes X))f^{(10)}+\f{1}{25}\mc{R}((S\otimes X)e_4(S\otimes X))f^{(10)}
    \end{align*}
    By Lemma \ref{lem: SSf10iszero}, the last term is zero. The other terms involve an $f^{(4)}$ on top of an $f^{(10)}$ which will cap it and become zero. 
    
\end{proof}

\begin{lem}
    In $\mc{P}$, $\mc{R}((Q_4 \otimes X)e_4(Q_4 \otimes X))$ is uncuppable and uncappable.
\end{lem}

\begin{proof}
    There are no strands on the top, so $\mc{R}((Q_4 \otimes X)e_4(Q_4 \otimes X))$ is obviously uncappable. Putting a cup between strands 1 and 2, 2 and 3, 3 and 4, 7 and 8, 8 and 9, and 9 and 10 are all zero since they cup $Q_4$. 
    
    Recall from Lemma \ref{lem: e7sufficientpartialtraces} that $E_4(Q_4)=\f{1}{2}f^{(3)}$. Putting a cup between strands 4 and 5 gives $\mc{R}(Q_4)-2\cdot \f{1}{2}\mc{R}(Q_4(f^{(3)}\otimes X))=\mc{R}((Q_4)-\mc{R}(Q_4))=0$, since $Q_4$ is uncuppable and uncappable. Nearly identically, a cup between strands 6 and 7 also gives zero. Putting a cup between strands 5 and 6 gives $2\mc{R}(Q_4)-2\mc{R}((Q_4)^2)=0$. Thus  $\mc{R}((Q_4 \otimes X)e_4(Q_4 \otimes X))$ is also uncuppable. 
\end{proof}

From this we now obtain a leaf relation.

\begin{lem}\label{lem: e7sufficientleafrelationq4}
    In $\mc{P}$, $(Q_4 \otimes X)-2(Q_4 \otimes X)e_4(Q_4\otimes X)=0$. 
\end{lem}

\begin{proof}
    Since $\mc{R}((Q_4 \otimes X)-2(Q_4 \otimes X)e_4(Q_4\otimes X))$ is uncuppable and uncappable, $\mc{R}((Q_4 \otimes X)-2(Q_4 \otimes X)e_4(Q_4\otimes X))=\mc{R}((Q_4 \otimes X)-2(Q_4 \otimes X)e_4(Q_4\otimes X))f^{(10)}=0$. 
\end{proof}

\begin{lem}
    $\mc{R}((P_6\otimes X))-2\mc{R}((P_6\otimes X)e_6(P_6\otimes X))$ is uncuppable and uncappable.
\end{lem}

\begin{proof}
    There are no strands at the top of the diagram so $\mc{R}((P_6\otimes X))-2\mc{R}((P_6\otimes X)e_6(P_6\otimes X))$ is uncappable. Since $P_6$ is uncuppable, putting a cup between strands 1 and 2, 2 and 3, 3 and 4, 4 and 5, 5 and 6, 9 and 10, 10 and 11, 11 and 12, 12 and 13, and 13 and 14 gives zero. 

    Recall that $E_6(P_6)=\f{1}{2}P_5$. Putting a cup between strands 6 and 7 gives $\mc{R}(P_6)-2\f{1}{2}\mc{R}(P_6(P_5 \otimes X))=\mc{R}(P_6)-\mc{R}(P_6)=0$. An identical argument works for putting a cup between strands 8 and 9. Putting a cup between strands 7 and 8 gives $2\mc{R}(P_6)-2\mc{R}((P_6)^2)=0$. Therefore, $\mc{R}((P_6\otimes X))-2\mc{R}((P_6\otimes X)e_6(P_6\otimes X))$ is uncuppable.
\end{proof}

\begin{lem}
    $\mc{R}((P_6\otimes X))f^{(14)}-2\mc{R}((P_6\otimes X)e_6(P_6\otimes X))f^{(14)}=0$.
\end{lem}

\begin{proof}
    Notice that $\mc{R}(P_6 \otimes X)f^{(14)}=0$ since the rainbowed $X$ caps $f^{(14)}$. Recall that $P_6=P_5 \otimes X - \f{3}{2}(P_5\otimes X)e_5(P_5 \otimes X)$. Using this, we expand the bottom $P_6\otimes X$ to get that $\mc{R}((P_6\otimes X)e_6(P_6\otimes X))f^{(14)}$ equals
    \begin{align*}
        \mc{R}((P_6\otimes X)e_6(P_5\otimes X^{\otimes 2}))f^{(14)}-\f{3}{2}\mc{R}((P_6\otimes X)e_6((P_5 \otimes X^{\otimes 2})e_5(P_5\otimes X^{\otimes 2}))
    \end{align*}
    The first term above is zero since $f^{(14)}$ will be capped. 
    
    We now focus on computing the second term. We will first expand both $P_6$ in $(P_6\otimes X)e_6(P_6\otimes X)$ down to a linear combination of $S$ and $f^{(4)}$. Along the way, we will cancel any term that, when rainbowed, will cap $f^{(14)}$. We will use the recursive absorption rule several times. Then $(P_6\otimes X)e_6(P_6\otimes X)$ equals
    \begin{align*}
        (P_6\otimes X)(X^{\otimes 5}\otimes \te{coev}_X)P_5(X^{\otimes 5}\otimes \te{ev}_X)-\f{3}{2}(P_6 \otimes X)(X^{\otimes 5} \otimes \te{coev}_X)P_5(X^{\otimes 4}\otimes \te{ev}_X \otimes X)(P_5 \otimes X^{\otimes 2})
    \end{align*}
    The first term will cap $f^{(14)}$ so we just need to continue expanding the second term. The $P_6$ can absorb the top $P_5$. So the second term equals $-\f{3}{2}(P_6 \otimes X)e_6e_5(P_5 \otimes X^{\otimes 2})$. Expanding $P_5$ gives that $-\f{3}{2}(P_6 \otimes X)e_6e_5(P_5 \otimes X^{\otimes 2})$ equals
    \begin{align*}
        -\f{3}{2}(P_6 \otimes X)e_6e_5(P_4 \otimes X^{\otimes 3})+2(P_6 \otimes X)e_6e_5(P_4 \otimes X^{\otimes 3})e_4(P_4 \otimes X^{\otimes 3})
    \end{align*}
    The first term will cap $f^{(14)}$ when rainbowed over. When $P_6$ absorbs the $P_4$ we obtain $2(P_6 \otimes X)e_6e_5e_4(P_4 \otimes X^{\otimes 3})$ for the second term. Using that $P_4=\f{3}{5}f^{(4)}-\f{1}{5}S$, we get that $2(P_6 \otimes X)e_6e_5e_4(P_4 \otimes X^{\otimes 3})$ equals
    \begin{align*}
        \f{6}{5}(P_6 \otimes X)e_6e_5e_4(f^{(4)} \otimes X^{\otimes 3})-\f{2}{5}(P_6 \otimes X)e_6e_5e_4(S \otimes X^{\otimes 3})
    \end{align*}
    The first term, when rainbowed, will cap $f^{(14)}$ when expanding the $f^{(4)}$. For the second term, we expand $P_6$ to get that $-\f{2}{5}(P_6 \otimes X)e_6e_5e_4(S \otimes X^{\otimes 3})$ equals
    \begin{align*}
        -\f{2}{5}(P_5 \otimes X^{\otimes 2})e_6e_5e_4(S \otimes X^{\otimes 3})+\f{3}{5}(P_5 \otimes X^{\otimes 2})e_5(P_5 \otimes X^{\otimes 2})e_6e_5e_4(S \otimes X^{\otimes 3})
    \end{align*}
    When rainbowed on top of $f^{(14)}$, the first term will be zero. Then the second term can be expanded to equal
    \begin{align*}
        \f{3}{5}(P_4 \otimes X^{\otimes 3})e_5(P_5 \otimes X^{\otimes 2})e_6e_5e_4(S \otimes X^{\otimes 3})-\f{4}{5}(P_4 \otimes X^{\otimes 3})e_4(P_4 \otimes X^{\otimes 3})e_5(P_5 \otimes X^{\otimes 2})e_6e_5e_4(S \otimes X^{\otimes 3})
    \end{align*}
    The first term will cap $f^{(14)}$ when rainbowed onto a $f^{(14)}$. The $P_5$ in the second term can absorb $P_4$ which will then be equal to $-\f{4}{5}(P_4\otimes X^{\otimes 3})e_4e_5e_6(P_5 \otimes X^{\otimes 2})e_5e_4(S\otimes X^{\otimes 3})$. Expanding $P_4$ gives that $-\f{4}{5}(P_4\otimes X^{\otimes 3})e_4e_5e_6(P_5 \otimes X^{\otimes 2})e_5e_4(S\otimes X^{\otimes 3})$ equals
    \begin{align*}
        -\f{12}{25}(f^{(4)}\otimes X^{\otimes 3})e_4e_5e_6(P_5 \otimes X^{\otimes 2})e_5e_4(S\otimes X^{\otimes 3}+\f{4}{25}(S\otimes X^{\otimes 3})e_4e_5e_6(P_5 \otimes X^{\otimes 2})e_5e_4(S\otimes X^{\otimes 3})
    \end{align*}
    The first term, when expanding the $f^{(4)}$, will cap the $f^{(14)}$ when rainbowed. For the second term, expanding $P_5$ gives that $\f{4}{25}(S\otimes X^{\otimes 3})e_4e_5e_6(P_5 \otimes X^{\otimes 2})e_5e_4(S\otimes X^{\otimes 3})$ equals
    \begin{align*}
        &\f{4}{25}(S\otimes X^{\otimes 3})e_4e_5e_6(P_4 \otimes X^{\otimes 3})e_5e_4(S\otimes X^{\otimes 3})-\f{16}{75}(S\otimes X^{\otimes 3})e_4e_5e_6(P_4 \otimes X^{\otimes 3})e_4(P_4 \otimes X^{\otimes 3})e_5e_4(S\otimes X^{\otimes 3})
    \end{align*}
    The first term, when rainbowed, will cap $f^{(14)}$ and thus be zero. The second term can have each of the $P_4$ expanded. This gives
    \begin{align*}
        &-\f{48}{625}(S\otimes X^{\otimes 3})e_4e_5e_6(f^{(4)} \otimes X^{\otimes 3})e_4(f^{(4)} \otimes X^{\otimes 3})e_5e_4(S\otimes X^{\otimes 3})+\f{16}{625}(S\otimes X^{\otimes 3})e_4e_5e_6(S \otimes X^{\otimes 3})e_4(f^{(4)} \otimes X^{\otimes 3})e_5e_4(S\otimes X^{\otimes 3})\\
        &+\f{16}{625}(S\otimes X^{\otimes 3})e_4e_5e_6(f^{(4)} \otimes X^{\otimes 3})e_4(S \otimes X^{\otimes 3})e_5e_4(S\otimes X^{\otimes 3})-\f{16}{1875}(S\otimes X^{\otimes 3})e_4e_5e_6(S \otimes X^{\otimes 3})e_4(S \otimes X^{\otimes 3})e_5e_4(S\otimes X^{\otimes 3})
    \end{align*}
    Consider each of these terms rainbowed over an $f^{(14)}$. Then these, up to scalar multiple, are of the form 
    \begin{align*}
        \adjustbox{scale=0.97}{\begin{overpic}[unit=1mm, scale=0.6]{41strandsdown.png}
         \put(5, 25){$S$}
         \put(23,25){$R$}
         \put(42, 25){$T$}
         \put(61,25){$S$}
            \put(34, 9){$f^{(14)}$}
           \put(3, 17){\tiny{$5$}}
           \put(15, 28){\tiny{$3$}}
         \put(22, 17){\tiny{$2$}}
         \put(52, 28){\tiny{$3$}}
         \put(42, 17){\tiny{$2$}}  
         \put(60,18){\tiny{$5$}} 
         \put(32.5, 28){\tiny{$3$}}
        \end{overpic}}
    \end{align*}
    where $R,T \in \{f^{(4)}, S\}$. When both $R$ and $T$ are $S$, then by Lemma \ref{lem: SSSSf14}, the diagram is zero. When both $R$ and $T$ are $f^{(4)}$, we can see by expanding $f^{(4)}$ that the only possible term not capping $f^{(14)}$ has the form with $c=1$ of Lemma \ref{eqn: SSf10}, which is zero. If exactly one of $R$ or $S$ is $f^{(4)}$ then the only options of diagrams not capping $f^{(14)}$ or $S$ must be have a cap between strands 3 and 4 and a cup between strands 3 and 4, which will be a scalar multiple of the diagram from either Lemma \ref{lem: SSSf14one} or \ref{lem: SSSf14two}, which are both zero. Therefore, $\mc{R}((P_6\otimes X)e_6(P_6 \otimes X))f^{(14)}=0$. 
\end{proof}

\begin{lem}\label{lem: e7sufficientleafrelationp6}
    In $\mc{P}$, $(P_6 \otimes X)-2(P_6\otimes X)e_6(P_6 \otimes X)=0$. 
\end{lem}

\begin{proof}
    We have that $\mc{R}((P_6 \otimes X)-2(P_6\otimes X)e_6(P_6 \otimes X))f^{(14)}=0$ and $\mc{R}((P_6 \otimes X)-2(P_6\otimes X)e_6(P_6 \otimes X))$ is uncuppable and uncappable, which gives the result. 
\end{proof}

\subsection{Closed Diagram Computations in the Planar Algebra}\label{sec: closed-diagrams}

Next, we compute a variety of closed diagrams in $\mc{P}$ that will become beneficial to know later on.

\begin{lem}\label{lem: e7sufficienttraceofSsquared}
    In $\mc{P}$, $\te{tr}(S^2)=30$. 
\end{lem}

\begin{proof}
    Since $S^2=S+6f^{(4)}$ and $S$ is U.C.C.S., $\te{tr}(S^2)=6\te{tr}(f^{(4)})=6\cdot 5 =30$.
\end{proof}

\begin{lem}\label{lem: fourboxesconnected}
    The following equality holds in $\mc{P}$. 
    \begin{align*}
        \adjustbox{scale=0.95}{\begin{overpic}[unit=1mm, scale=0.6]{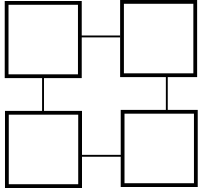}
         \put(24,5){$S$}
         \put(5, 5){$S$}
         \put(5, 23){$S$}
         \put(16, 25){\tiny{$3$}}
         \put(24, 23){$S$}
         \put(3,14){\tiny{$5$}}
         \put(27,14){\tiny{$5$}}
         \put(16, 1){\tiny{$3$}}
        \end{overpic}}
       \hspace{2mm}{\raisebox{1.2cm}{$=225$}} 
    \end{align*}  
\end{lem}

\begin{proof}
    First reduce the left set of vertical $S$-boxes connected by 5 strands using that $S^2=S+6f^{(4)}$. The term with an $S$-box will be zero since $S$ is U.C.C.S. The other term will have a partial trace of $f^{(4)}$, which was calculated in Lemma \ref{lem: e7sufficientpartialtraces} to be $\f{5}{4}f^{(3)}$. As $S$ is U.C.C.S., $S$ absorbs the $f^{(3)}$, resulting in $\f{5}{4}\cdot 6 \cdot \te{tr}(S^2)=225$. 
\end{proof}

\begin{lem}\label{lem: fourboxesconnectedequaling300}
    The following equality holds in $\mc{P}$. 
    \begin{equation}\label{eq: fourboxesconnectedequaling300}
        \adjustbox{scale=0.95}{\begin{overpic}[unit=1mm, scale=0.6]{fourboxesconnected.png}
         \put(24,5){$S$}
         \put(5, 5){$S$}
         \put(5, 23){$S$}
         \put(16, 25){\tiny{$2$}}
         \put(24, 23){$S$}
         \put(3,14){\tiny{$6$}}
         \put(27,14){\tiny{$6$}}
         \put(16, 1){\tiny{$2$}}
        \end{overpic}}
       \hspace{2mm}{\raisebox{1.2cm}{$=300$}} 
    \end{equation}  
\end{lem}

\begin{proof}
    First reduce the left set of vertical $S$-boxes using the $S^2$ relation. The term with $S$ will become zero as $S$ is U.C.C.S. The other term will have the leftmost box as $E_3(E_4(f^{(4)}))=\f{5}{4}\cdot \f{4}{3}f^{(2)}=\f{5}{3}f^{(2)}$ by Lemma \ref{lem: e7sufficientpartialtraces}. Since $S$ is U.C.C.S., diagram (\ref{eq: fourboxesconnectedequaling300}) will equal $6\cdot \f{5}{3}\te{tr}(S^2)=300$. 
\end{proof}

\begin{lem}\label{lem: fourboxesconnectedequaling450}
    The following equality holds in $\mc{P}$. 
    \begin{equation}\label{eq: fourboxesconnectedequaling450}
        \adjustbox{scale=0.95}{{\begin{overpic}[unit=1mm, scale=0.6]{fourboxesconnected.png}
         \put(24,5){$S$}
         \put(5, 5){$S$}
         \put(5, 23){$S$}
         \put(16, 25){\tiny{$1$}}
         \put(24, 23){$S$}
         \put(3,14){\tiny{$7$}}
         \put(27,14){\tiny{$7$}}
         \put(16, 1){\tiny{$1$}}
        \end{overpic}}}
       \hspace{2mm}{\raisebox{1.2cm}{$=450$}} 
    \end{equation}  
\end{lem}

\begin{proof}
    First reduce the left set of vertical $S$-boxes using the $S^2$ relation. The term with $S$ will become zero since $S$ is U.C.C.S. The other term will have the leftmost box as $E_2(E_3(E_4(f^{(4)})))=\f{5}{4}\cdot \f{4}{3}\cdot \f{3}{2}X=\f{5}{2}X$ by Lemma \ref{lem: e7sufficientpartialtraces}. Then diagram (\ref{eq: fourboxesconnectedequaling450}) will equal $6\cdot \f{5}{2}\te{tr}(S^2)=450$. 
\end{proof}

\begin{lem}\label{lem: fourboxesoneinmiddle}
    The following equality holds in $\mc{P}$.
     \begin{align*}
        \adjustbox{scale=0.95}{{\begin{overpic}[unit=1mm, scale=0.6]{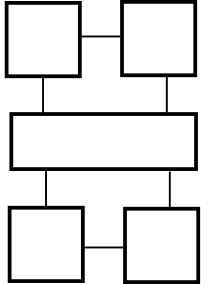}
         \put(24,5){$S$}
         \put(6, 5){$S$}
         \put(6, 38){$S$}
         \put(16, 40){\tiny{$4$}}
         \put(23, 38){$S$}
         \put(3,14){\tiny{$4$}}
         \put(28,14){\tiny{$4$}}
         \put(16, 1){\tiny{$4$}}
          \put(3, 30){\tiny{$4$}}
            \put(28, 30){\tiny{$4$}}
            \put(15, 22){$f^{(8)}$}
        \end{overpic}}}
       \hspace{2mm}{\raisebox{1.9cm}{$=30$}} 
    \end{align*}  
\end{lem}

\begin{proof}
    First, we use that $S^2=S+6f^{(4)}$ to reduce the two top and two bottom $S$-boxes. All but one of the terms in this expansion will involve an $f^{(4)}$ box on top or on bottom of $f^{(8)}$. All terms in the expansion of $f^{(4)}$ will cup or cap the $f^{(8)}$ box. Therefore, the only nonzero term is:
     \begin{align*}
        \adjustbox{scale=0.95}{{\begin{overpic}[unit=1mm, scale=0.6]{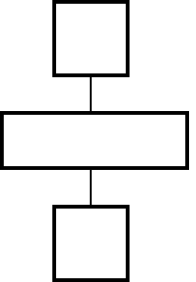}
         \put(13, 5){$S$}
         \put(13, 37){$S$}
         \put(10,14){\tiny{$8$}}
          \put(10, 29){\tiny{$8$}}
            \put(12, 21){$f^{(8)}$}
        \end{overpic}}}
    \end{align*} 

    Then we can expand $f^{(8)}$ using Frenkel and Khovanov's result from Proposition \ref{prop: MorrisonJWexpansion}. All the terms with a cup will be zero since $S$ is U.C.C.S. So only the first term in the expansion, which has coefficient 1 will be nonzero. This then results in $\te{tr}(S^2)$ which equals $30$.
\end{proof}

\subsection{Trivalent Graphical Notation}

In this section, we introduce the reader to notation from \cite{KL94} and give some computations that will later be needed.

\begin{defn}\label{defn: mnell}
    Let $a,b,$ and $c \in \mbb{Z}_{\geq 0}$. If $a+b-c$, $a+c-b$, and $b+c-a$ are all positive and even, we say $(a,b,c)$ is \tbf{admissible}. In this case, define
    \begin{align*}
        m \coloneq \f{a+b-c}{2}, \quad n \coloneq \f{b+c-a}{2},  \te{ and }\ell \coloneq \f{a+c-b}{2}.
    \end{align*}
\end{defn}

\begin{defn} (Kauffman-Lins \cite{KL94}, Chapter 4)
    An unlabelled box in the $a$-box space is $f^{(a)}$. That is,
    \begin{align*}
        \hspace{6mm}{\raisebox{0.3cm}{$f^{(a)}\coloneq$}}\hspace{2mm}
         \begin{overpic}[unit=1mm, scale=0.2]{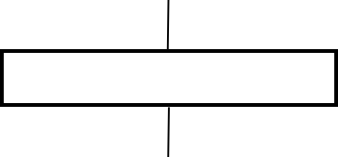}
         \put(10, 6.5){\tiny{$a$}}
        \end{overpic}.
    \end{align*}
    Suppose $(a,b,c)$ is admissible. A trivalent vertex is defined as follows
    \begin{align*}
         \begin{overpic}[unit=1mm, scale=0.5]{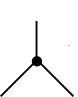}
         \put(2, 9){\tiny{$a$}}
          \put(0, 3){\tiny{$b$}}
           \put(8, 3){\tiny{$c$}}
        \end{overpic}
         {\raisebox{0.3cm}{$\coloneq$}}\hspace{2mm}
    {\raisebox{-1.2cm}{\begin{overpic}[unit=1mm, scale=0.2]{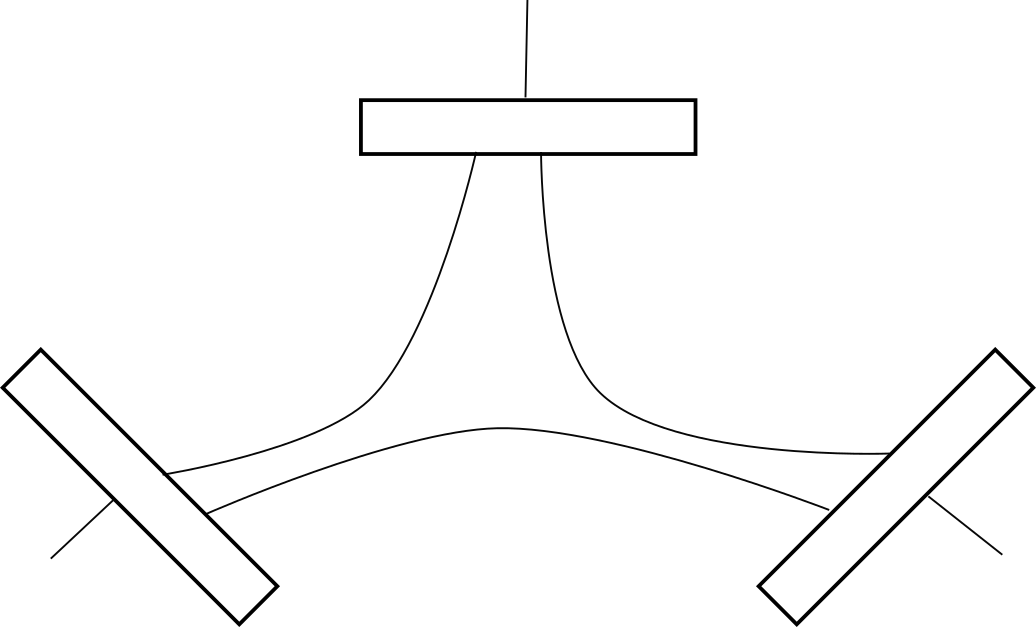}
         \put(24, 30){\tiny{$a$}}
          \put(1, 5){\tiny{$b$}}
           \put(52, 6){\tiny{$c$}}
            \put(18, 20){\tiny{$m$}}
          \put(35, 15){\tiny{$\ell$}}
           \put(24, 6){\tiny{$n$}}
        \end{overpic}}}.
    \end{align*}  
\end{defn}

In particular, this gives for $(a,b,k)$ admissible that 
 \begin{align*}
         \begin{overpic}[unit=1mm, scale=0.5]{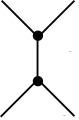}
         \put(1, 11){\tiny{$a$}}
          \put(8, 11){\tiny{$b$}}
          \put(1, 3){\tiny{$a$}}
          \put(6, 7){\tiny{$k$}}
           \put(8, 3){\tiny{$b$}}
        \end{overpic}
         \hspace{4mm}{\raisebox{0.3cm}{$\coloneq$}}\hspace{2mm}
    {\raisebox{-0.8cm}{\begin{overpic}[unit=1mm, scale=0.2]{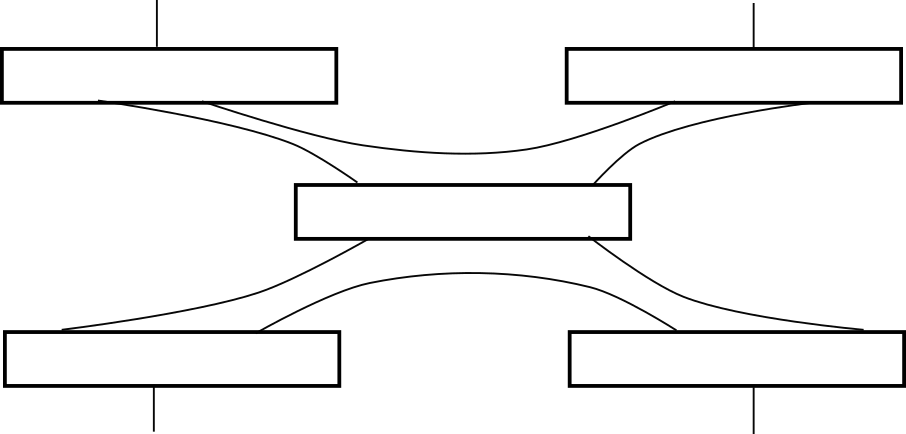}
         \put(4, 22){\tiny{$a$}}
          \put(42, 22){\tiny{$b$}}
           \put(42, 0){\tiny{$b$}}
           \put(4, 0){\tiny{$a$}}
            \put(24, 16){\tiny{$m$}}
          \put(35, 14){\tiny{$n$}}
          \put(40, 7){\tiny{$n$}}
           \put(24, 6){\tiny{$m$}}
           \put(5, 7){\tiny{$\ell$}}
            \put(10, 14){\tiny{$\ell$}}
        \end{overpic}}}.
    \end{align*}  
    where the middle box is $f^{(k)}$. 
\begin{defn} (Kauffman-Lins \cite{KL94}, Chapter 6)
    Let $(a,b,c)$ be admissible. Define: 
    \begin{align*}
         {\raisebox{-0.3cm}{$\theta(a,b,c)\coloneq$}}\hspace{6mm}
    {\raisebox{-0.5cm}{\begin{overpic}[unit=1mm, scale=0.6]{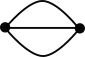}
        \put(6, 10){\tiny{$a$}}
          \put(6, 2.5){\tiny{$b$}}
           \put(6, -2){\tiny{$c$}}
        \end{overpic}}}
    \end{align*}  
    and 
    \begin{align*}
         {\raisebox{0.8cm}{$\te{Net}(m,n,\ell)\coloneq$}}\hspace{6mm}\begin{overpic}[unit=1mm, scale=0.2]{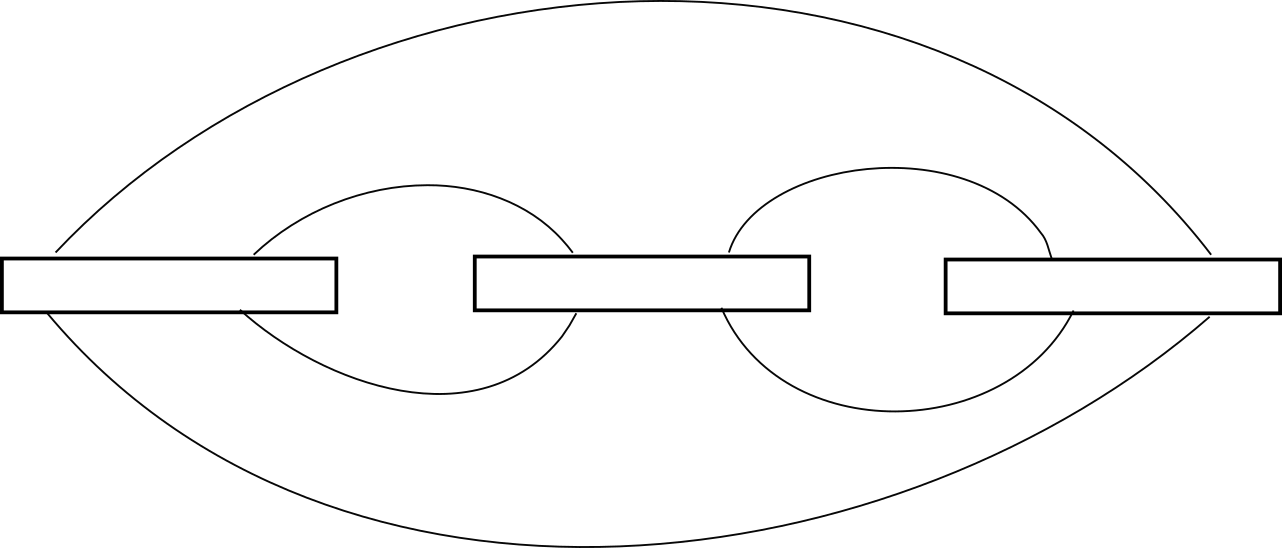}
         \put(24, 30){\tiny{$\ell$}}
          \put(42, 20){\tiny{$n$}}
           \put(15, 20){\tiny{$m$}}
        \end{overpic}.
    \end{align*}  
\end{defn}

The following equalities were found in \cite{KL94}.

\begin{prop} (Kauffman-Lins \cite{KL94}, Lemma 8 and Corollary 2)
    Let $a,b,c \in \mbb{Z}_{\geq 0}$ where $(a,b,c)$ is admissible and define $m,n,$ and $\ell$ as in Definition \ref{defn: mnell}. Then $\theta(a,b,c)=Net(m,n,\ell)$. Furthermore, \begin{align*}         \te{Net}(m,n,\ell)=\f{(-1)^{m+n+\ell}[m]![n]![\ell]![m+n+\ell+1]!}{[m+n]![n+\ell]![m+\ell]!} \end{align*} In particular, when $m+n$ is even, $\te{Net}(m,n,0)=[m+n+1]$.
\end{prop}

Define $\mathlarger{\mathlarger{\sum_{k=r}^s}} {\hspace{-0.3cm}\tiny{2}}\hspace{0.3cm}f(k)\coloneq f(r)+f(r+2)+...+f(s-2)+f(s)$. Then the following was proven by Chen in \cite{Che14}.

\begin{prop}\label{prop: chen} (Chen \cite{Che14}, Theorem 3.3.5) For $a,b \geq 0$,
    \begin{align*}
        f^{(a)}\otimes f^{(b)}= \mathlarger{\mathlarger{\sum_{k=|a-b|}}^{a+b}} {\hspace{-0.6cm}\tiny{2}} \hspace{0.5cm} \left( \f{[k+1]}{|\theta(a,b,k)|}  \hspace{1mm}{\raisebox{-0.5 cm}{\begin{overpic}[unit=1mm, scale=0.5]{twotrivalentvertices.png}
         \put(1, 11){\tiny{$a$}}
          \put(8, 11){\tiny{$b$}}
          \put(1, 3){\tiny{$a$}}
         \put(6, 7){\tiny{$k$}}
           \put(8, 3){\tiny{$b$}}
        \end{overpic}}}\right)
    \end{align*}
\end{prop}

\begin{rem*}
    Recall from Lemma \ref{Joneswenzlformulaupto3} that $[n]=n$ for index 4 planar algebras. 
\end{rem*} 

\begin{lem}\label{lem: thetavalues55}
    We can calculate the following $\theta$-values.
    \begin{enumerate}[label=(\roman*)]
        \item $|\theta(5,5,0)|=6$,
        \item $|\theta(5,5,2)|=\f{21}{5}$,
        \item $|\theta(5,5,4)|=\f{14}{5}$,
        \item $|\theta(5,5,6)|=\f{63}{25}$,
        \item $|\theta(5,5,8)|=\f{18}{5}$, and
        \item $|\theta(5,5,10)|=11$.
    \end{enumerate}
\end{lem}

\begin{proof}
    In all parts $a=b=5$. When $c=0$, this gives $m=5$, $n=0$, and $\ell=0$. Thus $|\theta(5,5,0)|=|\te{Net}(5,0,0)|=6$. We calculate the other values of $\theta$ based on their $c$-value below.
    \begin{center}
         \begin{tabular}{c|c|c}
            Value of $c$ & $(m,n,\ell)$ & $|\te{Net}(m,n,l)|$\\
            \hline
              \rule{0pt}{10pt} & &  \\[-9pt]
         $2$ &  $(4,1,1)$ & $\f{4!1!1!7!}{5!2!5!}=\f{21}{5}$ \\[4pt]
        $4$ & $(3,2,2)$ & $\f{3!2!2!8!}{5!5!4!}=\f{14}{5}$\\[4pt]
        $6$ & $(2,3,3)$ & $\f{2!3!3!9!}{5!5!6!}=\f{63}{25}$\\[4pt]
        $8$ & $(1,4,4)$ & $\f{1!4!4!10!}{5!5!8!}=\f{18}{5}$\\[4pt]
        $10$ & $(0,5,5)$ & $\f{0!5!5!11!}{5!5!10!}=11$,
    \end{tabular}
    \end{center}   
    as we wished.
\end{proof}

\begin{lem}\label{lem: thetavalues66}
    We can calculate the following $\theta$-values.
    \begin{enumerate}[label=(\roman*)]
        \item $|\theta(6,6,0)|=7$,
        \item $|\theta(6,6,2)|=\f{14}{3}$,
        \item $|\theta(6,6,4)|=\f{14}{5}$,
        \item $|\theta(6,6,6)|=\f{21}{10}$,
        \item $|\theta(6,6,8)|=\f{11}{5}$, 
        \item $|\theta(6,6,10)|=\f{11}{3}$, and
        \item $|\theta(6,6,12)|=13$.
    \end{enumerate}
\end{lem}

\begin{proof}
    In all parts $a=b=6$. When $c=0$, we get that $m=6$ and $n=\ell=0$ which gives that $|\theta(6,6,0)|=|\te{Net}(6,0,0)|=7$, We calculate the other values of $\theta$ based on their $c$-values below.
    \begin{center}
         \begin{tabular}{c|c|c}
            Value of $c$ & $(m,n,\ell)$ & $|\te{Net}(m,n,l)|$\\
            \hline 
            \rule{0pt}{10pt} & &  \\[-9pt]
         $2$ &  $(5,1,1)$ & $\f{5!1!1!8!}{6!2!6!}=\f{14}{3}$ \\[4pt]
        $4$ & $(4,2,2)$ & $\f{4!2!2!9!}{6!6!4!}=\f{14}{5}$\\[4pt]
        $6$ & $(3,3,3)$ & $\f{3!3!3!10!}{6!6!6!}=\f{21}{10}$\\[4pt]
        $8$ & $(2,4,4)$ & $\f{2!4!4!11!}{6!6!8!}=\f{11}{5}$\\[4pt]
        $10$ & $(1,5,5)$ & $\f{1!5!5!12!}{6!6!10!}=\f{11}{3}$\\[4pt]
        $12$ & $(0,6,6)$ & $\f{0!6!6!13!}{6!6!12!}=13$,
    \end{tabular}
    \end{center}  
    as we wished.
\end{proof}

\begin{lem}\label{lem: thetavalues77}
    We can calculate the following $\theta$-values.
    \begin{enumerate}[label=(\roman*)]
        \item $|\theta(7,7,0)|=8$,
        \item $|\theta(7,7,2)|=\f{36}{7}$,
        \item $|\theta(7,7,4)|=\f{20}{7}$,
        \item $|\theta(7,7,6)|=\f{66}{35}$,
        \item $|\theta(7,7,8)|=\f{396}{245}$,
        \item $|\theta(7,7,10)|=\f{286}{147}$,
        \item $|\theta(7,7,12)|=\f{26}{7}$, and
        \item $|\theta(7,7,14)|=15$.
    \end{enumerate}
\end{lem}

\begin{proof}
    In all parts $a=b=7$. When $c=0$, we get that $m=7$ and $n=\ell=0$ which gives that $|\theta(7,7,0)|=8$, We calculate the other values of $\theta$ based on their $c$-values below.
    \begin{center}
         \begin{tabular}{c|c|c}
            Value of $c$ & $(m,n,\ell)$ & $|\te{Net}(m,n,l)|$\\
            \hline
              \rule{0pt}{10pt} & &  \\[-9pt]
         $2$ &  $(6,1,1)$ & $\f{6!1!1!9!}{7!2!7!}=\f{36}{7}$ \\[4pt]
        $4$ & $(5,2,2)$ & $\f{4!2!2!10!}{7!7!4!}=\f{20}{7}$\\[4pt]
        $6$ & $(4,3,3)$ & $\f{4!3!3!11!}{7!7!6!}=\f{66}{135}$\\[4pt]
        $8$ & $(3,4,4)$ & $\f{3!4!4!12!}{7!7!8!}=\f{396}{245}$\\[4pt]
        $10$ & $(2,5,5)$ & $\f{2!5!5!13!}{7!7!10!}=\f{286}{147}$\\[4pt]
        $12$ & $(1,6,6)$ & $\f{1!6!6!14!}{7!7!12!}=\f{26}{7}$\\[4pt]
        $14$ & $(0,7,7)$ & $\f{0!7!7!15!}{7!7!14!}=15$,
    \end{tabular}
    \end{center}   
    as we wished.
\end{proof}

\subsection{The Jellyfish Algorithm}\label{e7jellyfish}

\begin{figure}[h!]
    \centering
   \begin{align*}
    \scalebox{0.39}{
        \begin{overpic}[unit=1mm, scale=0.7]{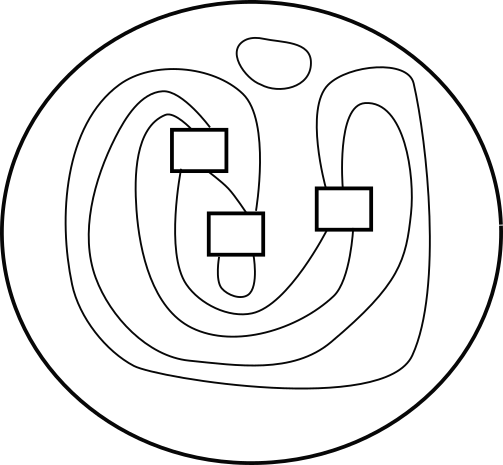}
          \put(36,57){$S$}
     \put(62,46){$S$}
     \put(42,41){$S$}
    \end{overpic} \quad\quad
    \begin{overpic}[unit=1mm, scale=0.7]{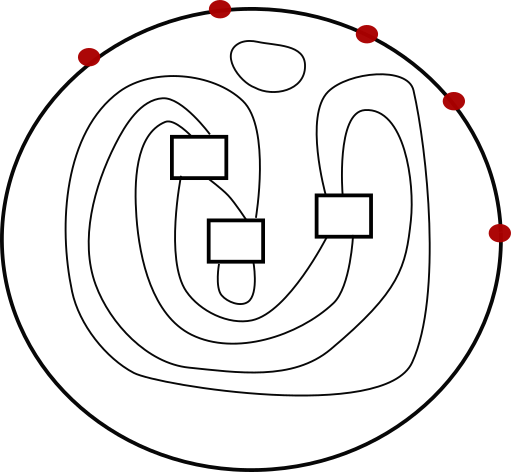}
    \put(36,57){$S$}
     \put(62,46){$S$}
     \put(42,41){$S$}
     \put(10, 82){$\mathlarger{\mathlarger{\mathlarger{m}}}$}
     \put(35,90){$\mathlarger{\mathlarger{\mathlarger{1}}}$}
     \put(65, 88){$\mathlarger{\mathlarger{\mathlarger{2}}}$}
     \put(85, 75){$\mathlarger{\mathlarger{\mathlarger{3}}}$}
     \put(96, 45){$\mathlarger{\mathlarger{\mathlarger{4}}}$}
    \end{overpic} \quad\quad
    \begin{overpic}[unit=1mm, scale=0.7]{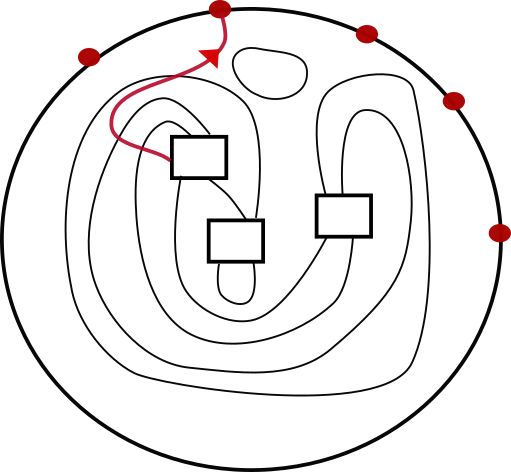}
    \put(36,57){$S$}
     \put(62,46){$S$}
     \put(42,41){$S$}
        \put(10, 82){$\mathlarger{\mathlarger{\mathlarger{m}}}$}
     \put(35,90){$\mathlarger{\mathlarger{\mathlarger{1}}}$}
     \put(65, 88){$\mathlarger{\mathlarger{\mathlarger{2}}}$}
     \put(85, 75){$\mathlarger{\mathlarger{\mathlarger{3}}}$}
     \put(96, 45){$\mathlarger{\mathlarger{\mathlarger{4}}}$}
    \end{overpic} \quad\quad
    }\\~\\
    \scalebox{0.39}{
    \begin{overpic}[unit=1mm, scale=0.7]{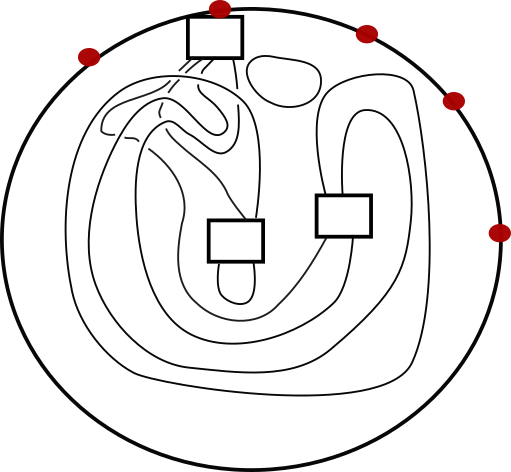}
    \put(38,79){$S$}
     \put(62,46){$S$}
     \put(42,41){$S$}
        \put(10, 82){$\mathlarger{\mathlarger{\mathlarger{m}}}$}
     \put(35,90){$\mathlarger{\mathlarger{\mathlarger{1}}}$}
     \put(65, 88){$\mathlarger{\mathlarger{\mathlarger{2}}}$}
     \put(85, 75){$\mathlarger{\mathlarger{\mathlarger{3}}}$}
     \put(96, 45){$\mathlarger{\mathlarger{\mathlarger{4}}}$}
    \end{overpic}\quad \quad
    \begin{overpic}[unit=1mm, scale=0.7]{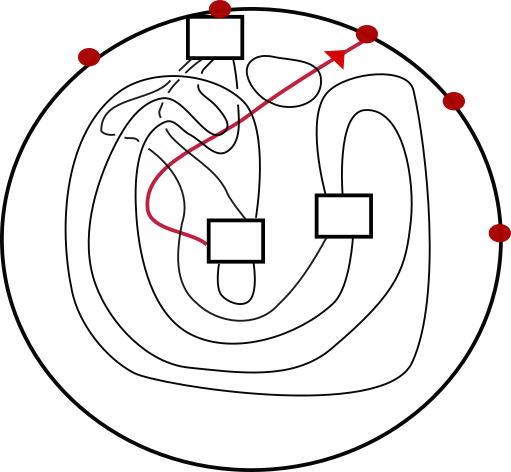}
    \put(38,79){$S$}
     \put(62,46){$S$}
     \put(42,41){$S$}
        \put(10, 82){$\mathlarger{\mathlarger{\mathlarger{m}}}$}
     \put(35,90){$\mathlarger{\mathlarger{\mathlarger{1}}}$}
     \put(65, 88){$\mathlarger{\mathlarger{\mathlarger{2}}}$}
     \put(85, 75){$\mathlarger{\mathlarger{\mathlarger{3}}}$}
     \put(96, 45){$\mathlarger{\mathlarger{\mathlarger{4}}}$}
    \end{overpic}  \quad\quad
    \begin{overpic}[unit=1mm, scale=0.7]{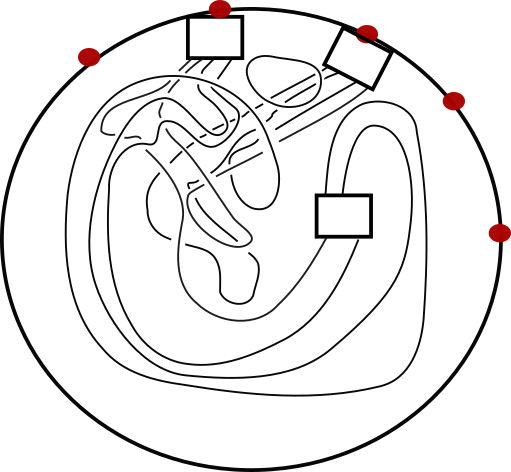} 
    \put(38,79){$S$}
     \put(62,46){$S$}
     \put(65,75){$S$}
        \put(10, 82){$\mathlarger{\mathlarger{\mathlarger{m}}}$}
     \put(35,90){$\mathlarger{\mathlarger{\mathlarger{1}}}$}
     \put(65, 88){$\mathlarger{\mathlarger{\mathlarger{2}}}$}
     \put(85, 75){$\mathlarger{\mathlarger{\mathlarger{3}}}$}
     \put(96, 45){$\mathlarger{\mathlarger{\mathlarger{4}}}$}
    \end{overpic}\quad\quad}\\~\\
    \scalebox{0.39}{\begin{overpic}[unit=1mm, scale=0.7]{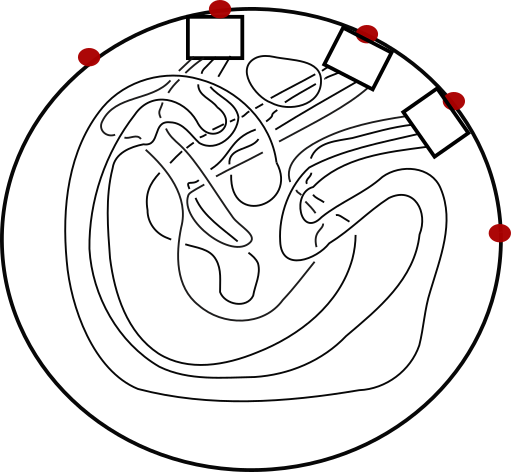}
    \put(38,79){$S$}
     \put(80,63){$S$}
     \put(65,75){$S$}
        \put(10, 82){$\mathlarger{\mathlarger{\mathlarger{m}}}$}
     \put(35,90){$\mathlarger{\mathlarger{\mathlarger{1}}}$}
     \put(65, 88){$\mathlarger{\mathlarger{\mathlarger{2}}}$}
     \put(85, 75){$\mathlarger{\mathlarger{\mathlarger{3}}}$}
     \put(96, 45){$\mathlarger{\mathlarger{\mathlarger{4}}}$}
    \end{overpic} \quad\quad
    \begin{overpic}[unit=1mm, scale=0.7]{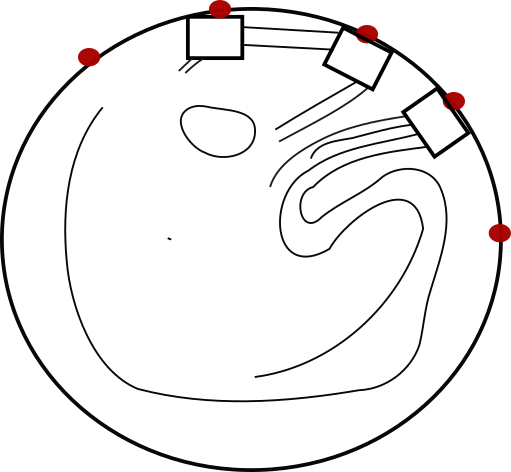}
     \put(38,79){$S$}
     \put(80,63){$S$}
     \put(65,75){$S$}
        \put(10, 82){$\mathlarger{\mathlarger{\mathlarger{m}}}$}
     \put(35,90){$\mathlarger{\mathlarger{\mathlarger{1}}}$}
     \put(65, 88){$\mathlarger{\mathlarger{\mathlarger{2}}}$}
     \put(85, 72){$\mathlarger{\mathlarger{\mathlarger{3}}}$}
     \put(96, 45){$\mathlarger{\mathlarger{\mathlarger{4}}}$}
     \put(31,66){$\begin{tikzpicture}
          \draw[cyan, thick] (-3,-1) ellipse (2.4cm and 1.4cm);
     \end{tikzpicture}$}
    \put(85, 78){$\mathlarger{\mathlarger{\mathlarger{\mathlarger{\mathlarger{\mathlarger{\mathlarger{S^2=6f^{(4)}+S}}}}}}}$}
    \end{overpic}\quad\quad\quad\quad\quad\quad\quad\quad\quad\quad\quad\quad\quad\quad\quad\quad\quad\quad}
    \end{align*}
    \caption{Schematic illustration of the jellyfish algorithm for the affine $E_7$ planar algebra (viewed in a circular tank).}
    \label{fig:jellyfishalgorithm}
\end{figure}

We now give a jellyfish algorithm for evaluating closed diagrams in $\mc{P}$. See figure \ref{fig:jellyfishalgorithm} for a schematic illustration. Later, we will show this algorithm is a well-defined function onto $\mbb{C}$. This will show that the 0-box space is one-dimensional. 

 \begin{tcolorbox}[breakable, pad at break*=0mm]
    Let $J:\mc{P}_0\to \mbb{C}$ be a function defined on diagrams in the following way. Let $D\in \mc{P}_0$. 
    \begin{enumerate}
        \item If $D$ has no $S$-boxes, skip to step 7(iv) and use $D=D_i$. 
        \item Assume $D$ has $m\geq 1$ $S$-boxes. Make $m$ evenly spaced points at the top of the diagram of $D$, starting with a point at the left side and ending at the right side. Label the points $1,...,m$. 
        \item Pick an $S$-box and label is $S_1$. From the distinguished star of $S_1$ to the point 1 at the top of $D$, draw an imaginary arc only traversing strands (the arc cannot pass through a crossing or box) and not passing through the same point more than twice. 
        \item Drag $S_1$ by stretching its strands along the imaginary arc, going under strands it passes along the way.  
        \item Repeat steps 3 and 4 $m$ times for the rest of the $S$ boxes and points. 
        \item Resolve all the crossings, which may make $D$ as sum of diagrams, say $D=\sum_{i=1}^{n_\ell} a_i D_i$, where $a_i \in \mbb{C}$. 
        \item Look individually at each $D_i$.
        \begin{enumerate}[label=(\roman*)]
            \item If $D_i$ has any $S$-box that is cupped, define $J(D_i)=0$. 
            \item If $D_i$ has no $S$-box that is cupped, then there must be a pair of adjacent $S$-boxes connected by 4 strands. Use the $S^2$ relation to rewrite $D_i$ as a linear combination of diagrams with at least one less $S$-box. 
            \item Repeat steps 7(i) and 7(ii) on all summands resulting from $D_i$ until no $S$-boxes are remaining. 
            \item Now each summand resulting from $D$ is in Temperley-Lieb. Define $J(D_i)$ to be its evaluation in $\mc{TL}$. 
        \end{enumerate}
    \end{enumerate}
    Extend the function linearly to be defined on all of $\mc{P}_0$. 
 \end{tcolorbox}
 
 We focus this section on proving the below theorem.

\begin{thm}\label{thm: e7jfaiswelldefined}
    The map $J:\mc{P}_0 \to \mbb{C}$ is well-defined, i.e., 
    \begin{enumerate}[label=(\roman*)]
        \item the evaluation algorithm does not depend on the order of $S^2$ to reduce in step 7(ii),
        \item the algorithm does not depend on the choice of arc in step 3,
        \item the algorithm does not depend on the order of the $S$-boxes chosen to float to the top in step 3, and  
        \item the algorithm respects the defining relations.
    \end{enumerate}
\end{thm}

\begin{rem*}
    At step 7 of the jellyfish algorithm, for each summand, we have all the $S$-boxes at the top of the tank and strands with no intersections on the bottom. Without loss of generality, we can instead assume our tank is a circle. Then, for example, we can view the below diagram as being the two top $S$-boxes adjacent in the circle and the rest of diagram being at the bottom of the circle, where $T$ is in Temperley-Lieb. If instead the jellyfish algorithm has only been completed up to step 5, we can assume that there are crossings in $T$.
    \begin{equation*}
        \adjustbox{scale=0.95}{\begin{overpic}[unit=1mm, scale=0.6]{fourboxesoneinmiddle.png}
         \put(25,5){$S$}
         \put(6, 5){$S$}
         \put(6, 38){$S$}
         \put(25, 38){$S$}
         \put(3,14){\tiny{$c$}}
         \put(28,14){\tiny{$d$}}
          \put(3, 30){\tiny{$a$}}
            \put(28, 30){\tiny{$b$}}
            \put(14, 21){$T$}
        \end{overpic}}
    \end{equation*}  
\end{rem*}

\subsection{The Jellyfish Algorithm: Order of $S^2$.}\label{jellyfish1}
\begin{lem}\label{lem: e7sufficientwelldefinedSsquared}
    The jellyfish algorithm is well defined with respect to the order in reducing $S^2$. 
\end{lem}

\begin{proof}
    Let $D$ be a diagram in $\mc{P}_0$ and suppose the jellyfish algorithm has been completed for $D$ up through step 6. Consider $D_i$ from the summand created in step 6. We will prove the jellyfish algorithm is well defined with respect to the order in reducing $S^2$ by induction on $n$, the number of $S$-boxes in $D_i$ at this point. 

    If $D_i$ has 0, 1, or 2 $S$-boxes, then there is no choice in the order of reduction using $S^2$, so the base case is true. Suppose that for some $n\geq 2$, if $D_i$ has $n$ $S$-boxes, the choice in the order of reducing $S^2$ does not matter. We want to show that for $1\leq j < k \leq \ell <p\leq n+1$, if $\{S_j, S_k\}$ and $\{S_\ell, S_p\}$ both form an $S^2$, reducing $\{S_j, S_k\}$ first is the same number as reducing $\{S_\ell, S_p\}$ first. (Notice that this covers the $n=3$ case since we allow $k=\ell$).

    We prove in cases. Case 1 is assuming that $k=\ell$. That is, we want to show that first reducing $\{S_j, S_k\}$ evaluates to the same as first reducing $\{S_k, S_p\}$. These boxes must all be adjacent, so at Step 6, we have the following picture:
    \begin{align*}
       \adjustbox{scale=0.95}{ \begin{overpic}[unit=1mm, scale=0.6]{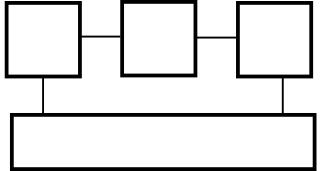}
         \put(24,20){$S_k$}
         \put(5, 20){$S_j$}
         \put(41, 20){$S_p$}
         \put(41, 12){\tiny{$4$}}
         \put(25, 3){$T$}
          \put(3, 12){\tiny{$4$}}
            \put(34, 22){\tiny{$4$}}
            \put(15, 22){\tiny{$4$}}
        \end{overpic}}
    \end{align*}  
    where $T$ is the rest of the diagram. When we first reduce $\{S_j, S_k\}$, we get
    \begin{align*}
        \adjustbox{scale=0.95}{  \begin{overpic}[unit=1mm, scale=0.6]{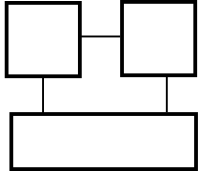}
         \put(28,12){\tiny{$4$}}
         \put(2, 11){\tiny{$4$}}
         \put(5, 20){$S$}
         \put(15, 19){\tiny{$4$}}
         \put(24, 20){$S$}
         \put(15, 3){$T$}
        \end{overpic}}
        \hspace{2mm}{\raisebox{1cm}{$+6$}}\hspace{2mm}
           \adjustbox{scale=0.95}{   \begin{overpic}[unit=1mm, scale=0.6]{twoboxesconnectedoneboxbelowclosed.png}
         \put(28,12){\tiny{$4$}}
         \put(2, 11){\tiny{$4$}}
         \put(4, 20){$f^{(4)}$}
         \put(15, 19){\tiny{$4$}}
         \put(24, 20){$S$}
         \put(15, 3){$T$}
        \end{overpic}}.
    \end{align*}  
    We now have two diagrams with at least one less $S$-box, so the choice in the order of reduction from here on out does not matter. Continuing step 7, we can notice that in the second summand, since $S$ is U.C.C.S., $S$ will absorb $f^{(4)}$. Then we can reduce the above diagram to 
       \begin{align*}
        \adjustbox{scale=0.95}{ \begin{overpic}[unit=1mm, scale=0.6]{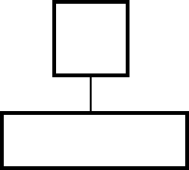}
         \put(13, 20){$S^2$}
         \put(17, 11){\tiny{$8$}}
         \put(13, 4){$T$}
        \end{overpic}}
        \hspace{2mm}{\raisebox{1cm}{$+6$}}\hspace{2mm}
            \adjustbox{scale=0.95}{  \begin{overpic}[unit=1mm, scale=0.6]{twoboxesstacked.png}
          \put(13, 20){$S$}
         \put(17, 11){\tiny{$8$}}
         \put(13, 4){$T$}
        \end{overpic}}
    \end{align*}  
    This is identical to what is obtained if instead $\{S_k, S_p\}$ is first reduced. The rest of the algorithm will then be identical for both cases.

    In case 2, assume $k\neq \ell$. Then at step 6 we have the following picture:
    \begin{align*}
        \adjustbox{scale=0.95}{ \begin{overpic}[unit=1mm, scale=0.6]{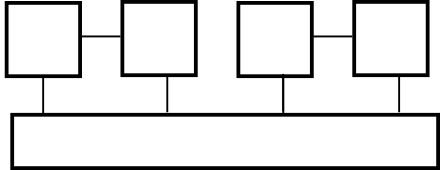}
         \put(5, 20){$S_j$}
         \put(23,20){$S_k$}
         \put(42, 20){$S_\ell$}
         \put(60,20){$S_p$}
            \put(34, 3){$T$}
           \put(3, 11){\tiny{$4$}}
           \put(15, 22){\tiny{$4$}}
         \put(22, 11){\tiny{$4$}}
         \put(52, 22){\tiny{$4$}}
         \put(42, 11){\tiny{$4$}}  
         \put(60,11){\tiny{$4$}} 
         \put(32.5, 22){\tiny{$...$}}
        \end{overpic}}
    \end{align*}  
    Reducing $\{S_j, S_k\}$ first gives
    \begin{align*}
     {\raisebox{1cm}{$6$}}\hspace{2mm}
         \adjustbox{scale=0.95}{   \begin{overpic}[unit=1mm, scale=0.6]{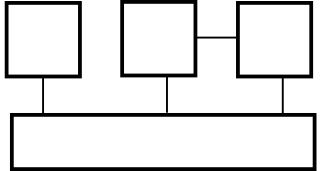}
         \put(4, 20){$f^{(4)}$}
         \put(23,20){$S_k$}
         \put(42, 20){$S_\ell$}
            \put(26, 3){$T$}
           \put(3, 11){\tiny{$8$}}
           \put(15, 22){\tiny{$...$}}
         \put(22, 11){\tiny{$4$}}
         \put(42, 11){\tiny{$4$}}  
         \put(32.5, 22){\tiny{$4$}}
        \end{overpic}}
          \hspace{2mm}  {\raisebox{1cm}{$+$}}\hspace{2mm}
                 \adjustbox{scale=0.95}{    \begin{overpic}[unit=1mm, scale=0.6]{3boxesin1pairand1below.png}
         \put(5, 20){$S$}
         \put(23,20){$S_k$}
         \put(42, 20){$S_\ell$}
            \put(26, 3){$T$}
           \put(3, 11){\tiny{$8$}}
           \put(15, 22){\tiny{$...$}}
         \put(22, 11){\tiny{$4$}}
         \put(42, 11){\tiny{$4$}}  
         \put(32.5, 22){\tiny{$4$}}
        \end{overpic}}
    \end{align*}
    Each of these components have one less $S$-box so the next choice to reduce in the pairing does not matter. So we choose next to pair $\{S_k, S_\ell\}$ and obtain
    \begin{align*}
    \scalebox{0.9}{{\raisebox{1cm}{$36$}}\hspace{2mm}
        \adjustbox{scale=0.95}{  \begin{overpic}[unit=1mm, scale=0.6]{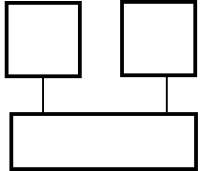}
         \put(28,12){\tiny{$4$}}
         \put(2, 11){\tiny{$4$}}
         \put(4, 20){$f^{(4)}$}
         \put(15, 19){\tiny{$...$}}
         \put(23, 20){$f^{(4)}$}
         \put(15, 3){$T$}
        \end{overpic}
        \hspace{2mm}{\raisebox{1cm}{$+6$}}\hspace{2mm}
             \begin{overpic}[unit=1mm, scale=0.6]{2unconnectedthen1.png}
         \put(28,12){\tiny{$4$}}
         \put(2, 11){\tiny{$4$}}
         \put(4, 20){$f^{(4)}$}
         \put(15, 19){\tiny{$...$}}
         \put(24, 20){$S$}
         \put(15, 3){$T$}
        \end{overpic}}
        \hspace{2mm}{\raisebox{1cm}{$+6$}}\hspace{2mm}
            \adjustbox{scale=0.95}{  \begin{overpic}[unit=1mm, scale=0.6]{2unconnectedthen1.png}
         \put(28,12){\tiny{$4$}}
         \put(2, 11){\tiny{$4$}}
         \put(5, 20){$S$}
         \put(15, 19){\tiny{$...$}}
         \put(23, 20){$f^{(4)}$}
         \put(15, 3){$T$}
        \end{overpic}}
        \hspace{2mm}{\raisebox{1cm}{$+$}}\hspace{2mm}
             \adjustbox{scale=0.95}{ \begin{overpic}[unit=1mm, scale=0.6]{2unconnectedthen1.png}
         \put(28,12){\tiny{$4$}}
         \put(2, 11){\tiny{$4$}}
         \put(5, 20){$S$}
         \put(15, 19){\tiny{$...$}}
         \put(24, 20){$S$}
         \put(15, 3){$T$}
        \end{overpic}}}
    \end{align*}  
    This will be the same result as first choosing to pair $\{S_\ell, S_p\}$, so the choice in the order of reducing $S^2$ where $j\neq k$ doesn't matter either.
\end{proof}

\begin{rem}\label{rem: aboutjellyfishandevaluation}
    It is now clear that if a diagram in $\mc{P}_0$ has been plugged into the jellyfish algorithm and the first six steps have been performed, then if one can evaluate the resulting diagram by only using the $S^2$ relation, the bubble relation, and that $S$ is U.C.C.S., then the evaluation of the diagram in $\mc{P}$ will equal its evaluation for the jellyfish algorithm. In particular, this gives that Lemmas \ref{lem: e7sufficienttraceofSsquared}, \ref{lem: fourboxesconnected}, \ref{lem: fourboxesconnectedequaling300}, \ref{lem: fourboxesconnectedequaling450},  and \ref{lem: fourboxesoneinmiddle} are also giving the evaluation of those diagrams in the jellyfish algorithm. Lemma \ref{lem: twoboxesconnectedoneboxbelow} also only uses the $S^2$ relation and that $S$ is U.C.C.S., so it is clear that if at step 7 of the jellyfish algorithm we have this diagram in some local neighborhood, then it will evaluate to zero.
\end{rem}

\subsection{The Jellyfish Algorithm: Choice of Arc}\label{jellyfish2}
Next, we show that the jellyfish algorithm is invariant under the choice of arc. 

First, notice that an imaginary arc is an oriented curve in $\mathbb{R}^2$ that may have self-intersections. Recall that an imaginary arc can pass through a point at most two times. When drawing an imaginary arc, if the arc intersects itself, we draw the intersection as a crossing. The lower strand of the crossing is where the curve intersected the point first, the upper strand is where the curve intersected second. 

With the above viewpoint of imaginary arcs, we can now view these arcs as neighborhoods of some oriented knot diagram, $K$, where $K$ is formed by gluing the boundary points of the arc together. Following the orientation of the imaginary arc, every time a crossing is reached, the arc first goes under the crossing. The knot $K$ is then called an \textbf{ascending knot}. This is always the unknot, as we can just pull two ends apart to unravel. Analogously, \textbf{descending knots} are also defined and are always the unknot. See \cite{Oza10} for further discussion on ascending and descending knots. 

Imaginary arcs on a diagram will be drawn in red. When the arc passes a strand of the planar algebra, we also use a crossing. The imaginary arc will always be drawn over the strands of the planar algebra for visual clarity. However, remember that when the $S$-box and its strands are dragged, the strands always pass under the strands it intersects.

 \begin{figure}[h!]
    \centering
    \begin{subfigure}{.45\textwidth}
    \centering
    \hspace{3mm} \adjustbox{scale=0.8}{ \begin{overpic}[unit=1mm, scale=0.7]{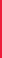}
    \end{overpic}}\quad {\raisebox{0.4cm}{{$\leftrightarrow$}}}\quad
    \adjustbox{scale=0.8}{ \begin{overpic}[unit=1mm, scale=0.7]{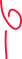}
    \end{overpic}}
    \caption{{$\Omega_{1a}$}}
    \label{fig:Omega1}
\end{subfigure}%
\begin{subfigure}{.45\textwidth}
    \centering
    \hspace{3mm} \adjustbox{scale=0.8}{ \begin{overpic}[unit=1mm, scale=0.7]{redstrand.png}
    \end{overpic}}\quad {\raisebox{0.4cm}{{$\leftrightarrow$}}}\quad
    \adjustbox{scale=0.8}{ \begin{overpic}[unit=1mm, scale=0.7]{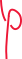}
    \end{overpic}}
    \caption{{$\Omega_{1b}$}}
    \label{fig:Omega1}
\end{subfigure}%
\\~\\
\begin{subfigure}{.45\textwidth}
  \centering
\hspace{3mm} \adjustbox{scale=0.8}{ \begin{overpic}[unit=1mm, scale=0.7]{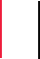}
    \end{overpic}}\quad {\raisebox{0.4cm}{{$\leftrightarrow$}}}\quad
   \adjustbox{scale=0.8}{  \begin{overpic}[unit=1mm, scale=0.7]{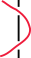}
    \end{overpic}}
  \caption{$\Omega_{2a}$}
\end{subfigure}%
\begin{subfigure}{.45\textwidth}
\centering
 \hspace{3mm} \adjustbox{scale=0.8}{ \begin{overpic}[unit=1mm, scale=0.7]{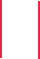}
    \end{overpic}}\quad {\raisebox{0.4cm}{{$\leftrightarrow$}}}\quad
    \adjustbox{scale=0.8}{ \begin{overpic}[unit=1mm, scale=0.7]{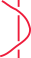}
    \end{overpic}}
  \caption{$\Omega_{2b}$}
\end{subfigure}%
\\~\\
\begin{subfigure}{.45\textwidth}
  \centering
\adjustbox{scale=0.8}{ \begin{overpic}[unit=1mm, scale=0.7]{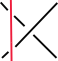}
    \end{overpic}} \quad {\raisebox{0.4cm}{{$\leftrightarrow$}}}\quad
    \adjustbox{scale=0.8}{ \begin{overpic}[unit=1mm, scale=0.7]{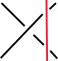}
    \end{overpic}}
  \caption{$\Omega_{3a}$}
\end{subfigure}%
\begin{subfigure}{.45\textwidth}
  \centering
 \adjustbox{scale=0.8}{  \begin{overpic}[unit=1mm, scale=0.7]{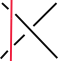}
    \end{overpic} }\quad {\raisebox{0.4cm}{{$\leftrightarrow$}}}\quad
    \adjustbox{scale=0.8}{ \begin{overpic}[unit=1mm, scale=0.7]{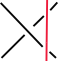}
    \end{overpic}}
  \caption{$\Omega_{3b}$}
\end{subfigure}%
\\~\\
\begin{subfigure}{.45\textwidth}
  \centering
  \adjustbox{scale=0.8}{ \begin{overpic}[unit=1mm, scale=0.7]{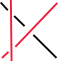}
    \end{overpic}} \quad {\raisebox{0.4cm}{{$\leftrightarrow$}}}\quad
    \adjustbox{scale=0.8}{ \begin{overpic}[unit=1mm, scale=0.7]{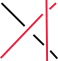}
    \end{overpic}}
  \caption{$\Omega_{3c}$}
\end{subfigure}%
\begin{subfigure}{.45\textwidth}
  \centering
  \adjustbox{scale=0.8}{ \begin{overpic}[unit=1mm, scale=0.7]{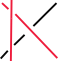}
    \end{overpic}} \quad {\raisebox{0.4cm}{{$\leftrightarrow$}}}\quad
    \adjustbox{scale=0.8}{ \begin{overpic}[unit=1mm, scale=0.7]{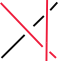}
    \end{overpic}}
  \caption{$\Omega_{3d}$}
\end{subfigure}%
\\~\\
\begin{subfigure}{.45\textwidth}
  \centering
  \adjustbox{scale=0.8}{ \begin{overpic}[unit=1mm, scale=0.7]{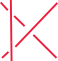}
    \end{overpic} }\quad {\raisebox{0.4cm}{{$\leftrightarrow$}}}\quad
    \adjustbox{scale=0.8}{ \begin{overpic}[unit=1mm, scale=0.7]{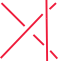}
    \end{overpic}}
  \caption{$\Omega_{3e}$}
\end{subfigure}%
\begin{subfigure}{.45\textwidth}
  \centering
  \adjustbox{scale=0.8}{ \begin{overpic}[unit=1mm, scale=0.7]{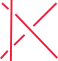}
    \end{overpic} }\quad {\raisebox{0.4cm}{{$\leftrightarrow$}}}\quad
    \adjustbox{scale=0.8}{ \begin{overpic}[unit=1mm, scale=0.7]{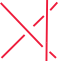}
    \end{overpic}}
  \caption{$\Omega_{3f}$}
\end{subfigure}%
\\~\\
\begin{subfigure}{.45\textwidth}
  \centering
   \adjustbox{scale=0.8}{ \begin{overpic}[unit=1mm, scale=0.4]{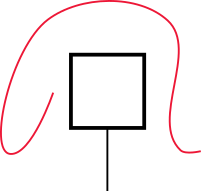}
  \put(10, 9,5){$S$}
  \put(4.2,9){$\fixedpointstar$}
    \end{overpic} }\quad {\raisebox{0.4cm}{{$\leftrightarrow$}}}\quad
   \adjustbox{scale=0.8}{  \begin{overpic}[unit=1mm, scale=0.4]{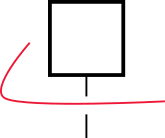}
    \put(8, 9.5){$S$}
    \put(1.7,9){$\fixedpointstar$}
    \end{overpic}}
  \caption{$\Omega_{4a}$}
\end{subfigure}%
\begin{subfigure}{.45\textwidth}
  \centering
  \adjustbox{scale=0.8}{ \begin{overpic}[unit=1mm, scale=0.4]{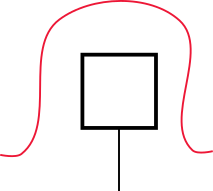}
  \put(11.5, 9){$S$}
    \end{overpic}} \quad {\raisebox{0.4cm}{{$\leftrightarrow$}}}\quad
    \adjustbox{scale=0.9}{ \begin{overpic}[unit=1mm, scale=0.4]{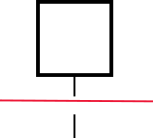}
    \put(6.5, 9.5){$S$}
    \end{overpic}}
  \caption{$\Omega_{4b}$}
\end{subfigure}%
\caption{The the $\Omega$-moves. A red strand indicates a portion of the imaginary arc, and a black strand is a portion of the strand of a diagram.}
\label{fig: omega-moves}
\end{figure}

\begin{defn}
    In figure \ref{fig: omega-moves} we define local relations involving imaginary arcs (drawn in red) and strands of the planar algebra (drawn in black). We define the collection of all the moves to be \textbf{$\Omega$-moves}. The set of $\Omega_{1a}$ and $\Omega_{1b}$-moves are called the \textbf{$\Omega_1$-moves}. Similarly, we define the \textbf{$\Omega_2$-moves}, \textbf{$\Omega_3$-moves} and \textbf{$\Omega_4$-moves}.
\end{defn}

 Notice that the $\Omega$-moves include all the classical Reidemeister-like moves where the overstrand must be an imaginary arc. The $\Omega_4$-moves allow an imaginary arc to be moved past an $S$-box both in the case where the arc belongs to the $S$-box drawn in the neighborhood, and the other where the arc may belong to a different $S$-box. (In fact, in Lemma \ref{lem: jkarespectsarc}, we prove that if two arcs with identical endpoints differ by an $\Omega_{4b}$-move in some neighborhood, $N$, then the imaginary arc belongs to a different $S$-box than the one in $N$.) Thus, we can make the following remark. 

\begin{rem*} Since the closure of every imaginary arc is an ascending knot, every imaginary arc can be unknotted by a finite sequence of $\Omega$-moves. In particular, crossing changes are not needed to unknot the arc with two fixed endpoints. The proof of Reidemeister's theorem \cite{Rei27} can thus be adapated to show that all imaginary arcs with the same fixed endpoints are connected by a finite sequence of $\Omega$-moves. 
\end{rem*}

We now prove that the jellyfish algorithm is invariant with respect to these $\Omega$-moves.

\begin{lem}\label{lem: jfaomega1}
    The jellyfish algorithm is invariant under $\Omega_{1}$-moves.
\end{lem}

\begin{proof}
    We only prove the lemma for $\Omega_{1a}$. The proof of $\Omega_{1b}$ follows a symmetric argument. Fix an $S$-box and a point $k$ on the top of the diagram created in step 2 of the jellyfish algorithm. Let $\gamma_1$ and $\gamma_2$ be two imaginary arcs identical except in a neighborhood  where $\gamma_1$ is the left-hand side of the $\Omega_{1a}$-move and $\gamma_2$ is the right-hand side of $\Omega_{1a}$. Then when dragging the $S$-box along $\gamma_1$ we have
    \begin{align*}
       \begin{overpic}[unit=1mm, scale=0.3]{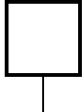}
     \put(2,4.5){$S$}
    \end{overpic} \hspace{2mm} {\raisebox{0.4cm}{{$\mapsto$}}} \hspace{2mm}
     \begin{overpic}[unit=1mm, scale=0.3]{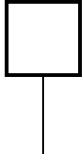}
     \put(2,8){$S$}
    \end{overpic}
    \end{align*}
    while dragging the $S$-box along $\gamma_2$ gives
    \begin{align*}
       \begin{overpic}[unit=1mm, scale=0.3]{Sbox1strand.png}
     \put(2,4.5){$S$}
    \end{overpic} \hspace{2mm} {\raisebox{0.4cm}{{$\mapsto$}}} \hspace{2mm}
    \raisebox{-0.3cm}{
     \begin{overpic}[unit=1mm, scale=0.3]{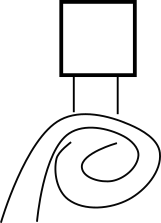}
     \put(6.5,13.5){$S$}
     \put(1.2, 2){\tiny{$...$}}
    \end{overpic}}
    \end{align*}
    Thus, to show the jellyfish algorithm is invariant under $\Omega_{1a}$ boils down to showing the below to diagrams evaluate the same in the jellyfish algorithm for any $T\in \te{Hom}(\emptyset, X^{\otimes 8})$.
    \begin{equation}\label{eqn: boxomega1box}
        \begin{overpic}[unit=1mm, scale=0.6]{twoboxesstacked.png}
         \put(17, 11){\tiny{$8$}}
         \put(13, 19){$S$}
         \put(13, 3){$T$}
        \end{overpic}\hspace{1mm}, \hspace{2mm} and \hspace{3mm}
         \begin{overpic}[unit=1mm, scale=0.6]{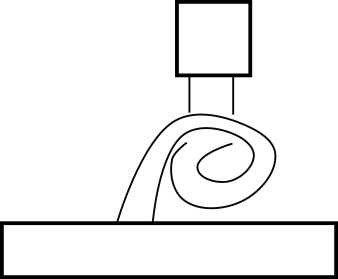}
     \put(32.5,37){$S$}
     \put(21, 12){\tiny{$...$}}
     \put(25, 4){$T$}
    \end{overpic}
    \end{equation}
    Assume that up through step 5 of the jellyfish algorithm has been completed. Suppose that at the current stage, the diagram looks like the right-hand side of (\ref{eqn: boxomega1box}). The next step is to resolve crossings. Use Lemma \ref{lem: cuponcrossingresolved} to resolve the innermost loop. This multiplies the diagram by $i$. Then use Lemma \ref{lem: boxwithcrossingofstrandtotheleft} to uncross the rightmost strand of $S$ crossed over the other 7 strands of $S$. This multiplies the diagram by $i^7$. Now the diagram is
    \begin{align*}
        \begin{overpic}[unit=1mm, scale=0.6]{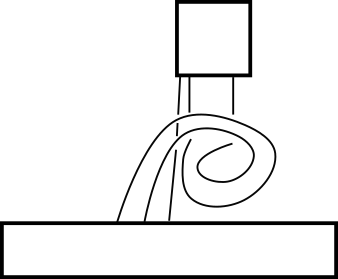}
     \put(32.5,37){$S$}
     \put(21, 12){\tiny{$...$}}
     \put(25, 4){$T$}
    \end{overpic}
    \end{align*}
    Repeating this for the seven other twists, we obtain the diagram on the left of (\ref{eqn: boxomega1box}), which gives the result.
\end{proof}

A nearly identical proof gives the following lemma.

\begin{lem}
    The jellyfish algorithm is invariant under $\Omega_{4a}$-moves.
\end{lem}

\begin{lem}
    The jellyfish algorithm is invariant under $\Omega_{2}$-moves.
\end{lem}

\begin{proof}
    The proof for $\Omega_{2a}$ and $\Omega_{2b}$ are the same. We just focus on $\Omega_{2a}$. Fix an $S$-box and a point on the boundary. Let $\gamma_1$ and $\gamma_2$ be two imaginary arcs from the $S$-box's star to the point on the boundary such that $\gamma_1$ and $\gamma_2$ are identical except in a neighborhood where $\gamma_1$ locally looks like the left-hand side of $\Omega_{2a}$ and $\gamma_2$ looks like the right-hand side. Complete steps 4 and 5 of the jellyfish algorithm identically in both cases. These two cases are identical except in a local neighborhood. In step 6, apply Reidemeister II to make the neighborhoods identical. 
\end{proof}

The same idea for proof as above except using Reidemeister III proves the following.
\begin{lem}
    The jellyfish algorithm is invariant under $\Omega_{3}$-moves.
\end{lem}

We can now prove the main lemma of this subsection.

\begin{lem}\label{lem: jkarespectsarc}
    The jellyfish algorithm is invariant under choice of arc.
\end{lem}

\begin{proof}
    Let $D$ and $D'$ be identical diagrams in $\mathcal{P}_0$, each with $m$ $S$-boxes. We prove the lemma by induction on the number of $S$-boxes left to be dragged to the top. Notice that, by the previous lemmas of this subsection, all that is left to show is that the jellyfish algorithm is invariant under $\Omega_{4b}$-moves. 

    Suppose that we have taken equivalent steps of the jellyfish algorithm for $D$ and $D'$ until there is only one $S$-box remaining. We argue that for this final $S$-box, $S_m$, the imaginary arc cannot differ by an $\Omega_{4b}$-move. 
    
    Suppose, for contradiction, that $\gamma_1$ and $\gamma_2$ are the imaginary arcs for the final $S$-box for $D$ and $D'$, respectively. Suppose $\gamma_1$ and $\gamma_2$ are identical except in a neighborhood $N$ where $\gamma_1$ looks like the left-hand side of the $\Omega_{4b}$-move and $\gamma_2$ looks like the right-hand side, where the $S$-box is the same in both diagrams. Since $S_m$ is the final $S$-box not at the top of the diagram, the $S$-box in $N$ must be $S_m$. However, then there would be a portion of the imaginary arc stemming from the $S$-box's star in $N$. Further, this arc, without crossing the portion of imaginary arc already drawn in $N$, would have to take different paths for $\gamma_1$ and $\gamma_2$. Therefore, when there is only one $S$-box left to be dragged to the top, the choice of arc doesn't matter. 

    Next, suppose that when we have take equivalent steps of the jellyfish algorithm for $D$ and $D'$ until there are $k$ $S$-boxes left to be dragged to the top, then different choices of arc for the rest of the $S$-boxes results in equivalent evaluations of the jellyfish algorithm for $D$ and $D'$.  

    Now suppose that we have taken equivalent steps of the jellyfish algorithm on $D$ and $D'$ until there are $k+1$ $S$-boxes remaining. Choose the same $S$-box in $D$ and $D'$ that will be dragged to the point $m-k$ at the top of the diagram. Choose imaginary arcs $\gamma_1$ and $\gamma_2$ to drag $S_{m-k}$ in $D$ and $D'$, respectively. Suppose these two arcs are identical except in a neighborhood $N$ where $\gamma_1$ looks like the left-hand side of $\Omega_{4b}$ and $\gamma_2$ looks like the right-hand side.

    By the same argument as in the base case, the $S$-box in $N$ is not $S_{m-k}$. Further, it is not any of the $m-(k-1)$ $S$-boxes already at the top of the diagram.  Call the $S$-box in $N$, $S_i$. Finish dragging $S_{m-k}$ to the top of $D$ and $D'$. By induction, the rest of the $S$-boxes being dragged to the top can use any choice of arc. Except for $S_i$, choose arcs dragging the rest of the $S$-boxes that avoid the neighborhood $N$. When dragging $S_i$ to the top, for $D$, choose an imaginary arc that begins by immediately going upwards, crossing the region of strands that were dragged following the imaginary arc in $N$ by $\gamma_1$. After leaving $N$ through this upwards path, never enter $N$ again. For $D'$, choose the same arc as chosen in $D$. After dragging $S_i$ to the top in both $D$ and $D'$ the diagrams are the same except in $N$, where the neighborhood of $D$ looks like the left diagram below and the neighborhood of $D'$ looks like the right diagram below:
        \begin{align*}
 \scalebox{0.5}{\hspace{2mm}
    \begin{overpic}[unit=1mm, scale=0.6]{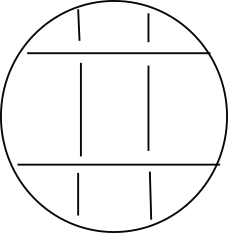}
          \put(3,20){$.$}
         \put(3,19){$.$}
         \put(3,18){$.$}
        \put(18, 18){{$...$}}
        \end{overpic}},
               \scalebox{0.5}{ \begin{overpic}[unit=1mm, scale=0.6]{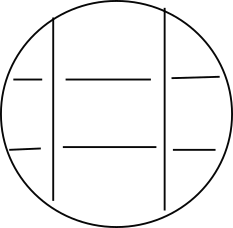}
         \put(3,20){$.$}
         \put(3,19){$.$}
         \put(3,18){$.$}
        \put(15, 18){$...$}
        \end{overpic}
         \hspace{2mm}},
    \end{align*}
    where there are both eight vertical and horizontal strands in each figure. Therefore, by Lemma \ref{lem: switchevennumberofcrossings}, these two diagrams will be equivalent after completing step 6. Thus, we can conclude that the jellyfish algorithm is invariant under the choice of arc. 
    
\end{proof}

\subsection{The Jellyfish Algorithm: Order of Floating $S$-boxes}\label{jellyfish3}

\begin{lem}\label{lem: twistanduntwistthesame1}
    Fix some $c\in\{0,...,8\}$ and $T\in \te{Hom}(\emptyset, X^{\otimes 2(8-c)})$. Assume the jellyfish algorithm has been completed up to step 5 on two diagram, which only differ as shown below in (\ref{eq: twoequivalentdiagramsuntwist}). Similarly, assume the jellyfish algorithm has been completed up to step 5 on two diagrams which only differ as shown below in (\ref{eq: twoequivalentdiagramsuntwistwithf}). Then following steps 6 and 7 on each of the diagrams in (\ref{eq: twoequivalentdiagramsuntwist}) results in the same number if and only if following steps 6 and 7 on each of the diagrams in (\ref{eq: twoequivalentdiagramsuntwistwithf}) results in the same number.
    
    \begin{equation}\label{eq: twoequivalentdiagramsuntwist}
         \begin{overpic}[unit=1mm, scale=0.6]{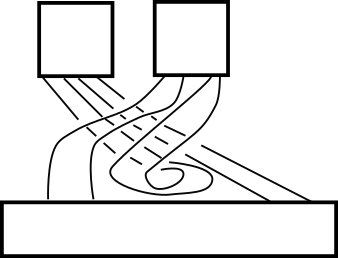}
        \put(10, 34){$S$}
        \put(29, 34){$S$}
        \put(10, 13){\tiny{$...$}}
        \put(20,3){$T$}
        \put(20, 12){\tiny{$...$}}
        \put(19, 13){\tiny{$(c)$}}
        \put(39,13){\tiny{$...$}}
    \end{overpic}\hspace{1mm},\hspace{4mm}
       \begin{overpic}[unit=1mm, scale=0.6]{twoboxesconnectedoneboxbelowclosed.png}
         \put(28,12){\tiny{$8-c$}}
         \put(-2, 11){\tiny{$8-c$}}
         \put(5, 20){$S$}
         \put(15, 19){\tiny{$c$}}
         \put(24, 20){$S$}
         \put(15, 3){$T$}
        \end{overpic}
    \end{equation}
    \begin{equation}\label{eq: twoequivalentdiagramsuntwistwithf}
         \begin{overpic}[unit=1mm, scale=0.6]{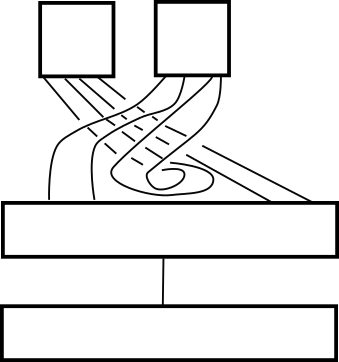}
        \put(10, 50){$S$}
        \put(29, 50){$S$}
        \put(10, 29){\tiny{$...$}}
        \put(24,3){$T$}
        \put(20, 29){\tiny{$...$}}
        \put(19, 30){\tiny{$(c)$}}
        \put(39,29){\tiny{$...$}}
         \put(15, 12){\tiny{$16-2c$}}
        \put(18, 20){{$f^{(16-2c)}$}}
    \end{overpic}\hspace{1mm},\hspace{4mm}
        \begin{overpic}[unit=1mm, scale=0.6]{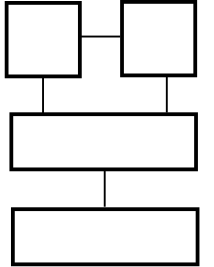}
         \put(14, 4){$T$}
         \put(5, 35){$S$}
         \put(24, 35){$S$}
         \put(-2,26){\tiny{$8-c$}}
        \put(15, 34){\tiny{$c$}}
        \put(20, 12){\tiny{$16-2c$}}
         \put(28,26){\tiny{$8-c$}}
         \put(15, 19){{$f^{(16-2c)}$}}
        \end{overpic}
    \end{equation}
\end{lem}

\begin{proof}
    The forward direction is clear. For the backwards direction, suppose that the jellyfish algorithm evaluates the two diagrams in (\ref{eq: twoequivalentdiagramsuntwistwithf}) to be the same numbers. First suppose $c=8$. In this case the result is clear since the Jones-Wenzl projection is $f^{(0)}=\emptyset$. So the diagrams in (\ref{eq: twoequivalentdiagramsuntwistwithf}) are equivalent to the diagrams in (\ref{eq: twoequivalentdiagramsuntwist}). 

    Suppose the result is true down to some $c+1$. Consider the case with $c$ twisted strands and suppose the diagrams in (\ref{eq: twoequivalentdiagramsuntwistwithf}) evaluate to the same number in the jellyfish algorithm. Expand $f^{(16-2c)}=X^{\otimes (16-2c)}+D$, where $D$ is a linear combination of diagrams with a cup and cap somewhere. For each diagram, since $S$ is U.C.C.S., the only nonzero terms in this expansion will have a cup in the $8-c$ and $8-c+1$ position, which are a sum of terms from the same diagram with $c+k$ strands in the middle where $k\geq 1$. The jellyfish algorithm is clearly linear, so by induction, the equivalence from the $c+k$ case for $k\geq 1$ will imply the equivalence of terms in the $c$ case. 
\end{proof}

A symmetric argument will show the following. 

\begin{lem}\label{lem: twistanduntwistthesame2}
     Fix $c\in\{0,...,8\}$ and $T\in \te{Hom}(\emptyset, X^{\otimes 2(8-c)})$. Assume the jellyfish algorithm has been completed up to step 5 on two diagrams which only differ as shown below in (\ref{eqn: twistotherway2}). Similarly, assume the jellyfish algorithm has been completed up to step 5 on two diagrams which only differ as shown below in (\ref{eqn: twistotherway}). Then following steps 6 and 7 on each of the diagrams in (\ref{eqn: twistotherway2}) results in the same number if and only if following steps 6 and 7 on each of the diagrams in (\ref{eqn: twistotherway}) results in the same number.
    \begin{equation}\label{eqn: twistotherway2}
         \begin{overpic}[unit=1mm, scale=0.6]{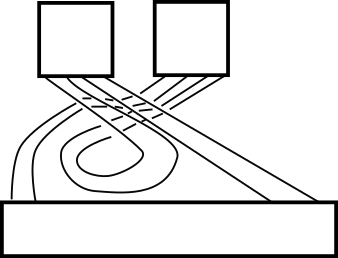}
        \put(10, 34){$S$}
        \put(29, 34){$S$}
        \put(2.5, 11){\tiny{$...$}}
        \put(24,3){$T$}
        \put(23, 15){\tiny{$...$}}
        \put(23.5, 17){\tiny{$(c)$}}
        \put(39,13){\tiny{$...$}}
    \end{overpic}\hspace{1mm},\hspace{4mm}
       \begin{overpic}[unit=1mm, scale=0.6]{twoboxesconnectedoneboxbelowclosed.png}
         \put(28,12){\tiny{$8-c$}}
         \put(-2, 11){\tiny{$8-c$}}
         \put(5, 20){$S$}
         \put(15, 19){\tiny{$c$}}
         \put(24, 20){$S$}
         \put(15, 3){$T$}
        \end{overpic}\hspace{1mm}.
    \end{equation}

    \begin{equation}\label{eqn: twistotherway}
         \begin{overpic}[unit=1mm, scale=0.6]{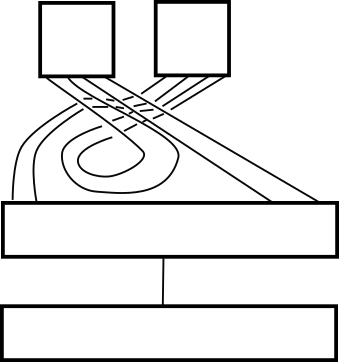}
        \put(10, 50){$S$}
        \put(29, 50){$S$}
        \put(2.5, 27){\tiny{$...$}}
        \put(24,3){$T$}
        \put(23, 31){\tiny{$...$}}
        \put(23.5, 33){\tiny{$(c)$}}
        \put(39,29){\tiny{$...$}}
         \put(14, 12){\tiny{$16-2c$}}
            \put(18, 20){{$f^{(16-2c)}$}}
    \end{overpic}\hspace{1mm},\hspace{4mm}
        \begin{overpic}[unit=1mm, scale=0.6]{fourboxestwoontop2bigonbottom.png}
         \put(14, 4){$T$}
         \put(5, 35){$S$}
         \put(24, 35){$S$}
         \put(-2,26){\tiny{$8-c$}}
        \put(15, 34){\tiny{$c$}}
        \put(20, 12){\tiny{$16-2c$}}
         \put(28,26){\tiny{$8-c$}}
            \put(15, 19){{$f^{(16-2c)}$}}
        \end{overpic}
    \end{equation}
\end{lem}

We now focus on proving that the jellyfish algorithm evaluates the diagrams in (\ref{eq: twoequivalentdiagramsuntwistwithf}) and (\ref{eqn: twistotherway}) equivalently by going through all the possible cases. 

\begin{lem}\label{lem: e7sufficientfourboxestwoontop2bigonbottomiszero}
    Assume the jellyfish algorithm has been completed up through step 5 (so the ``circular" diagram may have crossings in the middle). If at this point, the diagram is of the form:
    
    \begin{equation}\label{eq: fourboxestwoontop2bigonbottom}
        \begin{overpic}[unit=1mm, scale=0.6]{fourboxestwoontop2bigonbottom.png}
         \put(14, 4){$T$}
         \put(5, 35){$S$}
         \put(24, 35){$S$}
         \put(3,26){\tiny{$5$}}
        \put(15, 34){\tiny{$3$}}
        \put(18, 11){\tiny{$10$}}
         \put(28,26){\tiny{$5$}}
        \put(14, 19){$f^{(10)}$}
        \end{overpic}
    \end{equation}  

    \noindent where $T\in \te{Hom}(\emptyset, X^{\otimes 10})$, and $T$ may have crossings, then the jellyfish algorithm will evaluate this diagram to zero. 
\end{lem}

\begin{proof}
    We begin with step 6 of the jellyfish algorithm. It is well known in knot theory that the order in which we resolve crossings does not matter. So suppose all crossings of $T$ are now resolved and without loss of generality redefine $T$ to be any such summand in this resolution. As $T\in \te{Hom}(\emptyset, f^{(10)})$, we know that if $T\in \mc{TL}$ then it would cup $f^{(10)}$. Steps 7(i) through (iii) will not change the cupped $f^{(10)}$, so in step 7(iv), we get that the evaluation is zero. 

    Now suppose that $T$ has exactly one $S$-box. As $f^{(10)}$ and $T$ are connected by 10 strands, there still must be a cup somewhere on $f^{(10)}$ and again the diagram will evaluate to zero. 

    Suppose $T$ has exactly two $S$-boxes. Using Lemma \ref{lem: e7sufficientwelldefinedSsquared}, we can perform steps 7(i) through (iii) for the $S$-boxes in $T$. If there is some $S^2$, we reduce to the previous case. Otherwise, we can assume that $T$ has two $S$-boxes are connected by at most 3 strands. Then, as there are 10 points on the boundary, either some $S$-box is capped and thus by step 7(i) will evaluate the diagram to zero, or $T$ is, up to multiplication by a scalar, the below diagram,
     \begin{equation}\label{eq: twoboxesconnectedup}
        \begin{overpic}[unit=1mm, scale=0.6]{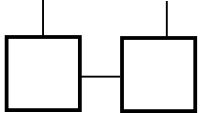}
         \put(24, 5){$S$}
         \put(15, 8){\tiny{$3$}}
         \put(5, 5){$S$}
        \end{overpic}. 
    \end{equation}  
    Suppose that $T$ has $n$ $S$-boxes. By Lemma \ref{lem: e7sufficientwelldefinedSsquared}, we can assume that we have used steps 7(i) and (ii) so that no $S^2$-boxes are remaining inside $T$.  Then the leftmost $S$ box has at most 3 strands connected to the adjacent $S$-box and at least 5 strands on the boundary. Similarly, the rightmost $S$-box has at most 3 strands connected to its neighbor and at least 5 strands on the boundary. Therefore, with the leftmost and rightmost $S$-boxes, we have accounted for all the strands on the boundary of $T$. If there are other $S$-boxes between the rightmost and leftmost $S$-boxes, there must be one connected only to its neighbors, forming an $S^2$, which contradicts that no $S^2$ remains in $T$. Therefore, the only case to consider is when $T$ is a scalar multiple of the diagram in (\ref{eq: twoboxesconnectedup}). We compute the diagram in (\ref{eq: fourboxestwoontop2bigonbottom}) using $T$ equal to the diagram in (\ref{eq: twoboxesconnectedup}).

    We use Chen's formula from Propostion \ref{prop: chen} to get a formula for $f^{(5)}\otimes f^{(5)}$. Using the $\theta$-values found in Lemma \ref{lem: thetavalues55} we obtain
     \begin{equation}\label{eqn: f5tensorf5expansion}
        f^{(5)}\otimes f^{(5)}= \f{1}{6} \hspace{1mm}{\raisebox{-0.5 cm}{\begin{overpic}[unit=1mm, scale=0.5]{twotrivalentvertices.png}
         \put(1, 11){\tiny{$5$}}
          \put(8, 11){\tiny{$5$}}
          \put(1, 3){\tiny{$5$}}
         \put(6, 7){\tiny{$0$}}
           \put(8, 3){\tiny{$5$}}
        \end{overpic}}}
        +\f{5}{7}
        {\raisebox{-0.5 cm}{\begin{overpic}[unit=1mm, scale=0.5]{twotrivalentvertices.png}
         \put(1, 11){\tiny{$5$}}
          \put(8, 11){\tiny{$5$}}
          \put(1, 3){\tiny{$5$}}
         \put(6, 7){\tiny{$2$}}
           \put(8, 3){\tiny{$5$}}
        \end{overpic}}}
        +\f{25}{14}
        {\raisebox{-0.5 cm}{\begin{overpic}[unit=1mm, scale=0.5]{twotrivalentvertices.png}
         \put(1, 11){\tiny{$5$}}
          \put(8, 11){\tiny{$5$}}
          \put(1, 3){\tiny{$5$}}
         \put(6, 7){\tiny{$4$}}
           \put(8, 3){\tiny{$5$}}
        \end{overpic}}}
        +\f{25}{9}
           {\raisebox{-0.5 cm}{\begin{overpic}[unit=1mm, scale=0.5]{twotrivalentvertices.png}
         \put(1, 11){\tiny{$5$}}
          \put(8, 11){\tiny{$5$}}
          \put(1, 3){\tiny{$5$}}
         \put(6, 7){\tiny{$6$}}
           \put(8, 3){\tiny{$5$}}
        \end{overpic}}}
        +\f{5}{2}
           {\raisebox{-0.5 cm}{\begin{overpic}[unit=1mm, scale=0.5]{twotrivalentvertices.png}
         \put(1, 11){\tiny{$5$}}
          \put(8, 11){\tiny{$5$}}
          \put(1, 3){\tiny{$5$}}
         \put(6, 7){\tiny{$8$}}
           \put(8, 3){\tiny{$5$}}
        \end{overpic}}}
        +
           {\raisebox{-0.5 cm}{\begin{overpic}[unit=1mm, scale=0.5]{twotrivalentvertices.png}
         \put(1, 11){\tiny{$5$}}
          \put(8, 11){\tiny{$5$}}
          \put(1, 3){\tiny{$5$}}
         \put(6, 7){\tiny{$10$}}
           \put(8, 3){\tiny{$5$}}
        \end{overpic}}}
    \end{equation}
    Further, notice that
    \begin{align*}
           {\raisebox{-0.5 cm}{\begin{overpic}[unit=1mm, scale=0.5]{twotrivalentvertices.png}
         \put(1, 11){\tiny{$5$}}
          \put(8, 11){\tiny{$5$}}
          \put(1, 3){\tiny{$5$}}
         \put(6, 7){\tiny{$10$}}
           \put(8, 3){\tiny{$5$}}
        \end{overpic}}}
        = 
        {\raisebox{-1.5 cm}{
         \begin{overpic}[unit=1mm, scale=0.4]{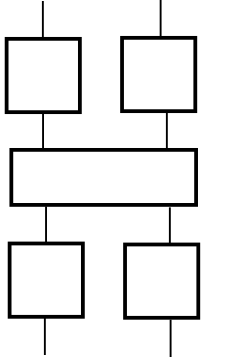}
         \put(16,1){\tiny{$5$}}
         \put(2.5, 7){$f^{(5)}$}
         \put(2.5, 29){$f^{(5)}$}
         \put(14.5, 29){$f^{(5)}$}
         \put(14.5, 7){$f^{(5)}$}
         \put(9, 18){$f^{(10)}$}
         \put(2, 13){\tiny{$5$}}
         \put(16, 13){\tiny{$5$}}
         \put(3, 23){\tiny{$5$}}
         \put(16, 23){\tiny{$5$}}
         \put(2,1){\tiny{$5$}}
          \put(15, 35){\tiny{$5$}}
         \put(2,35){\tiny{$5$}}
        \end{overpic}}}
        =f^{(10)}
    \end{align*}
    since $f^{(10)}$ is uncuppable and uncappable. Replace $f^{(10)}$ in the middle of (\ref{eq: fourboxestwoontop2bigonbottom}) using equation \ref{eqn: f5tensorf5expansion}. As all of equation \ref{eqn: f5tensorf5expansion} is in Temperley-Lieb, the jellyfish algorithm will not be affected. Then continue steps 7(i) through (iii) of the algorithm. Using that $S$ is U.C.C.S. gives the diagram (\ref{eq: fourboxestwoontop2bigonbottom}) will equal
    \begin{align*}
        {\scalebox{0.6}{%
         \begin{overpic}[unit=1mm, scale=0.6]{fourboxesconnected.png}
         \put(24,4){$S$}
         \put(5, 4){$S$}
         \put(5, 22){$S$}
         \put(16, 25){\tiny{$3$}}
         \put(24, 22){$S$}
         \put(3,14){\tiny{$5$}}
         \put(27,14){\tiny{$5$}}
         \put(16, 1){\tiny{$3$}}
        \end{overpic}
       \hspace{2mm}{\raisebox{1.2cm}{$-\f{1}{6}$}}
       \begin{overpic}[unit=1mm, scale=0.6]{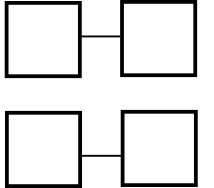}
         \put(24,4){$S$}
         \put(5, 4){$S$}
         \put(5, 22){$S$}
         \put(16, 25){\tiny{$8$}}
         \put(24, 22){$S$}
         \put(16, 1){\tiny{$8$}}
        \end{overpic}
        \hspace{2mm}{\raisebox{1.2cm}{$-\f{5}{7}$}}
        {\raisebox{-0.5cm}{
         \begin{overpic}[unit=1mm, scale=0.6]{fourboxesoneinmiddle.png}
         \put(6, 5){$S$}
         \put(6, 38){$S$}
             \put(24,5){$S$}
         \put(24, 38){$S$}
         \put(3,14){\tiny{$1$}}
            \put(3, 30){\tiny{$1$}}
         \put(16, 1){\tiny{$7$}}
          \put(16, 40){\tiny{$7$}}
       \put(28,14){\tiny{$1$}}
            \put(28, 30){\tiny{$1$}}
            \put(13, 21){$f^{(2)}$}
        \end{overpic}}}
        \hspace{2mm}{\raisebox{1.2cm}{$-\f{25}{14}$}}
        {\raisebox{-0.5cm}{
         \begin{overpic}[unit=1mm, scale=0.6]{fourboxesoneinmiddle.png}
         \put(6, 5){$S$}
         \put(6, 38){$S$}
             \put(24,5){$S$}
         \put(24, 38){$S$}
         \put(3,14){\tiny{$2$}}
            \put(3, 30){\tiny{$2$}}
         \put(16, 1){\tiny{$6$}}
          \put(16, 40){\tiny{$6$}}
       \put(28,14){\tiny{$2$}}
            \put(28, 30){\tiny{$2$}}
            \put(13, 21){$f^{(4)}$}
        \end{overpic}}}
             \hspace{2mm}{\raisebox{1.2cm}{$-\f{25}{9}$}}
        {\raisebox{-0.5cm}{
         \begin{overpic}[unit=1mm, scale=0.6]{fourboxesoneinmiddle.png}
         \put(6, 5){$S$}
         \put(6, 38){$S$}
             \put(24,5){$S$}
         \put(24, 38){$S$}
         \put(3,14){\tiny{$3$}}
            \put(3, 30){\tiny{$3$}}
         \put(16, 1){\tiny{$5$}}
          \put(16, 40){\tiny{$5$}}
       \put(28,14){\tiny{$3$}}
            \put(28, 30){\tiny{$3$}}
            \put(13, 21){$f^{(6)}$}
        \end{overpic}}}
         \hspace{2mm}{\raisebox{1.2cm}{$-\f{5}{2}$}}
        {\raisebox{-0.5cm}{
         \begin{overpic}[unit=1mm, scale=0.6]{fourboxesoneinmiddle.png}
         \put(6, 5){$S$}
         \put(6, 38){$S$}
             \put(24,5){$S$}
         \put(24, 38){$S$}
         \put(3,14){\tiny{$4$}}
            \put(3, 30){\tiny{$4$}}
         \put(16, 1){\tiny{$4$}}
          \put(16, 40){\tiny{$4$}}
       \put(28,14){\tiny{$4$}}
            \put(28, 30){\tiny{$4$}}
            \put(13, 21){$f^{(8)}$}
        \end{overpic}}}}}
    \end{align*}
    The choice of reducing $S^2$ doesn't matter, so from here we can use the evaluations found in Lemmas \ref{lem: fourboxesconnected}, \ref{lem: e7sufficienttraceofSsquared}, \ref{lem: twoboxesconnectedoneboxbelow}, and \ref{lem: fourboxesoneinmiddle} to evaluate the above to $225-\f{1}{6}\cdot 30^2-\f{5}{2}\cdot 30=0$. 
\end{proof}

\begin{lem}\label{lem: e7sufficienttwoSsonf14onlypossibilities}
    Suppose the jellyfish algorithm has been completed up through step 6 of the algorithm. Then suppose step 7 has been repeated to get rid of all $S^2$ boxes across the top and the bottom so that at this point it looks like the below diagram:
    \begin{equation}\label{eqn: twoSsonf14}
         \begin{overpic}[unit=1mm, scale=0.6]{fourboxestwoontop2bigonbottom.png}
         \put(14, 4){$T$}
         \put(5, 35){$S$}
         \put(24, 35){$S$}
         \put(3,26){\tiny{$7$}}
        \put(15, 34){\tiny{$1$}}
        \put(18, 11){\tiny{$14$}}
         \put(28,26){\tiny{$7$}}
        \put(14, 19){$f^{(14)}$}
        \end{overpic}
    \end{equation}
    where $T\in \te{Hom}(\emptyset, f^{(14)})$. If the diagram is nonzero, then, up to multiplication by a scalar, $T$ is of one of the following:
    \begin{enumerate}
        \item {\raisebox{-0.2cm}{$\begin{aligned}[t]
             \begin{overpic}[unit=1mm, scale=0.6]{twoboxesconnectedup.png}
         \put(23, 5){$S$}
         \put(15, 8){\tiny{$1$}}
         \put(5, 5){$S$}
        \end{overpic}
        \end{aligned}$}}
        \item  {\raisebox{-0.2cm}{$\begin{aligned}[t]
             \begin{overpic}[unit=1mm, scale=0.6]{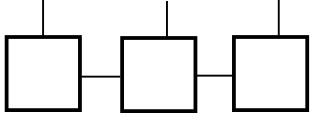}
         \put(23, 5){$S$}
         \put(15, 8){\tiny{$3$}}
         \put(5, 5){$S$}
         \put(33, 8){\tiny{$2$}}
         \put(42, 5){$S$}
        \end{overpic}
        \end{aligned}$}}
        \item  {\raisebox{-0.2cm}{$\begin{aligned}[t]
             \begin{overpic}[unit=1mm, scale=0.6]{3boxesconnectedstrandsup.png}
         \put(23, 5){$S$}
         \put(15, 8){\tiny{$2$}}
         \put(5, 5){$S$}
         \put(33, 8){\tiny{$3$}}
         \put(42, 5){$S$}
        \end{overpic}
        \end{aligned}$}}
         \item  {\raisebox{-0.2cm}{$\begin{aligned}[t]
             \begin{overpic}[unit=1mm, scale=0.6]{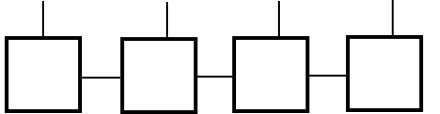}
         \put(23, 5){$S$}
         \put(15, 8){\tiny{$3$}}
         \put(5, 5){$S$}
         \put(33, 8){\tiny{$3$}}
         \put(42, 5){$S$}
         \put(60, 5){$S$}
         \put(51, 8){\tiny{$3$}}
        \end{overpic}
        \end{aligned}$}}
    \end{enumerate}
\end{lem}

\begin{proof}
    If the diagram is nonzero, then $f^{(14)}$ cannot be cupped. This then requires that all 14 points on the bottom boundary of $f^{(14)}$ are connected to an $S$-box. Then the lemma is clear as the four cases given in the lemma are the only choices to have $S$-boxes connected by less than 4 strands to each other and with 14 strands on the boundary. 
\end{proof}

\begin{lem}\label{lem: fourboxesoneinmiddle3750over77}
    The jellyfish algorithm evaluates the below diagram to $\f{3750}{77}$, which is its value in $\mc{P}$.
     \begin{equation}\label{eq: fourboxesoneinmiddle3750over77}
        \begin{overpic}[unit=1mm, scale=0.6]{fourboxesoneinmiddle.png}
         \put(24,5){$S$}
         \put(6, 5){$S$}
         \put(6, 38){$S$}
         \put(16, 40){\tiny{$2$}}
         \put(24, 38){$S$}
         \put(3,14){\tiny{$6$}}
         \put(28,14){\tiny{$6$}}
         \put(16, 1){\tiny{$2$}}
          \put(3, 30){\tiny{$6$}}
            \put(28, 30){\tiny{$6$}}
            \put(14, 21){$f^{(12)}$}
        \end{overpic}
    \end{equation}  
\end{lem}

\begin{proof}
    It is clear that steps 1-6 have already been completed for this diagram in the jellyfish algorithm. Further, the $S^2$-boxes have been reduced along the top and the bottom. It is then clear that evaluating this diagram in $\mc{P}$ will equal its evaluation in the jellyfish algorithm. By Chen's formula in Proposition \ref{prop: chen} and the $\theta$-values found in Lemma \ref{lem: thetavalues66}, we obtain
    \begin{equation}\label{eqn: f6tensorf6expansion}
        f^{(6)}\otimes f^{(6)}= \f{1}{7} \hspace{1mm}{\raisebox{-0.5 cm}{\begin{overpic}[unit=1mm, scale=0.5]{twotrivalentvertices.png}
         \put(1, 11){\tiny{$6$}}
          \put(8, 11){\tiny{$6$}}
          \put(1, 3){\tiny{$6$}}
         \put(6, 7){\tiny{$0$}}
           \put(8, 3){\tiny{$6$}}
        \end{overpic}}}
        +\f{9}{14}
        {\raisebox{-0.5 cm}{\begin{overpic}[unit=1mm, scale=0.5]{twotrivalentvertices.png}
         \put(1, 11){\tiny{$6$}}
          \put(8, 11){\tiny{$6$}}
          \put(1, 3){\tiny{$6$}}
         \put(6, 7){\tiny{$2$}}
           \put(8, 3){\tiny{$6$}}
        \end{overpic}}}
        +\f{25}{14}
        {\raisebox{-0.5 cm}{\begin{overpic}[unit=1mm, scale=0.5]{twotrivalentvertices.png}
         \put(1, 11){\tiny{$6$}}
          \put(8, 11){\tiny{$6$}}
          \put(1, 3){\tiny{$6$}}
         \put(6, 7){\tiny{$4$}}
           \put(8, 3){\tiny{$6$}}
        \end{overpic}}}
        +\f{10}{3}
           {\raisebox{-0.5 cm}{\begin{overpic}[unit=1mm, scale=0.5]{twotrivalentvertices.png}
         \put(1, 11){\tiny{$6$}}
          \put(8, 11){\tiny{$6$}}
          \put(1, 3){\tiny{$6$}}
         \put(6, 7){\tiny{$6$}}
           \put(8, 3){\tiny{$6$}}
        \end{overpic}}}
        +\f{45}{11}
           {\raisebox{-0.5 cm}{\begin{overpic}[unit=1mm, scale=0.5]{twotrivalentvertices.png}
         \put(1, 11){\tiny{$6$}}
          \put(8, 11){\tiny{$6$}}
          \put(1, 3){\tiny{$6$}}
         \put(6, 7){\tiny{$8$}}
           \put(8, 3){\tiny{$6$}}
        \end{overpic}}}
        +3
           {\raisebox{-0.5 cm}{\begin{overpic}[unit=1mm, scale=0.5]{twotrivalentvertices.png}
         \put(1, 11){\tiny{$6$}}
          \put(8, 11){\tiny{$6$}}
          \put(1, 3){\tiny{$6$}}
         \put(6, 7){\tiny{$10$}}
           \put(8, 3){\tiny{$6$}}
        \end{overpic}}}
        +
        {\raisebox{-0.5 cm}{\begin{overpic}[unit=1mm, scale=0.5]{twotrivalentvertices.png}
         \put(1, 11){\tiny{$6$}}
          \put(8, 11){\tiny{$6$}}
          \put(1, 3){\tiny{$6$}}
         \put(6, 7){\tiny{$12$}}
           \put(8, 3){\tiny{$6$}}
        \end{overpic}}}
    \end{equation}
    Further notice that 
     \begin{align*}
           {\raisebox{-0.5 cm}{\begin{overpic}[unit=1mm, scale=0.5]{twotrivalentvertices.png}
         \put(1, 11){\tiny{$6$}}
          \put(8, 11){\tiny{$6$}}
          \put(1, 3){\tiny{$6$}}
         \put(6, 7){\tiny{$12$}}
           \put(8, 3){\tiny{$6$}}
        \end{overpic}}}
        = 
        {\raisebox{-1.5 cm}{
         \begin{overpic}[unit=1mm, scale=0.4]{fourboxesoneinmiddlestrandsout.png}
         \put(16,1){\tiny{$6$}}
         \put(2, 13){\tiny{$6$}}
         \put(16, 13){\tiny{$6$}}
         \put(3, 23){\tiny{$6$}}
         \put(16, 23){\tiny{$6$}}
         \put(2,1){\tiny{$6$}}
          \put(15, 35){\tiny{$6$}}
         \put(2,35){\tiny{$6$}}
           \put(2.5, 7){$f^{(6)}$}
         \put(2.5, 29){$f^{(6)}$}
         \put(14.5, 29){$f^{(6)}$}
         \put(14.5, 7){$f^{(6)}$}
         \put(9, 18){$f^{(12)}$}
        \end{overpic}}}
        =f^{(12)}. 
        \end{align*}
        Replacing $f^{(12)}$ in the diagram (\ref{eq: fourboxesoneinmiddle3750over77}) using equation \ref{eqn: f6tensorf6expansion} and using the $S$ is U.C.C.S. gives that the diagram (\ref{eq: fourboxesoneinmiddle3750over77}) equals
         \begin{align*}
         &{\scalebox{0.7}{%
         \begin{overpic}[unit=1mm, scale=0.6]{fourboxesconnected.png}
         \put(23,4){$S$}
         \put(5, 4){$S$}
         \put(5, 22){$S$}
         \put(16, 25){\tiny{$2$}}
         \put(24, 22){$S$}
         \put(3,14){\tiny{$6$}}
         \put(27,14){\tiny{$6$}}
         \put(16, 1){\tiny{$2$}}
        \end{overpic}
       \hspace{2mm}{\raisebox{1.2cm}{$-\f{1}{7}$}}
       \begin{overpic}[unit=1mm, scale=0.6]{twopairsofboxes.png}
         \put(24,4){$S$}
         \put(5, 4){$S$}
         \put(5, 22){$S$}
         \put(16, 25){\tiny{$8$}}
         \put(23, 22){$S$}
         \put(16, 1){\tiny{$8$}}
        \end{overpic}
        \hspace{2mm}{\raisebox{1.2cm}{$-\f{9}{14}$}}
        {\raisebox{-0.5cm}{
         \begin{overpic}[unit=1mm, scale=0.6]{fourboxesoneinmiddle.png}
         \put(6, 5){$S$}
         \put(6, 38){$S$}
             \put(24,5){$S$}
         \put(24, 38){$S$}
         \put(3,14){\tiny{$1$}}
            \put(3, 30){\tiny{$1$}}
         \put(16, 1){\tiny{$7$}}
          \put(16, 40){\tiny{$7$}}
       \put(28,14){\tiny{$1$}}
            \put(28, 30){\tiny{$1$}}
            \put(13, 21){$f^{(2)}$}
        \end{overpic}}}
        \hspace{2mm}{\raisebox{1.2cm}{$-\f{25}{14}$}}
        {\raisebox{-0.5cm}{
         \begin{overpic}[unit=1mm, scale=0.6]{fourboxesoneinmiddle.png}
         \put(6, 5){$S$}
         \put(6, 38){$S$}
             \put(24,5){$S$}
         \put(24, 38){$S$}
         \put(3,14){\tiny{$2$}}
            \put(3, 30){\tiny{$2$}}
         \put(16, 1){\tiny{$6$}}
          \put(16, 40){\tiny{$6$}}
       \put(28,14){\tiny{$2$}}
            \put(28, 30){\tiny{$2$}}
            \put(13, 21){$f^{(4)}$}
        \end{overpic}}}}}\\
      & {\scalebox{0.7}{%  
       \hspace{2mm}{\raisebox{1.2cm}{$-\f{10}{3}$}}
      {\raisebox{-0.5cm}{
         \begin{overpic}[unit=1mm, scale=0.6]{fourboxesoneinmiddle.png}
         \put(6, 5){$S$}
         \put(6, 38){$S$}
             \put(24,5){$S$}
         \put(24, 38){$S$}
         \put(3,14){\tiny{$3$}}
            \put(3, 30){\tiny{$3$}}
         \put(16, 1){\tiny{$5$}}
          \put(16, 40){\tiny{$5$}}
       \put(28,14){\tiny{$3$}}
            \put(28, 30){\tiny{$3$}}
            \put(13, 21){$f^{(6)}$}
        \end{overpic}}}
       \hspace{2mm}{\raisebox{1.2cm}{$-\f{45}{11}$}}
        {\raisebox{-0.5cm}{
         \begin{overpic}[unit=1mm, scale=0.6]{fourboxesoneinmiddle.png}
         \put(6, 5){$S$}
         \put(6, 38){$S$}
             \put(24,5){$S$}
         \put(24, 38){$S$}
         \put(3,14){\tiny{$4$}}
            \put(3, 30){\tiny{$4$}}
         \put(16, 1){\tiny{$4$}}
          \put(16, 40){\tiny{$4$}}
       \put(28,14){\tiny{$4$}}
            \put(28, 30){\tiny{$4$}}
            \put(13, 21){$f^{(8)}$}
        \end{overpic}}}
               \hspace{2mm}{\raisebox{1.2cm}{$-3$}}
        {\raisebox{-0.5cm}{
         \begin{overpic}[unit=1mm, scale=0.6]{fourboxesoneinmiddle.png}
         \put(6, 5){$S$}
         \put(6, 38){$S$}
             \put(24,5){$S$}
         \put(24, 38){$S$}
         \put(3,14){\tiny{$5$}}
            \put(3, 30){\tiny{$5$}}
         \put(16, 1){\tiny{$3$}}
          \put(16, 40){\tiny{$3$}}
       \put(28,14){\tiny{$5$}}
            \put(28, 30){\tiny{$5$}}
            \put(13, 21){$f^{(10)}$}
        \end{overpic}}}}}.
    \end{align*}
    The third, fourth, and fifth term are zero by Lemma \ref{lem: twoboxesconnectedoneboxbelow}. Lemma \ref{lem: e7sufficientfourboxestwoontop2bigonbottomiszero} gives that the last term is zero. The second term is two copies of $S^2$. Lemma \ref{lem: fourboxesconnectedequaling300} evaluates the first term to 300. Lemma \ref{lem: fourboxesoneinmiddle} evaluates the sixth term to 30. Therefore, diagram \ref{eq: fourboxesoneinmiddle3750over77} equals $300-\f{1}{7}\cdot 30^2-\f{45}{11}\cdot 30=\f{3750}{77}$.
\end{proof}

\begin{lem}\label{lem: e7sufficientexpandSSontopoff14}
    Suppose the jellyfish algorithm has been completed up through step 6 of the algorithm. Then suppose step 7 has been repeated to get rid of all $S^2$ boxes across the top and the bottom so that at this point it looks like the below diagram:
    \begin{align*}
         \begin{overpic}[unit=1mm, scale=0.6]{fourboxestwoontop2bigonbottom.png}
         \put(14, 4){$T$}
         \put(5, 35){$S$}
         \put(24, 35){$S$}
         \put(3,26){\tiny{$7$}}
        \put(15, 34){\tiny{$1$}}
        \put(18, 11){\tiny{$14$}}
         \put(28,26){\tiny{$7$}}
        \put(14, 19){$f^{(14)}$}
        \end{overpic}
    \end{align*}
    where $T\in \te{Hom}(\emptyset, f^{(14)})$. Then this diagram evaluates under the jellyfish algorithm to the same as the following sum:
    \begin{equation}\label{eq: expansionSSontopoff14}
       \scalebox{0.9}{  \begin{overpic}[unit=1mm, scale=0.6]{twoboxesconnectedoneboxbelowclosed.png}
         \put(28,12){\tiny{$7$}}
         \put(2, 11){\tiny{$7$}}
         \put(4, 20){$S$}
         \put(15, 19){\tiny{$1$}}
         \put(24, 20){$S$}
         \put(15, 3){$T$}
        \end{overpic}
        {\raisebox{1.5cm}{$-\f{30}{8}$}}
        \begin{overpic}[unit=1mm, scale=0.6]{twoboxesstacked.png}
         \put(17, 11){\tiny{$14$}}
         \put(13, 4){$T$}
         \put(12.5, 20){$f^{(7)}$}
        \end{overpic}
         {\raisebox{1.6cm}{$-\f{245}{44}$}}
         {\raisebox{-2cm}{ \begin{overpic}[unit=1mm, scale=0.6]{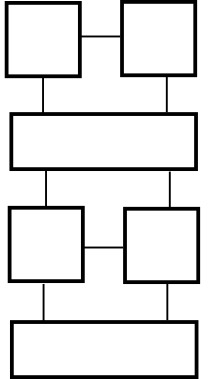}
         \put(15, 3){$T$}
         \put(5, 53){$S$}
        \put(24, 53){$S$}
         \put(16, 51){\tiny{$4$}}
         \put(4,11){\tiny{$7$}}
         \put(28,11){\tiny{$7$}}
         \put(16, 18){\tiny{$3$}}
          \put(4, 30){\tiny{$4$}}
            \put(28, 30){\tiny{$4$}}
              \put(4, 45){\tiny{$4$}}
            \put(28, 45){\tiny{$4$}}
            \put(23,20){$f^{(7)}$}
            \put(15,37){$f^{(8)}$}
            \put(5,20){$f^{(7)}$}
        \end{overpic}}}
         {\raisebox{1.6cm}{$-\f{7}{2}$}}
         {\raisebox{-2cm}{ \begin{overpic}[unit=1mm, scale=0.6]{sixboxes.png}
         \put(15, 3){$T$}
         \put(5, 53){$S$}
        \put(24, 53){$S$}
         \put(16, 51){\tiny{$2$}}
         \put(4,11){\tiny{$7$}}
         \put(28,11){\tiny{$7$}}
         \put(16, 18){\tiny{$1$}}
          \put(4, 30){\tiny{$6$}}
            \put(28, 30){\tiny{$6$}}
              \put(4, 45){\tiny{$6$}}
            \put(28, 45){\tiny{$6$}}
                 \put(23,20){$f^{(7)}$}
            \put(15,37){$f^{(12)}$}
            \put(5,20){$f^{(7)}$}
        \end{overpic}}}}
    \end{equation}
\end{lem}

\begin{proof}
    As done in the previous lemmas, it is enough to reduce using Chen's formula in Proposition \ref{prop: chen} and the $\theta$-values found in Lemma \ref{lem: thetavalues77}. This gives
    \begin{equation}\label{eqn: f7tensorf7expansion}
       f^{(7)}\otimes f^{(7)}= \f{1}{8} \hspace{1mm} {\raisebox{-0.5 cm}{\scalebox{0.9}{\begin{overpic}[unit=1mm, scale=0.5]{twotrivalentvertices.png}
         \put(1, 11){\tiny{$7$}}
          \put(8, 11){\tiny{$7$}}
          \put(1, 3){\tiny{$7$}}
         \put(6, 7){\tiny{$0$}}
           \put(8, 3){\tiny{$7$}}
        \end{overpic}}}}
        +\f{7}{12}
        {\raisebox{-0.5 cm}{\scalebox{0.9}{\begin{overpic}[unit=1mm, scale=0.5]{twotrivalentvertices.png}
         \put(1, 11){\tiny{$7$}}
          \put(8, 11){\tiny{$7$}}
          \put(1, 3){\tiny{$7$}}
         \put(6, 7){\tiny{$2$}}
           \put(8, 3){\tiny{$7$}}
        \end{overpic}}}}
        +\f{7}{4}
        {\raisebox{-0.5 cm}{\scalebox{0.9}{\begin{overpic}[unit=1mm, scale=0.5]{twotrivalentvertices.png}
         \put(1, 11){\tiny{$7$}}
          \put(8, 11){\tiny{$7$}}
          \put(1, 3){\tiny{$7$}}
         \put(6, 7){\tiny{$4$}}
           \put(8, 3){\tiny{$7$}}
        \end{overpic}}}}
        +\f{245}{66}
           {\raisebox{-0.5 cm}{\scalebox{0.9}{\begin{overpic}[unit=1mm, scale=0.5]{twotrivalentvertices.png}
         \put(1, 11){\tiny{$7$}}
          \put(8, 11){\tiny{$7$}}
          \put(1, 3){\tiny{$7$}}
         \put(6, 7){\tiny{$6$}}
           \put(8, 3){\tiny{$7$}}
        \end{overpic}}}}
        +\f{245}{44}
           {\raisebox{-0.5 cm}{\scalebox{0.9}{\begin{overpic}[unit=1mm, scale=0.5]{twotrivalentvertices.png}
         \put(1, 11){\tiny{$7$}}
          \put(8, 11){\tiny{$7$}}
          \put(1, 3){\tiny{$7$}}
         \put(6, 7){\tiny{$8$}}
           \put(8, 3){\tiny{$7$}}
        \end{overpic}}}}
        +\f{147}{26}
           {\raisebox{-0.5 cm}{\scalebox{0.9}{\begin{overpic}[unit=1mm, scale=0.5]{twotrivalentvertices.png}
         \put(1, 11){\tiny{$7$}}
          \put(8, 11){\tiny{$7$}}
          \put(1, 3){\tiny{$7$}}
         \put(6, 7){\tiny{$10$}}
           \put(8, 3){\tiny{$7$}}
        \end{overpic}}}}
        +\f{7}{2}
        {\raisebox{-0.5 cm}{\scalebox{0.9}{\begin{overpic}[unit=1mm, scale=0.5]{twotrivalentvertices.png}
         \put(1, 11){\tiny{$7$}}
          \put(8, 11){\tiny{$7$}}
          \put(1, 3){\tiny{$7$}}
         \put(6, 7){\tiny{$12$}}
           \put(8, 3){\tiny{$7$}}
        \end{overpic}}}}
        +
        {\raisebox{-0.5 cm}{\scalebox{0.9}{\begin{overpic}[unit=1mm, scale=0.5]{twotrivalentvertices.png}
         \put(1, 11){\tiny{$7$}}
          \put(8, 11){\tiny{$7$}}
          \put(1, 3){\tiny{$7$}}
         \put(6, 7){\tiny{$14$}}
           \put(8, 3){\tiny{$7$}}
        \end{overpic}}}}
    \end{equation}
    Replace $f^{(14)}$ with the linear combination of terms from equation \ref{eqn: f7tensorf7expansion}. The terms corresponding to the middle strand in equation \ref{eqn: f7tensorf7expansion} having values of $2,4,$ and $6$ will be zero by Lemma \ref{lem: twoboxesconnectedoneboxbelow}. The term corresponding to the middle strand value of 10 will be zero by Lemma \ref{lem: e7sufficientfourboxestwoontop2bigonbottomiszero}. The term corresponding to $f^{(7)}\otimes f^{(7)}$ then will be the first term of equation \ref{eq: expansionSSontopoff14}. The term corresponding to the middle strand in equation \ref{eqn: f7tensorf7expansion} having value 0 will give the second term of equation \ref{eq: expansionSSontopoff14}, where the factor of 30 comes from $\te{tr}(S^2)$. The terms corresponding to the middle strands being 8 and 12 respectively correspond to the third and fourth terms of equation \ref{eq: expansionSSontopoff14}, as we wished.
\end{proof}

From Remark \ref{rem: aboutjellyfishandevaluation}, it is enough to evaluate the diagram of Lemma \ref{lem: e7sufficientexpandSSontopoff14} in $\mc{P}$. From Lemma \ref{lem: e7sufficienttwoSsonf14onlypossibilities} we know, up to scalar multiples, the only possible diagrams $T$ in (\ref{eqn: twoSsonf14}) that are nonzero. We now use the previous lemma to evaluate all of these diagrams. 

\begin{lem}\label{lem: SSf14SSiszero}
    The jellyfish algorithm evaluates the below diagram to zero.
     \begin{equation}\label{eq: SSf14SSiszero}
        \scalebox{0.9}{\begin{overpic}[unit=1mm, scale=0.6]{fourboxesoneinmiddle.png}
         \put(24,5){$S$}
         \put(6, 5){$S$}
         \put(6, 38){$S$}
         \put(16, 40){\tiny{$1$}}
         \put(24, 38){$S$}
         \put(3,14){\tiny{$7$}}
         \put(28,14){\tiny{$7$}}
         \put(16, 1){\tiny{$1$}}
          \put(3, 30){\tiny{$7$}}
            \put(28, 30){\tiny{$7$}}
            \put(14, 21){$f^{(14)}$}
        \end{overpic}}
    \end{equation}  
\end{lem}

\begin{proof}
    From Lemma \ref{lem: e7sufficientexpandSSontopoff14} we can expand the diagram in (\ref{eq: SSf14SSiszero}) and use that $S$ is U.C.C.S. to see that it is equal to
    \begin{align*}
        {\scalebox{0.8}{%
         \begin{overpic}[unit=1mm, scale=0.6]{fourboxesconnected.png}
         \put(23,4){$S$}
         \put(5, 4){$S$}
         \put(5, 22){$S$}
         \put(16, 25){\tiny{$1$}}
         \put(24, 22){$S$}
         \put(3,14){\tiny{$7$}}
         \put(27,14){\tiny{$7$}}
         \put(16, 1){\tiny{$1$}}
        \end{overpic}
       \hspace{2mm}{\raisebox{1.2cm}{$-\f{30}{8}\cdot \te{tr}(S^2)$}}
        \hspace{2mm}{\raisebox{1.2cm}{$-\f{245}{44}$}}
        {\raisebox{-0.5cm}{
         \begin{overpic}[unit=1mm, scale=0.6]{fourboxesoneinmiddle.png}
         \put(6, 5){$S$}
         \put(6, 38){$S$}
             \put(24,5){$S$}
         \put(24, 38){$S$}
         \put(3,14){\tiny{$4$}}
            \put(3, 30){\tiny{$4$}}
         \put(16, 1){\tiny{$4$}}
          \put(16, 40){\tiny{$4$}}
       \put(28,14){\tiny{$4$}}
            \put(28, 30){\tiny{$4$}}
            \put(13, 21){$f^{(8)}$}
        \end{overpic}}}
        \hspace{2mm}{\raisebox{1.2cm}{$-\f{7}{2}$}}
        {\raisebox{-0.5cm}{
         \begin{overpic}[unit=1mm, scale=0.6]{fourboxesoneinmiddle.png}
         \put(6, 5){$S$}
         \put(6, 38){$S$}
             \put(24,5){$S$}
         \put(24, 38){$S$}
         \put(3,14){\tiny{$6$}}
            \put(3, 30){\tiny{$6$}}
         \put(16, 1){\tiny{$2$}}
          \put(16, 40){\tiny{$2$}}
       \put(28,14){\tiny{$6$}}
            \put(28, 30){\tiny{$6$}}
            \put(13, 21){$f^{(12)}$}
        \end{overpic}}}
             }}
    \end{align*}
    Lemmas \ref{lem: fourboxesconnectedequaling450}, \ref{lem: e7sufficienttraceofSsquared}, \ref{lem: fourboxesoneinmiddle}, and \ref{lem: fourboxesoneinmiddle3750over77} give the calculation of each term which gives that the diagram in (\ref{eq: SSf14SSiszero}) equals $450-\f{900}{8}-\f{245}{44}\cdot 30-\f{7}{2}\cdot \f{3750}{77}=0$.
\end{proof}

\begin{lem}\label{lem: twoboxes1box3boxesiszero}
The jellyfish algorithm evaluates the below diagram to zero. 
    \begin{align*}
        \scalebox{0.9}{\begin{overpic}[unit=1mm, scale=0.6]{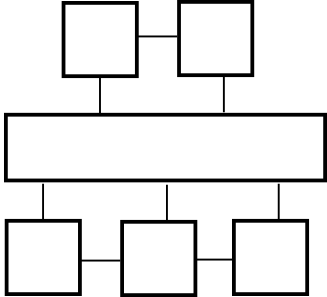}
         \put(6, 5){$S$}
             \put(24,5){$S$}
             \put(42,5){$S$}
         \put(13, 40){$S$}
         \put(33, 40){$S$}
            \put(24, 39){\tiny{$1$}}
            \put(13, 31){\tiny{$7$}}
        \put(32, 31){\tiny{$7$}}
        \put(42, 14){\tiny{$6$}}
        \put(4,14){\tiny{$5$}}
          \put(28,14){\tiny{$3$}}
         \put(16, 2){\tiny{$3$}}
         \put(34,2){\tiny{$2$}}
         \put(23, 23){$f^{(14)}$}
        \end{overpic}}
    \end{align*}
\end{lem}

\begin{proof}
    From Lemma \ref{lem: e7sufficientfourboxestwoontop2bigonbottomiszero} we know the below diagram is zero. 
    \begin{equation}\label{eq: 21rectanglesquare2}
         \scalebox{0.9}{\begin{overpic}[unit=1mm, scale=0.6]{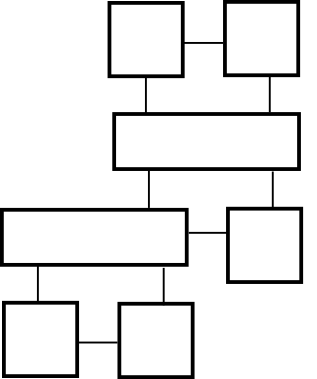}
         \put(5, 5){$S$}
         \put(24,5){$S$}
         \put(22, 53){$S$}
         \put(40, 53){$S$}
         \put(40, 20){$S$}
         \put(3,14){\tiny{$5$}}
         \put(15, 2){\tiny{$3$}}
          \put(20, 45){\tiny{$7$}}
          \put(40,45){\tiny{$7$}}
       \put(28,14){\tiny{$5$}}  
       \put(26, 30){\tiny{$8$}}
       \put(32, 25){\tiny{$2$}}
            \put(40, 30){\tiny{$6$}}
            \put(32, 51){\tiny{$1$}}
            \put(13, 21){$f^{(10)}$}
            \put(30, 37){$f^{(14)}$}
        \end{overpic}}
    \end{equation}
    On the other hand, we can expand the $f^{(10)}$ in the above term. The only terms that will not cup $S$ or $f^{(14)}$ on top would be all straight strands or a cup between strands 8 and 9. Likewise, on bottom, the only terms that would not cap an $S$-box would have a cup in the middle. Therefore, the only terms in the expansion of $f^{(10)}$ that will not cup or cap an $S$-box or $f^{(14)}$ are
    \begin{align*}
         \begin{overpic}[unit=1mm, scale=0.6]{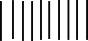}
        \end{overpic}\hspace{1mm},  \hspace{4mm} 
        \begin{overpic}[unit=1mm, scale=0.6]{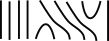}
        \end{overpic}\hspace{1mm},\hspace{4mm} 
          \begin{overpic}[unit=1mm, scale=0.6]{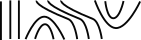}
        \end{overpic}
    \end{align*}
    By Frenkel and Khovanov's result in Proposition \ref{prop: MorrisonJWexpansion}, the first term has a coefficient of $1$ in the expansion of $f^{(10)}$. Let $a,b\in \mbb{Q}$ by the coefficients of the second and third term, respectively. Then expanding the diagram in (\ref{eq: 21rectanglesquare2}) gives
    \begin{equation}\label{eq: SSf14SSSexpansion}
        {\raisebox{1.9cm}{$0=$}}\hspace{1mm}
       {\scalebox{0.8}{%
       \begin{overpic}[unit=1mm, scale=0.6]{twoboxes1box3boxes.png}
         \put(6, 5){$S$}
             \put(24,5){$S$}
             \put(42,5){$S$}
         \put(14, 40){$S$}
         \put(33, 40){$S$}
            \put(24, 39){\tiny{$1$}}
            \put(13, 31){\tiny{$7$}}
        \put(32, 31){\tiny{$7$}}
        \put(42, 14){\tiny{$6$}}
        \put(4,14){\tiny{$5$}}
          \put(28,14){\tiny{$3$}}
         \put(16, 2){\tiny{$3$}}
         \put(34,2){\tiny{$2$}}
         \put(23, 23){$f^{(14)}$}
        \end{overpic}
        }}
        {\raisebox{1.9cm}{$+a$}}\hspace{1mm}
        {\scalebox{0.8}{%
         \begin{overpic}[unit=1mm, scale=0.6]{twoboxes1box3boxes.png}
         \put(6, 5){$S$}
             \put(24,5){$S$}
             \put(42,5){$S$}
         \put(14, 40){$S$}
         \put(33, 40){$S$}
            \put(24, 39){\tiny{$1$}}
            \put(13, 31){\tiny{$7$}}
        \put(32, 31){\tiny{$7$}}
        \put(42, 14){\tiny{$7$}}
        \put(4,14){\tiny{$4$}}
          \put(28,14){\tiny{$3$}}
         \put(16, 2){\tiny{$4$}}
         \put(34,2){\tiny{$1$}}
         \put(23, 23){$f^{(14)}$}
        \end{overpic}
        }}
        {\raisebox{1.9cm}{$+b$}}\hspace{2mm}
        {\scalebox{0.8}{%
        \begin{overpic}[unit=1mm, scale=0.6]{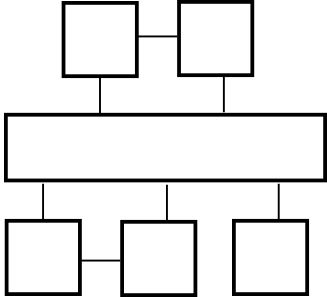}
         \put(6, 5){$S$}
             \put(24,5){$S$}
             \put(42,5){$S$}
         \put(14, 40){$S$}
         \put(33, 40){$S$}
            \put(24, 39){\tiny{$1$}}
            \put(13, 31){\tiny{$7$}}
        \put(32, 31){\tiny{$7$}}
        \put(42, 14){\tiny{$8$}}
        \put(4,14){\tiny{$3$}}
          \put(28,14){\tiny{$3$}}
         \put(16, 2){\tiny{$5$}}
         \put(23, 23){$f^{(14)}$}
        \end{overpic}
        }}
    \end{equation}

    Notice that the first term in the sum is what we want to evaluate for this lemma. We show the other two terms are zero to give the result. For the second term, use the $S^2$ relation to get a sum of two diagrams. The diagram corresponding to $S$ will be zero by Lemma \ref{lem: SSf14SSiszero}. The term corresponding to $f^{(4)}$ will also be zero since there is a single $S$-box on the bottom half which is not one of the cases given in Lemma \ref{lem: e7sufficienttwoSsonf14onlypossibilities}. 

    Now consider the third term in equation \ref{eq: SSf14SSSexpansion}. Use the $S^2$ relation to write as a sum of two diagrams. The diagram corresponding to $S$ will be zero since $S$ is U.C.C.S. The term corresponding to $f^{(4)}$ will also be zero since it will cup $f^{(14)}$. Therefore the second and third term of equation \ref{eq: SSf14SSSexpansion} is zero and we obtain the desired result.
\end{proof}

The next lemma is a symmetric argument to the previous result and will thus be omitted.

\begin{lem}\label{lem: SSf14SSSS2}
The jellyfish algorithm evaluates the below diagram to zero. 
    \begin{align*}
       \scalebox{0.8}{ \begin{overpic}[unit=1mm, scale=0.6]{twoboxes1box3boxes.png}
         \put(6, 5){$S$}
             \put(24,5){$S$}
             \put(42,5){$S$}
         \put(14, 40){$S$}
         \put(33, 40){$S$}
            \put(24, 39){\tiny{$1$}}
            \put(13, 31){\tiny{$7$}}
        \put(32, 31){\tiny{$7$}}
        \put(42, 14){\tiny{$5$}}
        \put(4,14){\tiny{$6$}}
          \put(28,14){\tiny{$3$}}
         \put(16, 2){\tiny{$2$}}
         \put(34,2){\tiny{$3$}}
         \put(23, 23){$f^{(14)}$}
        \end{overpic}}
    \end{align*}
\end{lem}

\begin{lem}\label{lem: SSf14SSSSiszero}
    The jellyfish algorithm evaluates the below diagram to zero. 
    \begin{align*}
        \begin{overpic}[unit=1mm, scale=0.6]{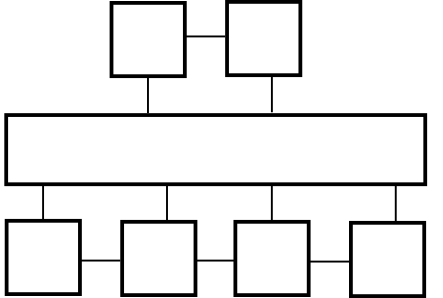}
         \put(6, 5){$S$}
             \put(23,5){$S$}
             \put(43,5){$S$}
         \put(22, 40){$S$}
         \put(40, 40){$S$}
         \put(60, 5){$S$}
            \put(33, 39){\tiny{$1$}}
            \put(20, 31){\tiny{$7$}}
        \put(39, 31){\tiny{$7$}}
        \put(40, 14){\tiny{$2$}}
        \put(4,14){\tiny{$5$}}
          \put(28,14){\tiny{$2$}}
         \put(16, 2){\tiny{$3$}}
         \put(34,2){\tiny{$3$}}
         \put(52, 2){\tiny{$3$}}
         \put(59,14){\tiny{$5$}}
         \put(30, 23){$f^{(14)}$}
        \end{overpic}
    \end{align*}
\end{lem}

\begin{proof}
    We follow a similar process to Lemma \ref{lem: twoboxes1box3boxesiszero}. That is, we use Lemma \ref{lem: e7sufficientfourboxestwoontop2bigonbottomiszero} to know that the below diagram is zero then expand $f^{(10)}$. 

     \begin{equation}\label{eq: 21rectanglesquare4}
        \scalebox{0.85}{ \begin{overpic}[unit=1mm, scale=0.6]{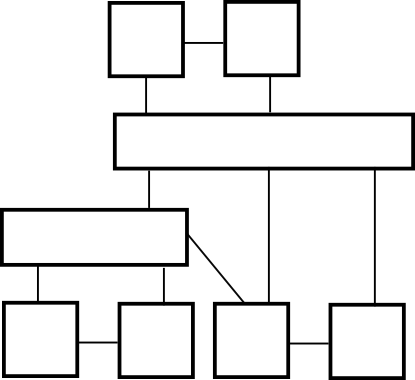}
         \put(5, 5){$S$}
         \put(23,5){$S$}
         \put(22, 52){$S$}
         \put(40, 52){$S$}
         \put(39, 5){$S$}
         \put(57, 5){$S$}
         \put(3,14){\tiny{$5$}}
         \put(15, 2){\tiny{$3$}}
          \put(20, 45){\tiny{$7$}}
          \put(40,45){\tiny{$7$}}
       \put(28,14){\tiny{$5$}}  
       \put(26, 28){\tiny{$7$}}
       \put(33, 22){\tiny{$3$}}
            \put(40, 28){\tiny{$2$}}
            \put(55,28){\tiny{$5$}}
            \put(48, 2){\tiny{$3$}}
            \put(32, 51){\tiny{$1$}}
            \put(12, 22){$f^{(10)}$}
            \put(38, 37){$f^{(14)}$}
        \end{overpic}}
    \end{equation}
    The only terms in the expansion of $f^{(10)}$ that will not cap or cup $f^{(14)}$ or an $S$-box will be the following
 \begin{align*}
        \adjustbox{scale=1.2}{ \begin{overpic}[unit=1mm, scale=0.6]{10straightstrands.png}
        \end{overpic}}\hspace{1mm},  \hspace{4mm} 
        \adjustbox{scale=1.2}{ \begin{overpic}[unit=1mm, scale=0.6]{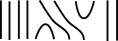}
        \end{overpic}}\hspace{1mm},\hspace{4mm} 
          \adjustbox{scale=1.2}{ \begin{overpic}[unit=1mm, scale=0.6]{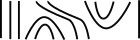}
        \end{overpic}}
        \hspace{1mm},\hspace{4mm} 
         \adjustbox{scale=1.2}{  \begin{overpic}[unit=1mm, scale=0.6]{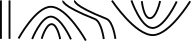}
        \end{overpic}}
    \end{align*}
    which will have coefficients 1, $a,b,c\in \mbb{Q}$, respectively. The term with coefficient 1 is the diagram we are trying to compute for the lemma. So it is sufficient to prove that the diagrams corresponding to the other coefficients are zero. These are the respective diagrams:
     \begin{align*}
        \scalebox{0.7}{%
        \begin{overpic}[unit=1mm, scale=0.6]{214.png}
         \put(6, 5){$S$}
             \put(24,5){$S$}
             \put(42,5){$S$}
         \put(22, 40){$S$}
         \put(41, 40){$S$}
         \put(60, 5){$S$}
            \put(33, 39){\tiny{$1$}}
            \put(20, 31){\tiny{$7$}}
        \put(39, 31){\tiny{$7$}}
        \put(40, 14){\tiny{$3$}}
        \put(4,14){\tiny{$4$}}
          \put(28,14){\tiny{$2$}}
         \put(16, 2){\tiny{$4$}}
         \put(34,2){\tiny{$2$}}
         \put(52, 2){\tiny{$3$}}
         \put(59,14){\tiny{$5$}}
         \put(29, 23){$f^{(14)}$}
        \end{overpic}}
        \hspace{1mm},
          \scalebox{0.7}{%
         \begin{overpic}[unit=1mm, scale=0.6]{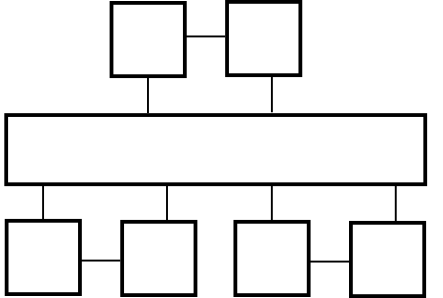}
         \put(6, 5){$S$}
             \put(24,5){$S$}
             \put(42,5){$S$}
         \put(22, 40){$S$}
         \put(41, 40){$S$}
         \put(60, 5){$S$}
            \put(33, 39){\tiny{$1$}}
            \put(20, 31){\tiny{$7$}}
        \put(39, 31){\tiny{$7$}}
        \put(40, 14){\tiny{$4$}}
        \put(4,14){\tiny{$3$}}
          \put(28,14){\tiny{$2$}}
         \put(16, 2){\tiny{$5$}}
         \put(34,2){\tiny{$1$}}
         \put(52, 2){\tiny{$3$}}
         \put(59,14){\tiny{$5$}}
         \put(29, 23){$f^{(14)}$}
        \end{overpic}}\hspace{1mm},
          \scalebox{0.7}{%
         \begin{overpic}[unit=1mm, scale=0.6]{214nostrand.png}
         \put(6, 5){$S$}
             \put(24,5){$S$}
             \put(42,5){$S$}
         \put(22, 40){$S$}
         \put(41, 40){$S$}
         \put(60, 5){$S$}
            \put(20, 31){\tiny{$7$}}
        \put(39, 31){\tiny{$7$}}
        \put(40, 14){\tiny{$5$}}
        \put(4,14){\tiny{$3$}}
          \put(28,14){\tiny{$2$}}
         \put(16, 2){\tiny{$6$}}
         \put(34,2){\tiny{$1$}}
         \put(52, 2){\tiny{$3$}}
         \put(59,14){\tiny{$5$}}
         \put(29, 23){$f^{(14)}$}
        \end{overpic}}
    \end{align*}
    All of these diagrams involve an $S^2$. Consider the first diagram above. When expanding the $S^2$, the $S$ term will become zero by Lemma \ref{lem: SSf14SSSS2}. The $f^{(4)}$ term will cup $f^{(14)}$. The second and third diagrams above also have an $S^2$ term. The resulting $S$ term will be sidecapped, resulting in zero. The resulting $f^{(4)}$ term will cup $f^{(14)}$. Thus all these diagrams are zero.
\end{proof}

\begin{lem}
    The jellyfish algorithm evaluates the two diagrams in (\ref{eq: twoequivalentdiagramsuntwistwithf}) equivalently for all $c\in \{0,...,8\}$. 
\end{lem}

\begin{proof}
    Suppose that the jellyfish algorithm has been completed up through step 5 for two diagrams which at this point of the jellyfish algorithm look like the diagrams from (\ref{eq: twoequivalentdiagramsuntwistwithf}). Recall that the proof of Lemma \ref{lem: canuntwisttwoSboxes} only relied on resolving crossing. So by Lemma \ref{lem: canuntwisttwoSboxes}, when $c$ is even this lemma is true. 

    When $c$ is odd, Lemma \ref{lem: canuntwisttwoSboxes} gives that the diagrams in (\ref{eq: twoequivalentdiagramsuntwistwithf}) evaluate to be the negatives of each other in the jellyfish algorithm. We show that the diagrams evaluate to zero. Continuing onto step 7 of the jellyfish algorihtm, we can assume that all $S^2$ boxes along the top and along the bottom have been reduced. Therefore, only $c=1$ and $c=3$ cases need to be considered. However, when $c=3$ and $c=1$, Lemmas \ref{lem: e7sufficientfourboxestwoontop2bigonbottomiszero}, \ref{lem: e7sufficienttwoSsonf14onlypossibilities}, \ref{lem: SSf14SSiszero}, \ref{lem: twoboxes1box3boxesiszero}, \ref{lem: SSf14SSSS2}, \ref{lem: SSf14SSSSiszero}, give that the jellyfish algorithm evaluates to zero, giving the result.
\end{proof}

The same proof can be used for the below lemma.

\begin{lem}
    The jellyfish algorithm evaluates the two diagrams in (\ref{eqn: twistotherway}) equivalently for all $c\in \{0,...,8\}$. 
\end{lem}

Using Lemmas \ref{lem: twistanduntwistthesame1} and \ref{lem: twistanduntwistthesame2}, now the following two lemmas are immediate.

\begin{lem}\label{lem: jfatwistsareequal1}
    The jellyfish algorithm evaluates the diagrams in (\ref{eq: twoequivalentdiagramsuntwist}) to be the same number.
\end{lem}

\begin{lem}\label{lem: jfatwistsareequal2}
   The jellyfish algorithm evaluates the diagrams in (\ref{eqn: twistotherway2}) to be the same number.   
\end{lem}

\begin{lem}\label{lem: jfarespectsorderofsboxes}
    The jellyfish algorithm is invariant under the order of the $S$-boxes chosen to drag to the top.
\end{lem}

\begin{proof}
    The symmetric group is generated by transpositions of the form $(i \te{ }i+1)$, so it is enough to prove that the jellyfish algorithm is invariant under switching the $S$-box getting dragged to position $i$ with the $S$-box getting dragged to position $i+1$. 

    Consider a closed diagram $D \in \mc{P}_0$ with $m$ $S$-boxes. We want to compute the jellyfish algorithm in two ways. The first way drags an $S$-box, say $S_i$ to point $i$ at the top and an $S$-box, $S_{i+1}$ to point $i+1$. The other way drags the $S_i$ box to point $i+1$ and $S_{i+1}$ to point $i$. All other steps of the jellyfish algorithm up to this point are identical. 

    It has already been proven that the choice of arc doesn't matter. There exists some neighborhood of the point $i$ and $i+1$ at the top where there are no strands. Call this neighborhood $N$. So in the first case, choose any arc, $\gamma_i$ to drag from $S_i$ to point $i$. For $S_{i+1}$, first drag $S_{i+1}$ to the position where $S_i$'s star originally was, then follow the same path as $\gamma_i$ up to the top until the $S_{i+1}$ box is inside the neighborhood $N$, then complete $\gamma_{i+1}$ by making an arc passing no more strands to the point $i+1$. (Technically, $\gamma_{i+1}$ should always be on the right, within some local neighborhood of $\gamma_i$ that crosses the same exactly strands as $\gamma_i$, since strands must cross transversely). 

    In the second case, make an imaginary arc $\gamma_i'$ from the star of $S_{i+1}$ that follows $\gamma_{i+1}$ until the neighborhood $N$, then follow the arc $\gamma_i$ to the point $i$. Drag $S_{i+1}$ to point $i$ using this arc. Then make arc $\gamma_{i+1}'$ from the star of $S_i$ following $\gamma_i$ until the neighborhood $N$, then follow $\gamma_{i+1}$ to the point $i+1$. 

    Showing that the first and second case are equivalent are is thus equivalent to showing that for all $c$ where $0\leq c \leq 8$, the jellyfish algorithm evaluates the two diagrams in (\ref{eq: twoequivalentdiagramsuntwist}) to the same number (or the two diagrams in (\ref{eqn: twistotherway2})). However, this was done in Lemmas \ref{lem: jfatwistsareequal1} and \ref{lem: jfatwistsareequal2}, so this proof is complete.
\end{proof}

\subsection{The Jellyfish Algorithm: Respecting the Defining Relations}\label{jellyfish4}

The last step in showing the jellyfish algorithm is a well-defined function is to show it respects the defining relations. 

\begin{lem}\label{lem: jfarespectsbubble}
    The jellyfish algorithm is invariant under the bubble relation.
\end{lem}

\begin{proof}
    Let $D\in \mc{P}_0$ be some closed diagram with $m\in \mbb{Z}_{\geq 0}$ $S$-boxes and somewhere a bubble. Let $D'$ be the identical to $D$ except without the bubble and instead multiplied by 2. As the jellyfish algorithm is invariant under the choice of arc, do steps 1 through 7(iii) identically and for the imaginary arcs, do not go through the bubble. In step 7(iv) pop the bubble in $D$ for a factor of 2. Now $D$ and $D'$ are identical and will give the same number in the jellyfish algorithm.
\end{proof}

\begin{lem}\label{lem: jfarespectsuccs}
    The jellyfish algorithm is invariant with respect to the $S$ is U.C.C.S. relation. 
\end{lem}

\begin{proof}
    Suppose $D\in \mc{P}_0$ is a closed diagram with some $S$-box that is capped, cupped, or sidecapped on the right. Suppose that this box is dragged to position $i$ at the top the diagram during the jellyfish algorithm. Then, since the jellyfish algorithm stretches strands up to the point, we get that at the end of step 5, we will have the $S_i$ at the top with a cap on it that is still in the same position. Use Reidemeister moves to move any strands intersecting the cap. Then at step 7(i) this diagram will evaluate to zero. 

    If instead $S$ is sidecapped on the left, then we are in the case of the below
    \begin{align*}
         \begin{overpic}[unit=1mm, scale=0.4]{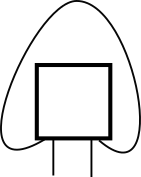}
        \put(6.5, 6.5){$S$}
        \put (6.3, 2){$...$}
    \end{overpic}.
    \end{align*}
    Call the open neighborhood ball that looks like the above picture $N$. In this case, since the order in dragging $S$ boxes to the top doesn't matter, choose this box to be $S_1$. Choose the imaginary arc $\gamma_1$ to 1 that is a straight line to the boundary of $N$, then goes to the point 1 by never re-entering the neighborhood. Then, when $S_1$ reaches the top it will locally look like the below.
    \begin{equation}\label{eqn: boxtwistsonbottom}
         \begin{overpic}[unit=1mm, scale=0.4]{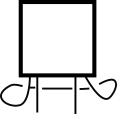}
        \put(4.7, 6.5){$S$}
        \put (4.5, 1){$...$}
    \end{overpic}.
    \end{equation}
    Finish step 5 by dragging the rest of the $S$ boxes up, avoiding the neighborhood that looks like the above. In step 6 resolve the crossings. Using Lemmas \ref{lem: cuponcrossingresolved} and \ref{lem: boxwithovercrossingonbottom} we get that the above local neighborhood will resolve the crossing to have a cap on an $S$. Thus $D$ will get sent to zero by the jellyfish algorithm.
\end{proof}

\begin{lem}\label{lem: jfarespectsssquared}
    The jellyfish algorithm is invariant under the $S^2$ relation.
\end{lem}

\begin{proof}
    Let $D\in \mc{P}_0$ be a closed diagram with two $S$-boxes forming an $S^2$. Let $D'\in \mc{P}_0$ be a sum of two closed diagrams where they are identical to $D$ except that in a local neighborhood, $B$, the $S^2$ has been replaced by an $S$-box in one and 6 times $f^{(4)}$ in the other. Make a point on the boundary of $B$. Suppose that $D$ has $m-1$ boxes and $D'$ has $m$ boxes. There is a neighborhood, $N$, in $\mc{D}$ containing points $m$ and $m-1$ at the top of all the diagrams that contains no strands. As the order of the $S$-boxes doesn't matter, we can choose to drag the the two $S$-boxes in $\mc{D}$ to the last two positions, $m-1$ and $m$. When dragging the first $m-2$ $S$-boxes in $\mc{D}$ and $\mc{D}'$, we choose identical imaginary arcs that avoid the open ball $B$, the distinguished point on the boundary of $B$, and the neighborhood $N$. 
    
    As the jellyfish algorithm is invariant under the choice of arc, we can first choose the arc $\gamma_{m-1}$ that is a straight line inside of the open ball until it hits the designated point. Then choose any arc from there to the point $m-1$ that never re-enters the open ball. For the arc $\gamma_m$, make a straight line from the star of the second $S$-box to the designated point on the boundary. Then follow the path of $\gamma_{m-1}$ until it reaches the neighborhood $N$, then go to the point $m$ by staying inside the neighborhood $N$. 

    For the summand of $D'$ with the single $S$-box in $B$, choose this box to go the $m-1$ position. In this case, this is the last point on the top of the diagram. For the arc $\gamma_{m-1}$, first choose an arc that goes from the star of $S$ to the designated point of the ball $B$, staying within the ball, and then follow the arc of $\gamma_{m-1}$ to the point $m-1$. Now the jellyfish algorithm can be assumed to be done for $D'$ up through step 5. In step 6, we can perform planar isotopies for the parts of the diagram in Temperley-Lieb. So for the $D'$ summand containing $6f^{(4)}$, first isotope the $f^{(4)}$ using the same arc as the other summand for $D'$ to the top of the diagram. 

    For $D$, we now have a diagram of the below form where $c\geq 4$ and $T\in \te{Hom}(\emptyset, X^{\otimes (16-2c})$ (and $T$ may include crossings):
    \begin{align*}
         \begin{overpic}[unit=1mm, scale=0.6]{twoboxesconnectedoneboxbelowclosed.png}
         \put(28,12){\tiny{$$8-c}}
         \put(0, 11){\tiny{$8-c$}}
         \put(5, 20){$S$}
         \put(15, 19){\tiny{$c$}}
         \put(24, 20){$S$}
         \put(15, 3){$T$}
        \end{overpic}\hspace{1mm}.
    \end{align*}
    On the other hand, the jellyfish form of $D'$ is 
    \begin{align*}
         \begin{overpic}[unit=1mm, scale=0.6]{twoboxesstacked.png}
         \put(13, 20){$S$}
         \put(17, 11){\tiny{$8$}}
         \put(13, 4){$T$}
        \end{overpic}
        \hspace{2mm}{\raisebox{1cm}{$+6$}}\hspace{2mm}
             \begin{overpic}[unit=1mm, scale=0.6]{twoboxesstacked.png}
          \put(12, 20){$f^{(4)}$}
         \put(17, 11){\tiny{$8$}}
         \put(13, 4){$T$}
        \end{overpic} 
    \end{align*}
    Complete steps 6 through 7(i) identically for $D$ and the summands of $D'$. At step 7(ii) we can reduce the $S^2$ in any order we choose. In $D$, first reduce the $S^2$ at points $m$ and $m-1$, which then will give the identical sum of diagrams in $D'$. Therefore, the jellyfish algorithm will evaluate $D$ and $D'$ to be the same number.
\end{proof}

\begin{lem}\label{lem: jfarespectsjfrelation}
    The jellyfish algorithm is invariant under the jellyfish relation. 
\end{lem}

\begin{proof}
    Let $D$ and $D'$ be closed diagrams in $\mc{P}_0$ except in some neighborhood, $N$, $D$ has the left side of the jellyfish relation (an $S$ box underneath a strand) and $D'$ has the right side of the jellyfish relation (an $S$-box with a strand crossed over its strands beneath it). Pick a distinguished point on the boundary of $N$. 

    Begin the jellyfish algorithm identically for $D$ and $D'$, avoiding having imaginary arcs going through the neighborhood $N$. Now consider the $S$-box in the neighborhood $N$, for $D$ choose a straight line from the star of the $S$ box, through the overstand to the distinguished point of $N$, then a path outside of $N$ to the point on the boundary that does not re-enter $N$. Choose the same imaginary arc for $D'$. Notice now, when dragging the $S$-box in diagram $D$, once it passes under the overstrand inside of $N$, it looks identical to the diagram $D'$. Now the diagrams are equal and the jellyfish algorithm will evaluate them to the same number.
\end{proof}

\begin{proof}[Proof of Theorem \ref{thm: e7jfaiswelldefined}]
   We showed the choice in reducing $S^2$ doesn't matter in Lemma \ref{lem: e7sufficientwelldefinedSsquared}. We showed the choice of arc doesn't matter in Lemma \ref{lem: jkarespectsarc}.   We showed that the order of choosing the $S$-boxes to float to the top doesn't matter in Lemma \ref{lem: jfarespectsorderofsboxes}. Lastly, we showed that the function respects the defining relations in Lemmas \ref{lem: jfarespectsbubble}, \ref{lem: jfarespectsuccs}, \ref{lem: jfarespectsssquared}, and \ref{lem: jfarespectsjfrelation}. 
\end{proof}

\begin{cor}\label{cor: e7sufficientp0onedim}
    $\mc{P}_0$ is one-dimensional. 
\end{cor}

\begin{proof}
    We have a well-defined map $J:\mc{P}_0\to \mbb{C}$ that sends the empty diagram to 1. So the dimension of $\mc{P}_0$ is greater than or equal to 1. On the other hand, this map is an evaluation algorithm, so the dimension of $\mc{P}_0$ is at most 1. 
\end{proof}

\begin{rem*}
    This has now shown that the jellyfish algorithm is an invariant on the vector space of closed diagrams for the affine $E_7$ subfactor planar algebra.
\end{rem*}

\subsection{Minimal Projections in $\mc{P}$}

In this section, we prove that the projections $\emptyset$, $P_1=f^{(1)}$, $P_2=f^{(2)}$, $P_3=f^{(3)}$, $P_4=\f{3}{5}f^{(4)}-\f{1}{5}S$, $P_5=P_4\otimes X-\f{4}{3}(P_4\otimes X)e_4(P_4 \otimes X)$, $P_6=P_5 \otimes X-\f{3}{2}(P_5 \otimes X)e_5(P_5 \otimes X)$, and $Q_4=\f{2}{5}f^{(4)}+\f{1}{5}S$ are nonisomorphic minimal projections and satisfy the tensor rules for the affine $E_7$ principal graph.

We define a modified jellyfish algorithm for $\mc{P}$ that can be used for any box space. The goal of this algorithm is to put the diagrams in a jellyfish form with $S$-boxes connected by at most 2 strands. 

 \begin{tcolorbox}[breakable, pad at break*=0mm]
 Define the \tbf{modified jellyfish algorithm} on diagrams in the following way. Let $D\in \mc{P}_k$ be a diagram.
 \begin{enumerate}
     \item Follow the steps 1 through 6 of the jellyfish algorithm for $\mc{P}$ except dragging $S$-boxes to the left side of the diagram instead of to the top and keeping the top and bottom boundary points fixed. This may make $D$ a sum of diagrams, say $D=\sum_{i=1}^{n_\ell}a_i D_i$, where $a_i \in \mbb{C}$.
     \item Look individually at each $D_i$.
     \begin{enumerate}
        \item If $D_i$ has any $S$-box that is cupped, evaluate $D_i$ to zero. 
         \item  If $D_i$ has no $S$-box that is cupped, then use the $S^2$ relation to reduce any diagram with an $S^2$ to a sum of diagrams with at least one less $S$ box. Repeat this process until all the resulting summands have adjacent $S$-boxes are connected by at most 3 strands. 
         \item Use equation \ref{eqn: writeS3SaslinearcombooflessS} to write any diagram with two $S$-boxes connected by 3 strands to be a summand of diagrams with at least one less $S$-box. Repeat this process until all the resulting summands have adjacent $S$-boxes connected by at most 2 strands.
     \end{enumerate}
     \item Extend the function linearly to be defined on all of $\mc{P}_k$. 
 \end{enumerate}
 \end{tcolorbox}

A diagram where all of its $S$-boxes are on the left side with no strands crossing and each box is connected by at most 2 strands is said to be in \tbf{modified jellyfish form}. 

We first show the below result which will be used frequently in the lemmas that follow. We thank the referee for guidance in the proof of this lemma.

\begin{lem}\label{lem: there-is-an-S2}
    Let $D$ be a diagram in $\mc{P}$ where the first step of the modified jellyfish algorithm has been taken. Label the $S$-boxes from top to bottom $S_1,...,S_m$. Suppose there is no $S$-box with a cup. If box $S_i$ connects to box $S_j$, where $i+1<j$, then between $S_i$ and $S_j$ there is a pair of $S$-boxes forming an $S^2$. 
\end{lem}

\begin{proof}
      Think of $S_i,...,S_j$ as vertices of a graph. Label the vertices by the label of their associated $S$-box. Connect vertices by edges that match how the $S$-boxes were connected. However, instead of $k$ edges between two vertices, just have 1 edge and label it $(k)$. Now we have formed a $(j-i+1)$-gon with some interior edges. 

      Every $n$-gon is triangulated by $n-3$ interior edges. If the graph we created is not yet triangulated, add interior edges colored blue until it is. By the Two Ears Theorem \cite{Mei75}, when $n\geq 3$, every triangulated $n$-gon has two nonadjacent ``ears", i.e., a vertex $a$ of a triangle $abc$ where edge $bc$ is interior to the $n$-gon. Since vertex $a$ is on the boundary of our graph, there cannot be any edge of $a$ connecting to any vertices other than $b$ or $c$ without introducing a crossing. Further, these edges connecting $a$ to $b$ and $c$ are not interior edges, so they were not colored blue. Thus the edges belonged to the original diagram.

      Since we have two nonadjacent ears, the ears cannot be both vertex $i$ and $j$. Choose an ear that is not $i$ or $j$. Then this vertex, and thus the corresponding $S$-box, can only connect to its neighbors. By the pigeonhole principle, this $S$-box must be connected by at least 4 strands with one of the neighbors, thus forming an $S^2$.
\end{proof}

 Now we begin prove the main results of this subsection. Notice that Corollary \ref{cor: e7sufficientp0onedim} showed that $\emptyset$ is a minimal projection. 

 \begin{lem}
     For $n\in \{1,2,3\}$, $\te{Hom}(P_n, P_n)$ is one-dimensional. 
 \end{lem}

 \begin{proof}
     Fix $n \in \{1,2,3\}$. Let $f \in \te{Hom}(P_n, P_n)$ and use the modified jellyfish algorithm on the middle of $f$. Without loss of generality, suppose that the middle of $f$ is a diagram in modified jellyfish form. Let $m$ be the number of $S$-boxes in $f$ in this form. Let $S_1$ be the $S$-box at the top left. $S_1$ can have at most 3 strands on the top boundary, depending on the value of $n$. If $m\geq 2$, then $S_1$ can have at most 2 strands connected to the $S$-box below it, say $S_2$. It then has at least 3 more strands to connect to either the bottom or another $S$-box. 

     Say a strand of $S_1$ connects to another $S$-box, say $S_k$. Then between $S_1$ and $S_k$ there must be some $S^2$ by Lemma \ref{lem: there-is-an-S2}. This contradicts that each of the $S$-boxes are connected by at most 2 strands. A similar argument shows that the strands of $S_1$ cannot connect to the bottom boundary. Therefore, $m<2$.

     Suppose $m=1$. As there are eight strands of $S$ and six points on the boundary, $S$ must have a cap, so $f=0$. Suppose $m=0$. Then as $f \in \mc{TL}$ and $P_1=f^{(1)}$, $P_2=f^{(2)}$, and $P_3=f^{(3)}$, we know from Temperley-Lieb being a subfactor planar algebra that $f$ must be a scalar multiple of the projection. Therefore, the dimension of $\te{Hom}(P_n,P_n)$ is at most 1. On the other hand, $\te{id}_{P_n}=P_n\in \te{Hom}(P_n,P_n)$. The trace of $P_n$ is nonzero and the 0-box space is one-dimensional, so $P_n\neq 0$. Therefore, the $\te{Hom}(P_n,P_n)$ is one-dimensional.
 \end{proof}

The below calculations will become useful later on.
\begin{lem}\label{lem: e7sufficientprojectiontimesS}
    The following equalities hold in $\mc{P}$:
    \begin{enumerate}[label=(\roman*)]
        \item $P_4S=SP_4=-2P_4$,
        \item $Q_4S=SQ_4=3Q_4$,
        \item $P_5(S\otimes X)=(S\otimes X)P_5=-2P_5$, and 
        \item $P_6(S\otimes X^{\otimes 2})=(S\otimes X^{\otimes 2})=-2P_6$.
    \end{enumerate}
\end{lem}

\begin{proof}
    For part (i) we use the $S^2$ relation as well as that $S$ absorbs $f^{(4)}$ to see that 
    \begin{align*}
        P_4S&=(\f{3}{5}f^{(4)}S-\f{1}{5}S^2)= \f{3}{5}S-\f{1}{5}(S+6f^{(4)})= \f{2}{5}S-\f{6}{5}f^{(4)}=-2(\f{3}{5}f^{(4)}-\f{1}{5}S)=-2P_4.
    \end{align*}
    By symmetry, we also have $SP_4=-2P_4$. 

    For part (ii) we use similar reasoning to get that 
    \begin{align*}
        Q_4S=SQ_4=(\f{2}{5}S+\f{1}{5}S^2)=\f{2}{5}S+\f{1}{5}(S+6f^{(4)})= \f{3}{5}S+\f{6}{5}f^{(4)}=3(\f{1}{5}S+\f{2}{5}f^{(4)})=3Q_4,
    \end{align*}
    as desired.

    For part (iii), we use part (i):
    \begin{align*}
        P_5(S\otimes X)&=(P_4 \otimes X)(S\otimes X)-\f{4}{3}(P_4 \otimes X)e_4(P_4 \otimes X)(S\otimes X)=-2(P_4 \otimes X)+\f{4}{3}\cdot 2(P_4\otimes X)e_4(P_4 \otimes X)
        = -2P_5.
    \end{align*}
    Similarly, $(S\otimes X)P_5=-2P_5$. For part (iv), we use part (iii):
    \begin{align*}
        P_6(S\otimes X^{\otimes 2})=(P_5 \otimes X)(S\otimes X^{\otimes 2})-\f{3}{2}(P_5 \otimes X)e_5 (P_5 \otimes X)(S\otimes X^{\otimes 2})
        =-2(P_5 \otimes X)+\f{3}{2}\cdot 2(P_5 \otimes X)e_5(P_5 \otimes X)
        =-2P_6.
    \end{align*}
    Similarly, $(S\otimes X^{\otimes 2})P_6=-2P_6$. 
\end{proof}

 \begin{lem}
     Both $\te{Hom}(P_4, P_4)$ and $\te{Hom}(Q_4, Q_4)$ are one-dimensional. 
 \end{lem}

 \begin{proof}
     Let $f$ be in either $\te{Hom}(P_4, P_4)$ or $\te{Hom}(Q_4, Q_4)$ and without loss of generality, suppose $f$ is a diagram whose middle is in modified jellyfish form. Label the $S$-boxes on the left, starting with the top: $S_1,...S_m$. 

     Suppose $S_1$ connects to the top boundary by $1\leq p \leq 4$ strands and that $m\geq 2$. Since $f$ is in modified jellyfish form, $S_1$ connects to $S_2$ by at most 2 strands. Then $S_1$ has at least two other strands to connect somewhere. However, by Lemma \ref{lem: there-is-an-S2} this will require somewhere to have an $S^2$, which cannot happen. Thus $m=1$ or $0$. 

     Suppose $m=1$. If $f$ is not zero, $S$ has no cap so has all of the strands are on the boundary. Thus the rest of $f$ is in Temperley-Lieb which can be evaluated to a scalar, say $a \in \mbb{C}$. So the middle of $f$ is $aS$. That is, $f=aP_4SP_4$, or $f=aQ_4SQ_4$. Thus, it is enough to calculate $P_4SP_4$ and $Q_4SQ_4$ in this case. Using Lemma \ref{lem: e7sufficientprojectiontimesS}, we get
     \begin{align*}
         P_4SP_4=-2P_4^2=-2P_4.
     \end{align*}
     Similarly, we also get
     \begin{equation*}
         Q_4SQ_4=3Q_4^2=3Q_4.
     \end{equation*}
     So in either case $f$ is a multiple of the projection. 
     
     If $m=0$ then the middle of $f$ is in Temperley-Lieb. $P_4$ and $Q_4$ are uncuppable and uncappable, so the only nonzero option would be the middle of $f$ to be straight strands. Thus $f$ is a scalar multiple of the projection. Thus the hom spaces are at most one-dimensional. The traces of $P_4$ and $Q_4$ are both nonzero, so the hom spaces are one-dimensional. 
 \end{proof}

 \begin{lem}\label{lem: e7sufficienthomp5p5}
     The dimension of $\te{Hom}(P_5, P_5)$ is 1.
 \end{lem}

 \begin{proof}
     Without loss of generality, let $f$ be a diagram whose middle is in modified jellyfish form. Enumerate the $S$-boxes on the left side: $S_1$, ...,$S_m$. 

     Suppose $m\geq 2$. Then $S_1$ connects to the top by at most 5 strands and to $S_2$ by at most 2 strands, leaving at least 1 strand that will need to pass by $S_2$. If the strand connects to another $S$-box, say $S_k$, then between $S_1$ and $S_k$ there will be an $S^2$ by Lemma \ref{lem: there-is-an-S2}, which would contradict $f$ being in modified jellyfish form. 

     Assume $m=1$. Then, up to some scalar multiple, we must have $f=P_5(S\otimes X)P_5=-2P_5$, by Lemma \ref{lem: e7sufficientprojectiontimesS}. 

     Assume $m=0$. Then $f$ is in Temperley-Lieb. Since $P_5$ is uncuppable and uncappable, the only option is for the middle of $f$ to be a scalar multiple of straight strands. Thus, $f=aP_5^2=P_5$ for some $a\in \mbb{C}$. Since the trace of $P_5$ is nonzero, the dimension of $\te{Hom}(P_5, P_5)=1$. 
 \end{proof}

 \begin{lem}
    The dimension of $\te{Hom}(P_6, P_6)$ is $1$. 
 \end{lem}

 \begin{proof}
     Without loss of generality, let $f$ be a diagram whose middle is in modified jellyfish form. Enumerate the $S$-boxes from top to bottom: $S_1$, ... ,$S_m$. 

     Suppose $m\geq 2$. $S_1$ connects to the top by at most $6$ strands. Suppose $S_1$ connects to the top by at most 5 strands. Then the same argument as in Lemma \ref{lem: e7sufficienthomp5p5} gives a contradiction. Suppose $S_1$ connects to the top by 6 strands. Then since $P_6$ absorbs $S$, $P_6$ will be capped, resulting in $f=0$. 

    Let $m=1$. Then, up to scalar multiple, $f=P_6(S\otimes X^{\otimes 2})P_6=-2P_6^2=-2P_6$. 

    Finally, suppose $m=0$. Then, as $P_6$ is uncappable and uncuppable, the middle of $f$ is a scalar multiple straight strands, so $f=aP_6$ for some $a\in \mbb{C}$. Therefore, $\te{Hom}(P_6, P_6)$ is at most one-dimensional. Since the trace of $P_6$ is nonzero, the hom space is one-dimensional.  
 \end{proof}

 Thus we have proved that $P_0=\emptyset$, $P_1$, ... , $P_6$, and $Q_4$ are all minimal projections. We now show that they are all nonisomorphic. 

 \begin{lem}
     The dimension of $\te{Hom}(P_i, P_j)=0$ for all $i\neq j$, $i,j \in \{0, ...,6\}$ (and $P_0\coloneq \emptyset$).
 \end{lem}

 \begin{proof}
     Fix $0\leq i < j \leq 6$. Without loss of generality, suppose $f\in \te{Hom}(P_i, P_j)$ is a diagram whose middle is in modified jellyfish form. Enumerate the $S$-boxes starting with the top as $S_1$, ..., $S_m$. 

     Suppose $m\geq 2$. As $j\leq 6$, $S_1$ is connected by at most 6 strands at the top. Since $P_6$ and $P_5$ absorbs $S$, if $S_1$ is connected by more than 4 strands then $f=0$. If $S_1$ is connected by 4 or less strands at the top then it is connected by at most 2 strands to $S_2$, leaving at least 2 strands that must go to the bottom boundary or another $S$-box. As seen in Lemma \ref{lem: there-is-an-S2}, this will force there to be some $S^2$ which contradicts the middle of $f$ being in modified jellyfish form. 

     Suppose $m=1$. Then $i+j \geq 8$ since all the strands of $S$ must be on the boundary. Thus, we have $P_j$ is either $P_5$ or $P_6$. Both $P_5$ and $P_6$ absorb $S$, so the middle of $f$ can then be recognized to be in Temperley-Lieb. We then consider the final case.

     Supppose $m=0$. Then the middle of $f$ is in Temperley-Lieb. Since the number of boundary points on the top and bottom are not equal, there must be a cap or cup. However, all $P_i, P_j$ are uncuppable and uncappable, so $f=0$. 
 \end{proof}

\begin{lem}
    For all $0\leq i \leq 6$, the dimension of $\te{Hom}(Q_4, P_i)=0$. 
\end{lem}

\begin{proof}
    By substituting $Q_4$ for $P_4$ in the previous lemma, we get the result for $i\neq 4$. Let $f \in \te{Hom}(Q_4, P_4)$ be a diagram whose middle is in modified jellyfish form. Since there are exactly 4 points on the top and bottom, using similar reasoning to the previous lemmas, the only option is that $f$ is a scalar multiple of $Q_4SP_4=3Q_4P_4$ (when the middle of $f$ has one $S$-box) or $Q_4P_4$ (when the middle of $f$ is in Temperley-Lieb). However,
    \begin{align*}
        Q_4P_4&=(\f{3}{5}f^{(4)}-\f{1}{5}S)(\f{2}{5}f^{(4)}+\f{1}{5}S)
        =\f{6}{25}f^{(4)}+\f{3}{25}S-\f{2}{25}S-\f{1}{25}S^2
        =\f{6}{25}f^{(4)}+\f{1}{25}S-\f{1}{25}(S+6f^{(4)})
        =0.
    \end{align*}
\end{proof}

In this section we have now proved the following.

\begin{cor}
    The projections $\emptyset$, $P_1$, $P_2$, $P_3$, $P_4$, $P_5$, $P_6$, and $Q_4$ are all minimal and nonisomorphic.
\end{cor}

Next we show that the minimal projections of $\mc{P}$ found in the previous corollary satisfy the properties to have the affine $E_7$ Dynkin diagram as its principal graph.

\begin{lem}\label{lem: e7sufficientprincipalgraphissatisfied}
    In $\mc{P}$, the minimal projections satisfy the principal graph given in (\ref{eq: generalE7graph}). 
\end{lem}

\begin{proof}
    Notice that $\emptyset \otimes X \cong P_1$, $P_1 \otimes X \cong \emptyset \oplus P_2$, $P_2 \otimes X \cong P_1 \oplus P_3$, $P_3 \otimes X \cong P_2 \oplus P_4 \oplus Q_4$, $P_4 \otimes X \cong P_3 \oplus P_5$, and $P_5 \otimes X \cong P_4 \oplus P_6$ were all shown in Lemma \ref{lem: e7sufficienttensorproductsofP}. What is left to show is that $Q_4 \otimes X \cong P_3$ and $P_6 \otimes X \cong P_5$. 

    In Lemma \ref{lem: e7sufficientleafrelationq4}, we found that $Q_4 \otimes X = 2 (Q_4 \otimes X)e_4(Q_4 \otimes X)$. Define $f=\sqrt{2}(Q_4 \otimes X)(X^{\otimes 3}\otimes \te{coev}_X)$ and $f^*=\sqrt{2}(X^{\otimes 3}\te{ev}_X)(Q_4 \otimes X)$. Then $ff^*=2(Q_4 \otimes X)e_4(Q_4\otimes X)=Q_4 \otimes X$ and $f^*f=2E_4(Q_4)=2\cdot \f{1}{2}f^{(3)}=P_3$. Thus $Q_4\otimes X\cong P_3$. 

    Similarly, in Lemma \ref{lem: e7sufficientleafrelationp6}, we proved that $P_6 \otimes X=2(P_6 \otimes X)e_6(P_6 \otimes X)$. Let $g=\sqrt{2}(P_6 \otimes X)(X^{\otimes 5}\otimes \te{coev}_X)$ and $g^*=\sqrt{2}(X^{\otimes 5}\otimes \te{ev}_X)(P_6 \otimes X)$. Then $gg^*=2(P_6 \otimes X)e_6(P_6 \otimes X)=P_6 \otimes X$ and $g^*g=2E_6(P_6)=2\cdot \f{1}{2}P_5=P_5$. Therefore, $P_6 \otimes X \cong P_5$, which completes the proof.
\end{proof}

We conclude this section with the proof of the main theorem. 

\begin{proof}[Proof of Theorem \ref{thm: E7sufficientmaintheorem}]
    In Corollary \ref{cor: e7sufficientp0onedim} we show that the space of closed diagrams is one-dimensional. To show $\mc{P}$ is spherical it suffices to show $\mc{P}_1=\te{Hom}(X,X)=\te{span}\{X\}$. However, this is clear as $X=f^{(1)}$ is a minimal projection. 

    To show positive-definiteness, we employ section 4.2 of \cite{MPS10}. In that section, they present an explicit positive orthogonal basis of the box spaces for any unshaded planar algebra whose corresponding category of projections is semisimple, with all minimal projections having positive trace. Then they conclude in Theorem 4.18 that their presentation results in a positive definite planar algebra. Since we already know the principal graph of $\mc{P}$ is affine $E_7$ from Lemma \ref{lem: e7sufficientprincipalgraphissatisfied} and Lemma \ref{lem: E7sufficienttraces} gives that all the minimal projections have positive trace, we obtain that $\mc{P}$ is positive definite. Now that we have an explicit basis, we get that $\mc{P}$ is evaluable.
\end{proof}

\section{Necessary Relations for Affine \texorpdfstring{$E_7$}{E7} Unshaded Subfactor Planar Algebras}\label{necessaryE7}

Let $\mc{P}$ be the (unshaded) subfactor planar algebra with principal graph $\til{E}_7$. In this section, we show that all of the relations given in the presentation of Theorem \ref{thm: E7sufficientmaintheorem} were necessary. Let $\mc{C}_{\mc{P}}$ be the category created from $\mc{P}$. Define $X$ to be a single, black, unoriented strand. Let $\mc{P}$ have a principal graph with the following labeling:

\begin{equation}\label{eq: generalE7graph2}
    \arbEsevengraph
\end{equation}

\begin{lem}\label{lem: e7necessarytracesofminproj}
    $\mc{P}$ has index 4. Further, the traces of all the minimal projections are:
    \begin{enumerate}[label=(\roman*)]
        \item $\te{tr}(\emptyset)=1$, 
        \item $\te{tr}(P_1)=2$,
        \item $\te{tr}(P_2)=3$,
        \item $\te{tr}(P_3)=4$,
        \item $\te{tr}(Q_4)=2$,
        \item $\te{tr}(P_4)=3$,
        \item $\te{tr}(P_5)=2$, and
        \item $\te{tr}(P_6)=1$.
    \end{enumerate}
\end{lem}
    
\begin{proof}
    By assumption $\te{tr}(\emptyset)=1$. Choose $\delta = 2$ and the map $\te{tr}$ from the Perron-Frobenius theorem to work as stated in the lemma. In this case, we can compute the traces of the other minimal projections as follows:
    \begin{align*}
        &2\te{tr}(\emptyset)=\te{tr}(P_1)=2,
        \quad 2\te{tr}(P_1)-\te{tr}(\emptyset)=\te{tr}(P_2)=3, \quad
        2\te{tr}(P_2)-\te{tr}(P_1)=\te{tr}(P_3)=4, \text{ and } \quad
        \te{tr}(Q_4)=\f{1}{2}\te{tr}(P_3)=2. \\
    \end{align*}
    
    Further,
    \begin{align*}
        &2\te{tr}(P_3)-\te{tr}(Q_4)-\te{tr}(P_2)=\te{tr}(P_4)=3, \quad
        2\te{tr}(P_4)-\te{tr}(P_3)=\te{tr}(P_5)=2, \te{ and } \quad 
        \te{tr}(P_6)=\f{1}{2}\te{tr}(P_5)=1.
    \end{align*}
    Since the trace formula holds for all vertices, by Perron-Frobenius these values must be correct. 
\end{proof}

\begin{lem}
    The following isomorphisms hold in $\mc{P}$:
    \begin{enumerate}[label=(\roman*)]
        \item $P_1 \cong f^{(1)}=X$,
        \item $P_2 \cong f^{(2)}$, and
        \item $P_3 \cong f^{(3)}$.
    \end{enumerate}
\end{lem}

\begin{proof}
    By the principal graph and Wenzl's relation, $\emptyset \otimes X \cong P_1 \cong f^{(1)}$. Similarly, $P_1 \otimes X \cong \emptyset \oplus P_2 \cong f^{(0)}\oplus f^{(2)}$, which gives that $f^{(2)}\cong P_2$. Using the same method, $P_2 \otimes X \cong P_1 \oplus P_3 \cong f^{(1)}\cong f^{(3)}$, so $f^{(3)}\cong P_3$.
\end{proof}

\begin{lem}
    $X^{\otimes 4}\cong \emptyset \oplus \emptyset \oplus P_2 \oplus P_2 \oplus P_2 \oplus Q_4 \oplus P_4$. 
\end{lem}

\begin{proof}
    Recall the result given by Jones in Proposition \ref{prop: JonesdirectsumoftensorsofX}. To compute $X^{\otimes 4}$ it is enough to find all endpoints of a length $4$ path starting at $\emptyset$. The lemma is then immediate.
\end{proof}

\begin{lem}\label{lem: E7uniquereps}
 There exist no representatives of $[Q_4],[P_4]\in \mc{P}_m$ where $0\leq m\leq 3$. Further, there exists a unique representative of $[Q_4]$ and $[P_4]$ in $\mc{P}_4$ such that $Q_4+P_4=f^{(4)}$.
\end{lem}

\begin{proof}
    Suppose $Q_4'$ (respectively, $P_4'$) is a representative of $Q_4$ (respectively, $P_4$) in $\mc{P}_m$ where $0\leq m\leq 3$. However, by Proposition \ref{prop: JonesdirectsumoftensorsofX}, $Q_4$ and $P_4$ cannot be in the tensor decomposition of $X^{\otimes m}$ since they are both distance 4 away from $\emptyset$. Thus $\te{dim}(\te{Hom}(X^{\otimes m}, Q_4))=\te{dim}(\te{Hom}(X^{\otimes m}, P_4))=0$, so $Q_4'=P_4'=0$, a contradiction.

   The rest of this proof is nearly identical to Lemma 3.2 of \cite{Mol24} so will be omitted.
\end{proof}

\begin{rem*}
    From now on, define $P_4$ and $Q_4$ to be the representatives found in the previous lemma. In particular, $P_4, Q_4 \in \mc{P}_4$ and $f^{(4)}=P_4+Q_4$. 
\end{rem*}

\begin{defn}\label{defn: e7necessarygenerator} (Generator of the affine $E_7$ Subfactor Planar Algebra)
    $\mc{P}$ with the principal graph in figure \ref{eq: generalE7graph2} has a special element $S\in \mc{P}_4$: 
    \begin{align*}
          \begin{overpic}[unit=1mm, scale=0.4]{S-4strand.png}
        \put(4.5, 6.4){$S$}
    \end{overpic}
    \end{align*}
    where $S\coloneq 3Q_4-2P_4$.
\end{defn}

We will prove at the end of this section that $S$ is indeed a generator of this planar algebra. In the work that follows, we will prove that if $\mc{P}$ is a subfactor planar algebra with affine $E_7$ principal graph, then $\mc{P}$ will have the following relations. (Notice we already proved relation (i) since $\mc{P}$ has index 4.)

\begin{defn}\label{defn: e7necessaryrelations} (Relations of the affine $E_7$ Subfactor Planar Algebra) Relations of the $\til{E}_7$ subfactor planar algebra include:
\begin{enumerate}[label=(\roman*)]
         \item (\textbf{the bubble relation})  $\bubble=2$,
          \item (\textbf{$S$ is U.C.C.S.})
      \quad {\makebox[0pt][l]{\raisebox{-0.6cm}{{$\begin{aligned}[t]  \begin{overpic}[unit=1mm, scale=0.4]{Scap1-4.png}
        \put(2.7, 7.4){$S$}
    \end{overpic} 
    \quad {\makebox[0pt][l]{\raisebox{0.6cm}{{$=...=$}}}}\quad\quad\quad
  \begin{overpic}[unit=1mm, scale=0.4]{Scap2-4.png}
        \put(2.9, 7.5){$S$}
    \end{overpic}
     \quad {\makebox[0pt][l]{\raisebox{0.6cm}{{$=$}}}}\quad
     \begin{overpic}[unit=1mm, scale=0.4]{Scup-4.png}
        \put(3.4, 7.5){$S$}
    \end{overpic} 
    \quad {\makebox[0pt][l]{\raisebox{0.6cm}{{$=...=$}}}}\quad\quad\quad
     \begin{overpic}[unit=1mm, scale=0.4]{Ssidecap-4.png}
        \put(3, 7.5){$S$}
    \end{overpic}  \quad {\makebox[0pt][l]{\raisebox{0.6cm}{{$=$}}}}
    \quad 
     \begin{overpic}[unit=1mm, scale=0.4]{Ssidecap2-4.png}
        \put(8.4, 7.5){$S$}
    \end{overpic} \quad {\makebox[0pt][l]{\raisebox{0.6cm}{{$=0$,}}}}
\end{aligned}$}}}}
    \item (\textbf{$S^2$ relation}) 
          {\makebox[0pt][l]{\raisebox{-0.6cm}{{$\begin{aligned}[t] \begin{overpic}[unit=1mm, scale=0.4]{S-4strand.png}
        \put(4.5, 6.4){$S^2$}
    \end{overpic}
     \quad {\raisebox{0.6cm}{{$=6$}}}
     \begin{overpic}[unit=1mm, scale=0.4]{S-4strand.png}
        \put(3.9, 6.4){$f^{(4)}$}
    \end{overpic}
    \quad  {\raisebox{0.6cm}{{$+$}}}
     \begin{overpic}[unit=1mm, scale=0.4]{S-4strand.png}
        \put(4.5, 6.4){$S$}
    \end{overpic}
    \end{aligned}$}}}}
    \item (\textbf{Jellyfish relation})
     {\makebox[0pt][l]{\raisebox{-0.6cm}{{$\begin{aligned}[t] \begin{overpic}[unit=1mm, scale=0.4]{capoverbox2-4.png}
       \put(8.5, 6.8){$S$}
    \end{overpic}\quad {\raisebox{0.6cm}{{$=$}}}
      \begin{overpic}[unit=1mm, scale=0.4]{boxwithbelowcrossings-4.png}
        \put(8.5, 12.9){$S$}
    \end{overpic}
    \end{aligned}$}}}}
    \end{enumerate}
\end{defn}

The rest of the section is devoted to proving the following theorem.

\begin{thm}\label{thm: e7necessarymaintheorem}
    If $\mc{P}$ is an unshaded subfactor planar algebra with principal graph $\til{E}_7$ with labeling given in figure \ref{eq: generalE7graph2}, then $\mc{P}$ is generated by the element $S$ given in Definition \ref{defn: e7necessarygenerator} with relations given in Definition \ref{defn: e7necessaryrelations}. 
\end{thm}

\begin{lem}
    We have the leaf relations $Q_4 \otimes X = 2(Q_4\otimes X)e_4(Q_4 \otimes X)$ and $P_6 \otimes X = 2(P_6 \otimes X)e_6(P_6 \otimes X)$. 
\end{lem}

\begin{proof}
    This follows from $Q_4$ and $P_6$ being leaves in the principal graph and Lemma \ref{lem: leafproperty}. 
\end{proof}

\begin{lem}
    $Q_4$ and $P_4$ are uncuppable and uncappable.
\end{lem}

\begin{proof}
    Since $f^{(4)}=Q_4+P_4$, we obtain that $f^{(4)}Q_4=Q_4^2+P_4Q_4$ and $f^{(4)}P_4=Q_4P_4+P_4^2$. Using that $P_4$ and $Q_4$ are nonisomorphic minimal projections we obtain that $f^{(4)}Q_4=Q_4$ and $f^{(4)}P_4=P_4$. A symmetric proof proves that $Q_4$ and $P_4$ are also uncuppable. 
\end{proof}

\begin{lem}\label{lem: E7necessaryminprojs}
    We can choose the following representatives for the minimal projections:
    \begin{enumerate}[label=(\roman*)]
        \item $Q_4\coloneq \f{2}{5}f^{(4)}+\f{1}{5}S$,
        \item $P_4\coloneq \f{3}{5}f^{(4)}-\f{1}{5}S$,
        \item $P_5\coloneq P_4 \otimes X -\f{4}{3} (P_4 \otimes X)e_4(P_4\otimes X)$, and
        \item $P_6 \coloneq P_5 \otimes X - \f{3}{2}(P_5 \otimes X)e_5 (P_5 \otimes X)$.
    \end{enumerate}
\end{lem}

\begin{proof}
    Since $Q_4$ and $P_4$ defined in the lemma are in $\mc{P}_4$ it is enough for parts (i) and (ii) to show that $Q_4$ and $P_4$ are identical to the representatives chosen from Lemma \ref{lem: E7uniquereps}. Using that $f^{(4)}=Q_4+P_4$ and $S=3Q_4-2P_4$ we obtain
    \begin{align*}
         \f{2}{5}f^{(4)}+\f{1}{5}S&= \f{2}{5}(Q_4+P_4)+\f{1}{5}(3Q_4-2P_4)=Q_4, \te{ and }\\
         \f{3}{5}f^{(4)}-\f{1}{5}S &= \f{3}{5}(P_4+Q_4)-\f{1}{5}(3Q_2-2P_4)=P_4.
    \end{align*}

    By the principal graph, $P_4\otimes X \cong P_3 \oplus P_5$. Therefore, there exists an $a\in \mbb{C}$ such that $P_5=P_4\otimes X-a(P_4 \otimes X)(X^{\otimes 3}\otimes \te{coev}_X)P_3(X^{\otimes 3}\otimes \te{ev}_X)(P_4\otimes X)$. Since $P_3=f^{(3)}$ and $P_4$ is uncuppable and uncappable, $P_5=P_4 \otimes X - a(P_4 \otimes X)e_4(P_4 \otimes X)$. Taking the trace gives that $\te{tr}(P_5)=2=2\te{tr}(P_4)-a\te{tr}(P_4)=(2-a)\cdot3$ so $a=\f{4}{3}$.

    Before proving part (iv), notice that $P_5(P_4 \otimes X)=(P_4 \otimes X)^2-\f{4}{3}(P_4 \otimes X)e_4(P_4 \otimes X)^2=P_5$ since $P_4^2=P_4$. Then, using the principal graph, there exists a $b\in \mbb{C}$ such that $P_6=P_5 \otimes X - b (P_5 \otimes X)(X^{\otimes 4}\otimes \te{coev}_X)P_4(X^{\otimes 4}\otimes \te{ev}_X)(P_5 \otimes X)=P_5 \otimes X - b (P_5 \otimes X)e_5(P_5 \otimes X)$. Taking the trace gives $b=\f{3}{2}$.
\end{proof}

\begin{lem}\label{lem: e7necessarySsquared}
    In $\mc{P}$, $S^2=6f^{(4)}+S$. 
\end{lem}

\begin{proof}
    Since $Q_4$ and $P_4$ are nonisomorphic minimal projections, $Q_4^2=Q_4$, $P_4^2=P_4$, and $Q_4P_4=P_4Q_4=0$. Thus 
    \begin{align*}
        S^2=(3Q_4-2P_4)^2=9Q_4+4P_4=6(Q_4+P_4)+(3Q_4-2P_4)=6f^{(4)}+S
    \end{align*}
    as we claimed.
\end{proof}

\begin{lem}
    $S$ is self-adjoint.
\end{lem}

\begin{proof}
    $S$ is a linear combination of projections, so the lemma follows.
\end{proof}

\begin{lem}\label{lem: e7necessarycapoverSf10is0}
    $\mc{O}(\mc{R}(S))f^{(10)}=0$. 
\end{lem}

\begin{proof}
    If the dimension of $\te{Hom}(f^{(10)}, \emptyset)=0$ then the lemma is true. We then focus on proving this. Recall that $f^{(1)}\cong P_1$, $f^{(2)}\cong P_2$, $f^{(3)}\cong P_3$, and $f^{(4)}\cong P_4 \oplus Q_4$. Wenzl's relation and the principal graph gives 
    \begin{align*}
        f^{(4)}\otimes X \cong (P_4 \oplus Q_4)\otimes X \cong P_3 \oplus P_5 \oplus P_3 \cong f^{(3)}\oplus f^{(5)}
    \end{align*}
    so $f^{(5)}\cong P_3 \oplus P_5$. Similarly,
    \begin{align*}
        f^{(5)}\otimes X \cong (P_3\oplus P_5)\otimes X \cong P_2 \oplus Q_4 \oplus P_4 \oplus P_4 \oplus P_6 \cong f^{(4)}\oplus f^{(6)}\cong P_4 \oplus Q_4 \oplus f^{(6)},
    \end{align*}
    which gives that $f^{(6)}\cong P_2 \oplus P_4 \oplus P_6$. Continuing this process, we get
    \begin{align*}
        f^{(6)}\otimes X \cong (P_2 \oplus P_4 \oplus P_6)\otimes X 
        \cong P_1 \oplus P_3 \oplus P_3 \oplus P_5 \oplus P_5 \cong f^{(5)}\oplus f^{(7)}\cong P_3\oplus P_5 \oplus f^{(7)}
    \end{align*}
    so $f^{(7)}\cong P_1 \oplus P_3 \oplus P_5$. Next, we find that
    \begin{align*}
        f^{(7)}\otimes X \cong (P_1 \oplus P_3 \oplus P_5) \otimes X\cong \emptyset \oplus P_2 \oplus P_2 \oplus Q_4 \oplus P_4 \oplus P_4 \oplus P_6 \cong f^{(6)}\oplus f^{(8)}\cong P_2 \oplus P_4 \oplus P_6 \oplus f^{(8)}
    \end{align*}
    giving that $f^{(8)}\cong \emptyset \oplus P_2 \oplus Q_4 \oplus P_4$. Similarly, we have 
    \begin{align*}
        f^{(8)}\otimes X \cong (\emptyset \oplus P_2 \oplus Q_4 \oplus P_4)\otimes X \cong P_1 \oplus P_1 \oplus P_3 \oplus P_3 \oplus P_3 \oplus P_5 \cong f^{(7)}\oplus f^{(9)} \cong P_1 \oplus P_3 \oplus P_5 \oplus f^{(9)}
    \end{align*}
    which gives that $f^{(9)}\cong P_1 \oplus P_3 \oplus P_3$. 
    
    Finally,
    \begin{align*}
        f^{(9)}\otimes X \cong (P_1 \oplus P_3 \oplus P_3) \otimes X \cong \emptyset \oplus P_2 \oplus P_2 \oplus P_4 \oplus Q_4 \oplus P_2 \oplus Q_4 \oplus P_4 \cong f^{(8)}\oplus f^{(10)} \cong \emptyset \oplus P_2 \oplus Q_4 \oplus P_4 \oplus f^{(10)}
    \end{align*}
    Therefore $f^{(10)}\cong P_2 \oplus P_2 \oplus Q_4 \oplus P_4$, which gives the result.
\end{proof}

\begin{lem}
    In $\mc{P}$ we have the following partial traces.
    \begin{enumerate}[label=(\roman*)]
        \item $E_4(P_4)=E_4^\ell(P_4)=\f{3}{4}P_3$, and
        \item $E_4(Q_4)=E_4^\ell(Q_4)=\f{1}{2}P_3$.
    \end{enumerate}
\end{lem}

\begin{proof}
    Since $P_4$ is uncuppable and uncappable, so is $E_4(P_4)$. Thus we recognize that $E_4(P_4)\in \te{Hom}(f^{(3)}, f^{(3)})$ since $E_4(P_4)=f^{(3)}E_4(P_4)f^{(3)}$. Since $f^{(3)}$ is a minimal projection, there exists an $a \in \mbb{C}$ such that $E_4(P_4)=af^{(3)}$. Taking the trace of both sides gives that $a=\f{3}{4}$. Similarly, there exists a $b \in \mbb{C}$ such that $E_4^\ell(P_4)=f^{(3)}E_4^\ell(P_4)f^{(3)}=bf^{(3)}$. Take the trace of both $E_4^\ell(P_4)$ and $f^{(3)}$. Use sphericality to notice that $\te{tr}(E_4^\ell(P_4))=\te{tr}(P_4)=3$. Therefore, $b=\f{3}{4}$. 

    We do nearly the same argument for $Q_4$ to get part (ii). 
\end{proof}

\begin{lem}\label{lem: e7necessarySisuccs}
    $S$ is U.C.C.S.
\end{lem}

\begin{proof}
    $S$ is a linear combination of uncuppable and uncappable elements so $S$ itself is uncuppable and uncappable. Further,
    \begin{align*}
        E_4(S)=3E_4(Q_4)-2E_4(P_4)=3\cdot \f{1}{2}P_3-2\cdot \f{3}{4}P_3=0
    \end{align*}
    and similarly, $E_4^\ell(S)=0$. 
\end{proof}

\begin{lem}
    $S$ is nonzero.
\end{lem}

\begin{proof}
    Since $S$ is U.C.C.S., $\te{tr}(S)=0$. So, $\te{tr}(S^2)=6\te{tr}(f^{(4)})+\te{tr}(S)=6\cdot 5=30$. Thus $S^2$ is nonzero and so is $S$. 
\end{proof}

\begin{lem}
    There exists an $a\in \mbb{C}^\times$ such that $\mc{F}(S)=aS$. 
\end{lem}

\begin{proof}
  Since $S$ is U.C.C.S., $\mc{R}(S)f^{(8)}=\mc{R}(S)\neq 0$ and $\mc{R}(\mc{F}(S))f^{(8)}=\mc{R}(\mc{F}(S))$. Recall in Lemma \ref{lem: e7necessarycapoverSf10is0} that $f^{(8)}=\emptyset \oplus P_2 \oplus P_4 \oplus Q_4$, so the dimension of $\te{Hom}(f^{(8)}, \emptyset)$ is $1$. Therefore, there exists a scalar $a\in \mbb{C}$ such that $\mc{R}(S)f^{(8)}=\mc{R}(S)=a\mc{R}(\mc{F}(S))f^{(8)}=a\mc{R}(\mc{F}(S))$. Since $S$ is nonzero, so is $a$. 
\end{proof}

\begin{lem}
    In $\mc{P}$, $\mc{F}(S)=S$. 
\end{lem}

\begin{proof}
    Recall that  $Q_4 \otimes X = 2(Q_4 \otimes X)e_4(Q_4 \otimes X)$. We first take the left partial trace of the leaf property: $E_5^\ell(Q_4 \otimes X)=2E_5^\ell((Q_4\otimes X)e_4(Q_4\otimes X))$. Then we expand out using $Q_4=\f{2}{5}f^{(4)}+\f{1}{5}S$ to obtain on the left-hand side:
    \begin{align*}
        E_5^\ell(Q_4 \otimes X) = \f{2}{5}E_5^\ell(f^{(4)}\otimes X)+\f{1}{5}E_5^\ell(S \otimes X)=\f{2}{5}E_5^\ell(f^{(4)}\otimes X)=\f{2}{5}(E_4^\ell(f^{(4)})\otimes X)
    \end{align*}
    since $S$ is U.C.C.S. We also obtain $2E_5^\ell((Q_4\otimes X)e_4(Q_4\otimes X))$ equals
    \begin{align*}
        \f{8}{25}E_5^\ell((f^{(4)}\otimes X)e_4(f^{(4)}\otimes X))+\f{4}{25}E_5^\ell((f^{(4)}\otimes X)e_4(S\otimes X))+\f{4}{25}E_5^\ell((S\otimes X)e_4(f^{(4)}\otimes X))+\f{2}{25}E_5^\ell((S\otimes X)e_4(S\otimes X)).
    \end{align*}
    Now we rainbow the expanded versions of $E_5^\ell(Q_4\otimes X)$ and $2E_5^\ell((Q_4\otimes X)e_4(Q_4\otimes X))$ onto $f^{(8)}$. We first obtain
    \begin{align*}
        \mc{R}(E_5^\ell(Q_4 \otimes X))f^{(8)}=\f{2}{5}(E_4^\ell(f^{(4)})\otimes X)f^{(8)}=0
    \end{align*}
    since $f^{(8)}$ is uncappable. We compute $2\mc{R}(E_5^\ell((Q_4\otimes X)e_4(Q_4\otimes X)))f^{(8)}$ term-by-term. First, we notice that $\f{8}{25}\mc{R}(E_5^\ell((f^{(4)}\otimes X)e_4(f^{(4)}\otimes X)))f^{(8)}$ is zero since on top of $f^{(8)}$ involves only Temperley-Lieb.
    
    The second term, $\f{4}{25}\mc{R}(E_5^\ell((f^{(4)}\otimes X)e_4(S\otimes X)))f^{(8)}$, expands to 
    \begin{align*}
        &\f{4}{25}(\mc{R}(E_5^\ell((f^{(3)}\otimes X^{\otimes 2})e_4(S\otimes X)))f^{(8)}-\f{3}{4}\mc{R}(E_5^\ell((f^{(3)}\otimes X^{\otimes 2})e_3(f^{(3)}\otimes X^{\otimes 2})e_4(S\otimes X)))f^{(8)})\\
        &=-\f{3}{4}\mc{R}(E_5^\ell((f^{(3)}\otimes X^{\otimes 2})e_3(f^{(3)}\otimes X^{\otimes 2})e_4(S\otimes X)))f^{(8)}\\
        &=-\f{3}{25}\mc{R}(E_5^\ell((f^{(3)}\otimes X^{\otimes 2})e_3e_4(S\otimes X)))f^{(8)}
    \end{align*}
    where the first equality comes from $f^{(8)}$ being uncuppable and uncappable and the second equality comes from $S$ being uncuppable and uncappable. Expanding $f^{(3)}$ in $-\f{3}{25}\mc{R}(E_5^\ell((f^{(3)}\otimes X^{\otimes 2})e_3e_4(S\otimes X)))f^{(8)}$ gives
    \begin{align*}
        &-\f{3}{25}(\mc{R}(E_5^\ell((f^{(2)}\otimes X^{\otimes 3})e_3e_4(S\otimes X)))f^{(8)}-\f{2}{3}\mc{R}(E_5^\ell((f^{(2)}\otimes X^{\otimes 3})e_2(f^{(2)}\otimes X^{\otimes 3})e_3e_4(S\otimes X)))f^{(8)})\\
        &=\f{2}{25}\mc{R}(E_5^\ell((f^{(2)}\otimes X^{\otimes 3})e_2(f^{(2)}\otimes X^{\otimes 3})e_3e_4(S\otimes X)))f^{(8)}\\
        &=\f{2}{25}\mc{R}(E_5^\ell((f^{(2)}\otimes X^{\otimes 3})e_2e_3e_4(S\otimes X)))f^{(8)}
    \end{align*}
    where the first equality comes from $f^{(8)}$ being uncuppable and uncappable and the second equality comes from $S$ being uncuppable and uncappable. Expanding $f^{(2)}$ in $ \f{2}{25}\mc{R}(E_5^\ell((f^{(2)}\otimes X^{\otimes 3})e_2e_3e_4(S\otimes X)))f^{(8)}$ gives
    \begin{align*}
       &\f{2}{25}(\mc{R}(E_5^\ell((X^{\otimes 5})e_2e_3e_4(S\otimes X)))f^{(8)}
        -\f{1}{2}\mc{R}(E_5^\ell(e_1e_2e_3e_4(S\otimes X)))f^{(8)})\\
        &=-\f{1}{25}\mc{R}(E_5^\ell(e_1e_2e_3e_4(S\otimes X)))f^{(8)}),
        \end{align*}
        where the equality comes from $f^{(8)}$ being uncuppable and uncappable. This equals
        
        \begin{equation*}
        -\f{1}{25}\mc{R}(\mc{F}(S))f^{(8)}=-\f{a}{25}\mc{R}(S)f^{(8)}=-\f{a}{25}\mc{R}(S),
    \end{equation*}
     since $S$ is U.C.C.S. 

    The third term, $\f{4}{25}\mc{R}(E_5^\ell((S\otimes X)e_4(f^{(4)}\otimes X)))f^{(8)}$, becomes
    \begin{align*}
        &\f{4}{25}(\mc{R}(E_5^\ell((S\otimes X)e_4(f^{(3)}\otimes X^{\otimes 2})))f^{(8)}-\f{3}{4}\mc{R}(E_5^\ell((S\otimes X)e_4(f^{(3)}\otimes X^{\otimes 2})e_3(f^{(3)}\otimes X^{\otimes 2})))f^{(8)})\\
        &= -\f{3}{25}\mc{R}(E_5^\ell((S\otimes X)e_4(f^{(3)}\otimes X^{\otimes 2})e_3(f^{(3)}\otimes X^{\otimes 2})))f^{(8)}\\
        &=-\f{3}{25}\mc{R}(E_5^\ell((S\otimes X)e_4e_3(f^{(3)}\otimes X^{\otimes 2})))f^{(8)}\\
    \end{align*}
    where the first equality comes from $f^{(8)}$ being uncuppable and uncappable and the second equality comes from $S$ being uncuppable and uncappable. Then $ -\f{3}{25}\mc{R}(E_5^\ell((S\otimes X)e_4e_3(f^{(3)}\otimes X^{\otimes 2})))f^{(8)}$ equals
    \begin{align*}
       &-\f{3}{25}(\mc{R}(E_5^\ell((S\otimes X)e_4e_3(f^{(2)}\otimes X^{\otimes 3})))f^{(8)}-\f{2}{3}\mc{R}(E_5^\ell((S\otimes X)e_4e_3(f^{(2)}\otimes X^{\otimes 3})e_2(f^{(2)}\otimes X^{\otimes 3})))f^{(8)})\\
        &=\f{2}{25}\mc{R}(E_5^\ell((S\otimes X)e_4e_3(f^{(2)}\otimes X^{\otimes 3})e_2(f^{(2)}\otimes X^{\otimes 3})))f^{(8)}\\
        &=\f{2}{25}\mc{R}(E_5^\ell((S\otimes X)e_4e_3e_2(f^{(2)}\otimes X^{\otimes 3})))f^{(8)}
    \end{align*}
    where the second equality comes from $f^{(8)}$ being uncuppable and uncappable and the third equality comes from $S$ being uncuppable and uncappable. Finally, expanding $f^{(2)}$ gives
    \begin{align*}
        \f{2}{25}\mc{R}(E_5^\ell((S\otimes X)e_4e_3e_2(f^{(2)}\otimes X^{\otimes 3})))f^{(8)}&=\f{2}{25}(\mc{R}(E_5^\ell((S\otimes X)e_4e_3e_2(X^{\otimes 5})))f^{(8)}-\mc{R}(E_5^\ell((S\otimes X)e_4e_3e_2e_1))f^{(8)})\\
        &=-\f{1}{25}\mc{R}(E_5^\ell((S\otimes X)e_4e_3e_2e_1))f^{(8)})\\
        &=-\f{1}{25}\mc{R}(\mc{F}^{-1}(S))f^{(8)}\\
        &=-\f{a^{-1}}{25}\mc{R}(S)f^{(8)}\\
        &=-\f{a^{-1}}{25}\mc{R}(S),
    \end{align*}
    where the second equality comes from $f^{(8)}$ being uncuppable and uncappable and the fifth equality comes from $S$ being U.C.C.S. 

    The fourth and final term expands out in the following way
    \begin{align*}
        \f{2}{25}\mc{R}(E_5^\ell((S\otimes X)e_4(S\otimes X)))f^{(8)}&=\f{2}{25}\mc{R}(S^2)f^{(8)}
        =\f{2}{25}(\mc{R}(S)f^{(8)}+\mc{R}(f^{(4)})f^{(8)})
        =\f{2}{25}\mc{R}(S)f^{(8)}
        =\f{2}{25}\mc{R}(S),
    \end{align*}
    where the first equality comes from clicking the top $S$-box counterclockwise and the bottom $S$-box clockwise, the third equality comes from $f^{(8)}$ being uncuppable and uncappable, and the fourth equality comes from $S$ being U.C.C.S. 

    Therefore, $\mc{R}(E_5^\ell(Q_4\otimes X)=2\mc{R}(E_5^\ell((Q_4\otimes X)e_4(Q_4\otimes X)))f^{(8)}$ is equivalent to $0=(-\f{a}{25}-\f{a^{-1}}{25}+\f{2}{25})\mc{R}(S)$. Since $\mc{R}(S)\neq 0$ this gives $a=1$.
\end{proof}

\begin{lem}\label{lem: e7necessaryjellyfishrelation}
    $\mc{P}$ satisfies the jellyfish relation of Definition \ref{defn: e7necessaryrelations}. 
\end{lem}

\begin{proof}
    Define the below diagram to be $T$.
    \begin{align*}
          \begin{overpic}[unit=1mm, scale=0.4]{boxwithbelowcrossings-4.png}
        \put(8.5, 12.9){$S$}
    \end{overpic}
    \end{align*}
    Recall in Lemma \ref{lem: e7necessarycapoverSf10is0} we showed $\mc{O}(\mc{R}(S))f^{(10)}=0$. The same proof gives that $Tf^{(10)}=0$. We can expand $f^{(10)}$ in both $\mc{O}(\mc{R}(S))f^{(10)}$ and $Tf^{(10)}$. The same terms in $f^{(10)}$ will cap $S$ and become zero. The nonzero terms will be the identity on 10 strands and any element with a cup on the first two or last two strands. 

    Suppose $T$ is multiplied by an element in the expansion of $f^{(10)}$ that has a cup on the first two strands. Recall, from Lemma \ref{lem: cuponcrossingresolved} and Lemma \ref{lem: boxwithovercrossingonbottom}, that 
    \begin{align*}
        \begin{overpic}[unit=1mm, scale=0.5]{boxwithovercrossingonbottomtwist.png}
        \put(7,-2){\tiny{$(7)$}}
        \put(7,2){\tiny{$...$}}
        \put(6,17.5){$S$} 
        \end{overpic}
        {\raisebox{1cm}{$=-i$}}
        \begin{overpic}[unit=1mm, scale=0.5]{boxwithovercrossingonbottom2.png}
        \put(7,-2){\tiny{$(7)$}}
        \put(7,2){\tiny{$...$}}
        \put(6,17.5){$S$}
    \end{overpic}
    {\raisebox{1cm}{$=\mc{R}(S)$.}}
    \end{align*}
    Meanwhile, a cup on the first two strands of $\mc{O}(\mc{R}(S))=\mc{F}(\mc{R}(S))=\mc{R}(S)$. Similarly, we can see that a cup on the last two strands for $T$ and $\mc{O}(\mc{R}(S))$ both equal $\mc{R}(S)$. Therefore, when expanding $f^{(10)}$ for both $\mc{O}(\mc{R}(S))f^{(10)}$ and $Tf^{(10)}$ we get that the terms with cups between the first two or last two strands are equal. Therefore, all the terms in the expansion of $f^{(10)}$ have been shown to be equal except the identity on ten strands. Since $\mc{O}(\mc{R}(S))f^{(10)}=Tf^{(10)}=0$, this gives that the $\mc{O}(\mc{R}(S))X^{\otimes 10}=TX^{\otimes 10}$. Thus $\mc{O}(\mc{R}(S))=T$, as we wished. 
\end{proof}

\begin{proof}[Proof of Theorem \ref{thm: e7necessarymaintheorem}]
    In Lemma \ref{lem: e7necessarytracesofminproj} we showed the bubble relation. In Lemma \ref{lem: e7necessarySisuccs} we showed that $S$ is U.C.C.S. In Lemma \ref{lem: e7necessarySsquared} we showed that $S^2=6f^{(4)}+S$. In Lemma \ref{lem: e7necessaryjellyfishrelation} we showed the jellyfish relation.

    What is left to show is that $S$ generates $\mc{P}$. We know that $S$ generates $\mc{P}'$, a subfactor planar algebra with principal graph $\til{E}_7$, as defined in Theorem \ref{thm: E7sufficientmaintheorem}. Therefore, it is enough to show that the dimension of all the box-spaces are equal. However, this is clear as a positive definite basis of each box-space  created in \cite{MPS10} relies only on the principal graph and the traces of the minimal projections. Since the principal graphs and the trace of minimal projections are equal for $\mc{P}$ and $\mc{P}'$, we get that $S$ must generate $\mc{P}$. 
\end{proof}

\begin{cor}\label{onlyonee7}
    There is exactly one unshaded subfactor planar algebra with principal graph $\til{E}_7$. 
\end{cor}

\begin{proof}
    We have at least one by Theorem \ref{thm: E7sufficientmaintheorem}. Further, if we have a subfactor planar algebra with principal graph $\til{E}_7$, Theorem \ref{thm: e7necessarymaintheorem} showed it must have the same presentation as the one found in Theorem \ref{thm: E7sufficientmaintheorem}.
\end{proof}

\bibliographystyle{alpha}
\bibliography{Cite}

\end{document}